\documentclass[11pt,a4paper]{article}
\usepackage[utf8]{inputenc}

\usepackage[hyphens]{url}
\usepackage{xurl}
\usepackage{breakurl}

\usepackage{mathrsfs,amsthm,xcolor,verbatim,bbm,amsmath,amsfonts,amssymb,nicefrac,enumitem,bm,mathtools,xparse,etoolbox,textcomp}
\usepackage[margin=1in]{geometry}
\usepackage[capitalise]{cleveref} 

 \usepackage{hyperref}

\crefname{enumi}{item}{items}
\crefname{equation}{}{}
\crefname{subsection}{Subsection}{Subsections}
\crefname{cor}{Corollary}{Corollaries}
\crefname{theorem}{Theorem}{Theorems}
\crefname{prop}{Proposition}{Propositions}

\numberwithin{equation}{section}

\hypersetup{
    colorlinks,
    linkcolor={red!80!black},
    citecolor={green},
    urlcolor={blue!80!black}
}

\theoremstyle{plain}
\newtheorem{theorem}{Theorem}[section]
\newtheorem{lemma}[theorem]{Lemma}
\newtheorem{prop}[theorem]{Proposition}
\newtheorem{cor}[theorem]{Corollary}

\newtheorem{setting}[theorem]{Setting}

\theoremstyle{remark}

\theoremstyle{definition}
\newtheorem{definition}[theorem]{Definition}

\DeclareFontEncoding{LS1}{}{}
\DeclareFontSubstitution{LS1}{stix}{m}{n}
\DeclareMathAlphabet{\mathscr}{LS1}{stixscr}{m}{n}

\renewcommand{\P}{\mathbb{P}}

\newcommand{\R}{\mathbb{R}}
\newcommand{\N}{\mathbb{N}}

\newcommand{\bbA}{\mathbb{A}}
\newcommand{\bbB}{\mathbb{B}}

\newcommand{\dens}{\mathfrak{p}}
\newcommand{\mdim}{k}
\newcommand{\tbound}{\fD}

\newcommand{\projm}{\scrp}
\newcommand{\tubn}{\bfP}
\newcommand{\dist}{\mathscr{d}}

\newcommand{\w}[1]{\mathfrak{w}^{#1}}
\renewcommand{\b}[1]{\mathfrak{b}^{#1}}
\renewcommand{\v}[1]{\mathfrak{v}^{#1}}
\newcommand{\q}[1]{\mathfrak{q}^{#1}}
\renewcommand{\c}[1]{\mathfrak{c}^{#1}}

\newcommand{\smallsum}{\textstyle\sum}
\newcommand{\tint}{\textstyle\int}

\newcommand{\with}{\curvearrowleft}

\newcommand{\indicator}[1]{\mathbbm{1}_{\smash{#1}}}

\newcommand{\rank}{\operatorname{rank}}
\newcommand{\Hs}{\operatorname{Hess}}

\newcommand{\realization}[1] {\mathscr{N} ^{ #1  }}
\newcommand{\width}{H}

\newcommand{\Rect}{\mathfrak{R}}

\renewcommand{\d}{ \mathrm{d}}

\newcommand{\cB}{\mathcal{B}}

\newcommand{\cF}{\mathcal{F}}
\newcommand{\cG}{\mathcal{G}}

\newcommand{\cL}{\mathcal{L}}
\newcommand{\cM}{\mathcal{M}}

\newcommand{\cQ}{\mathcal{Q}}

\newcommand{\cT}{\mathcal{T}}

\newcommand{\cV}{\mathcal{V}}

\newcommand{\bfc}{\mathbf{c}}

\newcommand{\bfk}{\mathbf{k}}

\newcommand{\bfv}{\mathbf{v}}

\newcommand{\bfP}{\mathbf{P}}

\newcommand{\scrL}{\mathscr{L}}

\newcommand{\scrN}{\mathscr{N}}

\newcommand{\scrP}{\mathscr{P}}

\newcommand{\scrp}{\mathscr{p}}

\newcommand{\fC}{\mathfrak{C}}
\newcommand{\fD}{\mathfrak{D}}

\newcommand{\fL}{\mathfrak{L}}

\newcommand{\fV}{\mathfrak{V}}

\newcommand{\fb}{\mathfrak{b}}
\newcommand{\fc}{\mathfrak{c}}
\newcommand{\fd}{\mathfrak{d}}

\newcommand{\fg}{\mathfrak{g}}

\newcommand{\fm}{\mathfrak{m}}

\newcommand{\fq}{\mathfrak{q}}

\newcommand{\fs}{\mathfrak{s}}

\newcommand{\fu}{\mathfrak{u}}
\newcommand{\fv}{\mathfrak{v}}
\newcommand{\fw}{\mathfrak{w}}
\newcommand{\fx}{\mathscr{x}}
\newcommand{\fy}{\mathscr{y}}

\renewcommand{\emptyset}{\varnothing}

\DeclarePairedDelimiter{\norm}{\lVert}{\rVert}
\DeclarePairedDelimiter{\abs}{\lvert}{\rvert}
\DeclarePairedDelimiter{\rbr}{(}{)}
\DeclarePairedDelimiter{\br}{[}{]}
\DeclarePairedDelimiter{\cu}{\{}{\}}
\DeclarePairedDelimiter{\spro}{\langle}{\rangle}

\newcommand{\qandq}{\qquad\text{and}\qquad}

\ExplSyntaxOn

\bool_new:N \g_noteobserve

\NewDocumentCommand{\nobs}{}{
  \bool_if:nTF { \g_noteobserve } {
    \bool_gset_false:N \g_noteobserve 
    note~
  } {
    \bool_gset_true:N \g_noteobserve 
    observe~
  }
}

\NewDocumentCommand{\Nobs}{}{
  \bool_if:nTF { \g_noteobserve } {
    \bool_gset_false:N \g_noteobserve 
    Note~
  } {
    \bool_gset_true:N \g_noteobserve 
    Observe~
  }
}

\seq_new:N \g_cflist_loaded
\seq_new:N \g_cflist_pending

\NewDocumentCommand{\cfadd}{ m }
{
  \seq_if_in:NnF \g_cflist_loaded { #1 } {
    \seq_if_in:NnF \g_cflist_pending { #1 } {
      \seq_gput_right:Nn \g_cflist_pending { #1 }
    }
  }
}

\NewDocumentCommand{\cfconsiderloaded}{ m }{
  \seq_gput_right:Nn \g_cflist_loaded {#1}
}

\NewDocumentCommand{\cfremove}{ m }
{
  \seq_gremove_all:Nn \g_cflist_pending { #1 }
}

\NewDocumentCommand{\cfload}{ o }
{
  \seq_if_empty:NTF \g_cflist_pending {\unskip} {
    (cf.\ \cref{\seq_use:Nn \g_cflist_pending {,}})\IfValueTF{#1}{#1~}{\unskip}
    \seq_gconcat:NNN \g_cflist_loaded \g_cflist_loaded \g_cflist_pending
    \seq_gclear:N \g_cflist_pending
  }
}

\NewDocumentCommand{\cfclear} {} {
  \seq_gclear:N \g_cflist_loaded
  \seq_gclear:N \g_cflist_pending
}

\NewDocumentCommand{\cfout}{ o }
{
  \seq_if_empty:NTF \g_cflist_pending {\unskip} {
    (cf.\ \cref{\seq_use:Nn \g_cflist_pending {,}})\IfValueTF{#1}{#1~}{\unskip}
    \seq_gclear:N \g_cflist_pending
  }
}

\NewDocumentCommand{\ifnocf} { m } {
  \seq_if_empty:NT \g_cflist_pending { #1 }
}

\ExplSyntaxOff

\NewDocumentEnvironment{cproof}{m}
{\begin{proof}[Proof of \cref{#1}]}%
{\noindent The proof of~\cref{#1} is thus complete.
\end{proof}}

\NewDocumentEnvironment{cproof2}{m}
{\begin{proof}[Proof of \cref{#1}]}%
{\noindent This completes the proof of~\cref{#1}.
\end{proof}}

\title{A proof of convergence for the gradient descent\\ optimization method 
with random initializations\\
in the training of neural networks with ReLU\\
activation for piecewise linear target functions 
}

\author{Arnulf Jentzen$^{1, 2}$ and
Adrian Riekert$^3$
\bigskip
\\
\small{$^1$ Applied Mathematics: Institute for Analysis and Numerics,}
\vspace{-0.1cm}\\
\small{University of Münster, Germany, e-mail: \texttt{ajentzen}\textcircled{\texttt{a}}\texttt{uni-muenster.de}}
\smallskip
\\
\small{$^2$ School of Data Science and Shenzhen Research Institute of Big Data,}
\vspace{-0.1cm}\\
\small{The Chinese University of Hong Kong, Shenzhen, China, e-mail: \texttt{ajentzen}\textcircled{\texttt{a}}\texttt{cuhk.edu.cn}}
\smallskip
\\
\small{$^3$ Applied Mathematics: Institute for Analysis and Numerics,}
\vspace{-0.1cm}\\
\small{University of Münster, Germany, e-mail: \texttt{ariekert}\textcircled{\texttt{a}}\texttt{uni-muenster.de}}}

\date{\today}

\begin{document}

\maketitle

\begin{abstract}
\noindent
    Gradient descent (GD) type optimization methods are the standard instrument to train artificial neural networks (ANNs) with rectified linear unit (ReLU) activation.
    Despite the great success of GD type optimization methods in numerical simulations for the training of ANNs with ReLU activation, it remains -- even in the simplest situation of the plain vanilla GD optimization method with random initializations and ANNs with one hidden layer --
an open problem to prove (or disprove) the conjecture that the risk of the GD optimization method converges in the training of such ANNs to zero as the width of the ANNs, the number of independent random initializations, and the number of GD steps increase to infinity.
In this article we prove this conjecture in the situation where the probability distribution of the input data is equivalent to the continuous uniform distribution on a compact interval, where the probability distributions for the random initializations of the ANN parameters are standard normal distributions, and where the target function under consideration is continuous and piecewise affine linear. 
Roughly speaking, the key ingredients in our mathematical convergence analysis are 
(i) to prove that suitable sets of global minima of the risk functions are \emph{twice continuously differentiable submanifolds of the ANN parameter spaces},
(ii) to prove that the Hessians of the risk functions on these sets of global minima satisfy an appropriate \emph{maximal rank condition}, and, thereafter,
(iii) to apply the machinery in [Fehrman, B., Gess, B., Jentzen, A., Convergence rates for the stochastic gradient descent method for non-convex objective functions. J. Mach. Learn. Res. 21(136): 1--48, 2020] to establish convergence of the GD optimization method with random initializations.

\end{abstract}

\tableofcontents

\section{Introduction}

Gradient descent (GD) type optimization methods are the standard schemes to train artificial neural networks (ANNs) with rectified linear unit (ReLU) activation; cf., e.g., Goodfellow et al.~\cite[Chapter 5]{GoodfellowBengioCourville2016}. Even though GD type optimization methods seem to perform very effectively in numerical simulations, until today in general there is no mathematical convergence analysis in the literature which explains the success of GD optimization methods in the training of ANNs with ReLU activation. 

There are, however, several promising mathematical analysis approaches for GD optimization methods in the scientific literature.
In the case of convex objective functions,
the convergence of GD type optimizations methods to the global minimum in different settings was shown, e.g., in~\cite{BachMoulines2013, JentzenKuckuckNeufeldVonWurstemberger2021, BachMoulines2011, Nesterov2015, Nesterov2004, Rakhlin2012, SchmidtleRoux2013}.

Typically,
the objective functions occurring in the training of ANNs with ReLU activation are  non-convex and, instead, 
admit infinitely many non-global local minima and saddle points.
In view of this,
it becomes important to study the landscapes of the risk functions in the training of ANNs and to develop an understanding of the appearance of critical points (such as non-global local extrema and saddle points) of the risk functions.
Recently, in the article Cheridito et al.~\cite{CheriditoJentzenRossmannek2021} a characterization of the saddle points and non-global
local minima of the risk function was obtained
for the case of affine target functions.
Sufficient conditions which ensure that the convergence of GD type optimization methods to saddle points can be excluded have been revealed, e.g., in~\cite{Ge2015,LeePanageasRecht2019, LeeJordanRecht2016, PanageasPiliouras2017, Panageas2019firstorder}.

Another promising direction of research is to study the convergence of GD type optimization methods for the training of ANNs in the so-called overparametrized regime, where the number of ANN parameters has to be sufficiently large when compared to the number of used input-output data pairs.
In this situation the risks of GD type optimization methods can be shown to converge to zero with high probability;
see, e.g., \cite{AroraDuHuLiWang2019, DuZhaiPoczosSingh2018, EMaWu2020, JentzenKroeger2021, LiLiang2019, RotskoffEijnden2018, ZhangMartensGrosse2019} for the case of ANNs with one hidden layer
and see, e.g., \cite{AllenzhuLiLiang2019, AllenzhuLiSong2019, DuLeeLiWangZhai2019, SankararamanDeXuHuangGoldstein2020, ZouCaoZhouGu2019} for the case of ANNs with more than one hidden layer. The results in these articles apply to the empirical risk, which is measured with respect to a finite set of input-output data pairs.

For convergence results for GD type optimization schemes 
without convexity but
under {\L}ojasiewicz type
assumptions we point, e.g., to \cite{AbsilMahonyAndrews2005,AttouchBolte2009,DereichKassing2021,Karimi2020linear,LeiHuLiTang2020,Wojtowytsch2021stochastic,XuYin2013}.
Further abstract convergence results for GD type optimization schemes
in the non-convex setting can be found, e.g., in \cite{AkyildizSabanis2021, BertsekasTsitsiklis2000, DereichMuller_Gronbach2019, FehrmanGessJentzen2020, LovasSabanis2020, Patel2021stopping} and the references mentioned therein.
In particular,
the article Fehrman et al.~\cite{FehrmanGessJentzen2020} shows convergence towards the global minimum value of
some GD type optimization algorithms with random initilizations,
provided that the set of global minima of the objective function is locally
a suitable submanifold of the parameter space
and provided that the Hessian of the objective function satisfies a certain maximal rank condition at these global minima.
A key contribution of this work is to demonstrate that these regularity assumptions are satisfied in the training of ANNs with one hidden layer and
ReLU activation provided that the target function is piecewise affine linear.

We also refer, e.g., to \cite{CheriditoJentzenRossmannek2020, JentzenvonWurstemberger2020, LuShinSuKarniadakis2020, Shamir2019} 
for lower bounds and divergence results for GD type optimization methods. 
For more detailed overviews and further literature on GD type optimization schemes 
we point, e.g., to \cite{BercuFort2013}, \cite{Bottou2018optimization}, 
\cite{EMaWojtowytschWu2020}, 
\cite[Section 1.1]{FehrmanGessJentzen2020}, 
\cite[Section 1]{JentzenKuckuckNeufeldVonWurstemberger2021}, and \cite{Ruder2017overview}.

There are different variants of GD type optimization methods in the scientific literature, such as the plain vanilla GD optimization method, GD optimization methods with momentum, and adaptive GD optimization methods (cf., e.g., Ruder~\cite{Ruder2017overview}),
and the plain vanilla GD optimization method with independent random initializations is maybe the GD based ANN training scheme which is most accessible for a mathematical convergence analysis. Despite the above mentioned promising mathematical analysis approaches in the literature, it remains -- even in the simple situation of the plain vanilla GD optimization method with independent random initializations and ANNs with one hidden layer and ReLU activation -- an open problem to prove (or disprove) the conjecture that the risk of the GD optimization method converges to the risk of the global minima of the risk function in the training of such ANNs. It is one of the key contributions of this article to prove this conjecture 
for the plain vanilla GD optimization method with independent random initializations and ANNs with one hidden layer and ReLU activation in the situation 
where the probability distribution of the input data is equivalent 
to the continuous uniform distribution on a compact interval with a Lipschitz continuous
density, where the probability distributions for the random initializations of the
ANN parameters are standard normal distributions, 
and where the target function under consideration is continuous and piecewise affine linear. 
The precise formulation of this statement is given in \cref{theo:intro:random:init} below 
within this introductory section.

In \cref{theo:intro:random:init} the target function (the function which describes the relationship between the input and the output data in the considered supervised learning problem) is described through the function $ f \colon [a,b] \to \R $ from the compact interval $ [a,b] $ to the real numbers $ \R $ where $ a, b \in \R $ are real numbers with $ a < b $. 
In \cref{theo:intro:random:init} this target function $ f \in C( [a,b], \R ) $ is assumed to be an 
element of the set $ C( [a,b], \R ) $ of continuous functions from $ [a,b] $ to $ \R $. 
In addition, in \cref{theo:intro:random:init} the target function $ f \colon [a,b] \to \R $ is assumed 
to be piecewise affine linear in the sense that there exist $ N \in \N $, $ \fx_0, \fx_1, \ldots, \fx_N \in \R $ with 
\begin{equation}
  a = \fx_0 < \fx_1 < ... < \fx_N = b
\end{equation}
so that for all $ i \in \{ 1, 2, \ldots , N \} $ 
we have that the target function 
$ [ \fx_{ i - 1 }, \fx_i ] \ni x \mapsto f(x) \in \R $ 
restricted to the subinterval $[ \fx_{ i - 1 }, \fx_i ]$ 
is affine linear; see above \cref{theo:intro:random:init:eq1} in \cref{theo:intro:random:init} below.

The risk functions associated to ANNs 
with ReLU activation fail to be continuously differentiable 
due to the lack of differentiability of the ReLU activation 
function 
$ \R \ni x \mapsto \max\{ x, 0 \} \in \R $
and, in view of this, 
one needs to introduce appropriate generalized gradients 
of the risk function which mathematically describe 
the behave of GD steps in implementations 
in numerical simulations 
to mathematically formulate the GD optimization method 
for the training of ANNs with ReLU activation. 
To accomplish this, we approximate as in \cite[(7) in Setting 2.1]{JentzenRiekert2021} and \cite[Theorem 1.1 and Proposition 2.3]{CheriditoJentzenRiekert2021}
the ReLU activation function 
$ \R \ni x \mapsto \max\{ x, 0 \} \in \R $
through appropriate continuously differentiable activation functions 
and then specify the generalized gradients as the limits 
of the usual gradients of the approximated risk functions; 
see \cref{eq:loss:gradient} in \cref{prop:loss:approximate:gradient} in \cref{subsection:setting} below. 
Specifically, in \cref{theo:intro:random:init} below the continuously differentiable 
functions 
$ \Rect_r \colon \R \to \R $, $ r \in \N $, 
serve as approximations for the ReLU activation function 
$ \Rect_{ \infty } \colon \R \to \R $ 
in the sense that for all $x \in \R$ it holds that 
$ \Rect_{ \infty }( x ) = \max\{ x, 0 \} $
and $\limsup\nolimits_{r \to \infty}  \rbr*{ \abs { \Rect_r ( x ) - \max \{ x , 0 \} } + \abs { (\Rect_r)' ( x ) - \indicator{(0, \infty)} ( x ) } } = 0$;
see \cref{theo:intro:random:init:eq1} in \cref{theo:intro:random:init} below.

In \cref{theo:intro:random:init} we also assume that the probability distribution of the input data
in the supervised learning problem considered in \cref{theo:intro:random:init} below 
is equivalent to the standard uniform distribution 
on $ [a,b] $ with a Lipschitz continuous density. 
More specifically, the Lipschitz continuous function 
$\dens \colon [a,b] \to (0,\infty)$ 
in \cref{theo:intro:random:init} is assumed to be an unnormalized density of the probability 
distribution of the input data with respect to the Lebesgue measure 
restricted to $ [a,b] $.

In \cref{theo:intro:random:init:eq2} in \cref{theo:intro:random:init} we consider 
fully connected feedforward ANNs with ReLU activation and three layers: 
one input layer with $ 1 $ neuron on the input layer (1-dimensional input), 
one hidden layer with $ \width \in \N $ neurons on the hidden layer ($ \width $-dimensional 
hidden layer), and one output layer with $ 1 $ neuron on the output layer 
(1-dimensional output). 
In particular, for every number $ \width \in \N $ of neurons 
on the hidden layer and every approximation parameter $ r \in \N \cup \{ \infty \} $
(see \cref{theo:intro:random:init:eq1} below)
we describe in \cref{theo:intro:random:init:eq2} below the risk function 
$ \cL^\width_r \colon \R^{3 \width + 1 } \to \R $
associated to the supervised learning problem considered in \cref{theo:intro:random:init}. 
The functions $\cG^\width \colon \R^{3 \width + 1 } \to \R^{3 \width + 1 }$, $ \width \in \N $, 
in \cref{theo:intro:random:init} specify generalized gradient functions 
of the risk functions $\cL^\width_{ \infty } \colon \R^{3 \width + 1 } \to \R$, $\width \in \N$, in \cref{theo:intro:random:init:eq2}.

For every number $ \width \in \N $ of neurons on the hidden layer, 
every natural number $ k \in \N $, 
and every learning rate $ \gamma \in \R $ 
we have that the 
random variables $\Theta^{ \width , k , \gamma }_n \colon \Omega \to \R^{3 \width + 1 }$, $n \in \N_0$, 
in \cref{theo:intro:random:init:eq3}
describe the GD process 
with learning rate $ \gamma $. 
Observe that the assumption in \cref{theo:intro:random:init} that for all $\width \in \N$, $\gamma \in \R$
it holds that 
$ \Theta^{ \width , k , \gamma }_0 \colon \Omega \to \R^{3 \width + 1 } $, $ k \in \N $,
are i.i.d.\ random variables 
ensures that
for all $ \width \in \N $, $ n \in \N_0 $, $ \gamma \in \R $
we have that the random variables 
$ \Theta^{ \width , k , \gamma }_n \colon \Omega \to \R^{3 \width + 1 }$,
$ k \in \N$,
are i.i.d.\ random variables. 
Loosely speaking, for every number $\width \in \N$ of neurons on the hidden layer, 
every natural number $ k \in \N $, 
every learning rate $ \gamma \in \R $, 
and every number $ n \in \N $ of GD steps we have that
the random variable ${\bfk}^{ \width, k, \gamma }_n \colon \Omega \to \N$ 
in \cref{theo:intro:random:init:eq4} selects an independent 
random initialization with the smallest risk.

Roughly speaking, in \cref{theo:intro:random:init:eqresult} in \cref{theo:intro:random:init} we prove that there exists a sufficiently small 
strictly positive real number $ \fg \in (0,\infty) $ 
such that for every learning rate $ \gamma \in (0,\fg ] $ which is smaller or equal 
than the strictly positive real number $ \fg $ 
we have as the number $ K \in \N $ of independent random realizations  
and the number $ \width \in \N $ of neurons on the hidden layer
increase to infinity 
convergence to one of the probability 
that the risk of the GD optimization method 
with independent standard normal random initializations converges to zero. 
We now present the precise statement of \cref{theo:intro:random:init} 
in a self-contained style and, thereafter, 
we outline how we prove \cref{theo:intro:random:init}.

\begin{samepage}
\begin{theorem} \label{theo:intro:random:init}
Let $ N \in \N$, $\fx_0, \fx_1, \ldots, \fx_N, a \in \R$, $b \in (a, \infty)$, $f \in C ( [a , b] , \R)$ satisfy $a = \fx_0 < \fx_1 < \cdots < \fx_N = b$, 
assume for all $i \in \{1, 2, \ldots, N\}$ that $f |_{[\fx_{i-1}, \fx_i ]}$ is affine linear,
let $\Rect_r \in C ( \R , \R )$, $r \in \N \cup \{ \infty\} $, satisfy for all $x \in \R$ that $( \bigcup_{r \in \N} \cu{ \Rect_r } ) \subseteq C^1( \R , \R)$, $\Rect_\infty ( x ) = \max \{ x , 0 \}$,
 $\sup_{r \in \N} \sup_{y \in [- \abs{x}, \abs{x} ] } \abs{ ( \Rect_r)'(y)} < \infty$, and
\begin{equation} \label{theo:intro:random:init:eq1}
    \limsup\nolimits_{r \to \infty}  \rbr*{ \abs { \Rect_r ( x ) - \Rect _\infty ( x ) } + \abs { (\Rect_r)' ( x ) - \indicator{(0, \infty)} ( x ) } } = 0,
\end{equation}
let $\dens \colon [a,b] \to (0, \infty)$ be Lipschitz continuous, 
let $\cL_r^\width \colon \R^{3 \width + 1} \to \R$, $r \in \N \cup \{ \infty\}$, $\width \in \N$,
satisfy for all $r \in \N \cup \{ \infty \}$, $\width \in \N$, $\theta = (\theta_1, \ldots, \theta_{3 \width + 1}) \in \R^{3 \width + 1}$ that
\begin{equation} \label{theo:intro:random:init:eq2}
    \cL_r^\width ( \theta ) = \int_a^b \rbr[\big]{ f ( x ) - \theta_{\fd} - \smallsum_{j=1}^\width \theta_{2 \width + j} \br{ \Rect_r ( \theta_{j} x + \theta_{\width + j } ) }  }^2 \dens ( x ) \, \d x,
\end{equation}
let $\cG^\width  \colon \R^{3 \width + 1} \to \R^{3 \width + 1}$, $\width \in \N$, satisfy for all $\width \in \N$,
$\theta \in \{ \vartheta \in \R^{3 \width + 1} \colon ( ( \nabla \cL^\width_r ) ( \vartheta ) ) _{r \in \N} \text{ is } \allowbreak \text{convergent} \}$
that $\cG^\width ( \theta ) = \lim_{r \to \infty} (\nabla \cL^\width_r) ( \theta )$,
let $(\Omega, \cF, \P)$ be a probability space,
let $\Theta^{\width , k , \gamma}_n  \colon \Omega \to  \R^{3 \width + 1 }$, $\width , k  \in \N$, $\gamma \in \R$, $n \in \N_0$,
and $\bfk^{\width , k , \gamma}_n \colon \Omega \to \N$, $\width , k  \in \N$, $\gamma \in \R$, $n \in \N_0$, be random variables,
assume for all $\width \in \N$, $\gamma \in \R$ that $\Theta_0^{\width , k , \gamma}$, $k \in \N$, are independent standard normal random vectors,
and assume for all $\width , k \in \N$, $\gamma \in \R$, $n \in \N_0$, $\omega \in \Omega$ that
\begin{equation} \label{theo:intro:random:init:eq3}
    \Theta_{n+1}^{\width , k , \gamma} ( \omega ) = \Theta_{n}^{\width , k , \gamma} ( \omega ) - \gamma \cG^\width (\Theta_{n}^{\width , k , \gamma } ( \omega ) ) 
\end{equation}
and
\begin{equation}  \label{theo:intro:random:init:eq4}
    \bfk^{\width , k , \gamma}_n ( \omega) \in \arg \min \nolimits_{\ell \in \{1, 2, \ldots, k \} } \cL_\infty^\width ( \Theta_n^{\width , \ell , \gamma } ( \omega ) ).
\end{equation}
Then there exists $\fg \in (0, \infty)$ such that for all $\gamma \in (0, \fg]$ it holds that
\begin{equation}  \label{theo:intro:random:init:eqresult}
    \liminf\nolimits_{\width \to \infty}  \liminf\nolimits_{K \to \infty} \P \rbr*{ \limsup\nolimits_{n \to \infty} \cL_\infty^\width  \rbr[\big]{ \Theta^{\width ,  \bfk^{\width , K , \gamma}_n, \gamma}_n } = 0  } = 1.
\end{equation}
\end{theorem}
\end{samepage}

\cref{theo:intro:random:init} is a direct consequence of \cref{cor:gd:random:init} below.
\cref{cor:gd:random:init}, in turn, follows from \cref{theo:gd:convergence} in \cref{subsection:gd:convergence} below, 
which is the main result of this article. Loosely speaking, 
\cref{theo:gd:convergence} establishes 
in the case of ANNs with three layers (1-dimensional input layer, 
$\width $-dimensional hidden layer, and 1-dimensional output layer)
and in the case of a continuous and piecewise 
affine linear target function $ f \colon [a,b] \to \R $ 
with $ N \in \N \cap [ 1 , \width ] $ grid points 
that there exists an appropriate open subset 
$ U \subseteq \R^{ \fd } $ 
of the ANN parameter space 
$ \R^{ \fd } = \R^{ 3 \width + 1 } $
such that for every sufficiently small 
learning rate $ \gamma \in (0,\infty) $ 
and every initial value $ \theta \in U $ 
it holds that 
the risk of the plain vanilla deterministic GD optimization method 
with initial value $ \theta $ and learning rate $ \gamma $ 
(see \cref{theo:gd:eq:def:descent} in \cref{theo:gd:convergence} in \cref{subsection:gd:convergence}) 
converges in the training of the considered ANNs exponentially quick to zero.

To make the statement of \cref{theo:gd:convergence} more accessible to the reader
within this introductory section,
we illustrate \cref{theo:gd:convergence}
by means of another consequence of \cref{theo:gd:convergence}
which is also of independent interest. Specifically, 
in \cref{theo:intro} below in this introductory section we prove 
in the case of ANNs with three layers (1-dimensional input layer, 
$\width$-dimensional hidden layer, and 1-dimensional output layer)
and 
in the case of a continuous and piecewise 
affine linear target function $ f \colon [a,b] \to \R $ 
with $ N \in \N \cap [ 1 , \width ] $ grid points 
that for every sufficiently small learning rate $ \gamma $ 
we have that 
the risk of the plain vanilla GD optimization method 
with learning rate $ \gamma $ and one standard normal random initialization 
(see \cref{theo:intro:eq:def:descent} in \cref{theo:intro}) 
converges exponentially to zero 
with strictly positive probability 
(see \cref{theo:intro:eq:conclusion} in \cref{theo:intro}). 
We now present the precise statement of \cref{theo:intro} and, thereafter, 
we briefly sketch how we prove 
\cref{theo:gd:convergence} in \cref{subsection:gd:convergence} and \cref{theo:intro}, respectively.

\begin{samepage}
\begin{theorem} \label{theo:intro}
Let $ \width, \fd \in \N$, $N \in \N \cap [ 1 , \width ]$, $\fx_0, \fx_1, \ldots, \fx_N, a \in \R$, $b \in (a, \infty)$, $f \in C ( [a , b] , \R)$ satisfy $\fd = 3 \width + 1$ and $a = \fx_0 < \fx_1 < \cdots < \fx_N = b$, 
assume for all $i \in \{1, 2, \ldots, N\}$ that $f |_{[\fx_{i-1}, \fx_i ]}$ is affine linear,
let $\Rect_r \in C ( \R , \R )$, $r \in \N \cup \{ \infty\} $, satisfy for all $x \in \R$ that $( \bigcup_{r \in \N} \cu{ \Rect_r } ) \subseteq C^1( \R , \R)$, $\Rect_\infty ( x ) = \max \{ x , 0 \}$,
 $\sup_{r \in \N} \sup_{y \in [- \abs{x}, \abs{x} ] }  \abs{ ( \Rect_r)' ( y ) } < \infty$, and
\begin{equation}
    \limsup\nolimits_{r \to \infty}  \rbr*{ \abs { \Rect_r ( x ) - \Rect _\infty ( x ) } + \abs { (\Rect_r)' ( x ) - \indicator{(0, \infty)} ( x ) } } = 0,
\end{equation}
let $\dens \colon [a,b] \to (0, \infty)$ be Lipschitz continuous,
let $\cL_r \colon \R^\fd \to \R$, $r \in \N \cup \{ \infty\}$,
satisfy for all $r \in \N \cup \{ \infty \}$, $\theta = (\theta_1, \ldots, \theta_\fd) \in \R^{\fd}$ that
\begin{equation} \label{theo:intro:def:risk}
    \cL_r ( \theta ) = \int_a^b \rbr[\big]{ f ( x ) - \theta_{\fd} - \smallsum_{j=1}^\width \theta_{2 \width + j} \br{ \Rect_r ( \theta_{j} x + \theta_{\width + j } ) }  }^2 \dens ( x ) \, \d x,
\end{equation}
let $(\Omega, \cF, \P)$ be a probability space,
let $\Theta^\gamma_n \colon \Omega \to  \R^{\fd }$, $\gamma \in \R$, $n \in \N_0$, be random variables,
assume for every $\gamma \in \R$ that $\Theta_0^\gamma $ is standard normally distributed, 
let $\cG \colon \R^\fd \to \R^\fd$ satisfy for all
$\theta \in \{ \vartheta \in \R^\fd \colon ( ( \nabla \cL_r ) ( \vartheta ) ) _{r \in \N} \text{ is convergent} \}$
that $\cG ( \theta ) = \lim_{r \to \infty} (\nabla \cL_r) ( \theta )$,
and assume for all $\gamma \in \R$, $n \in \N_0$, $\omega \in \Omega$ that
\begin{equation} \label{theo:intro:eq:def:descent}
    \Theta_{n+1}^\gamma ( \omega ) = \Theta_{n}^\gamma ( \omega ) - \gamma \cG (\Theta_{n}^\gamma ( \omega ) ) .
\end{equation}
Then there exist $\fc, \fC \in (0, \infty)$ such that for all $\gamma \in (0, \fc]$ it holds that 
\begin{equation} \label{theo:intro:eq:conclusion}
   \P \rbr[\big]{ \limsup \nolimits_{n \to \infty} \cL_\infty ( \Theta_n^\gamma) = 0 } \geq  \P \rbr[\big]{ \forall \, n \in \N_0 \colon \cL _\infty ( \Theta_n^\gamma ) \leq \fC \exp ( - \fc \gamma n ) } \geq \fc > 0.
\end{equation}
\end{theorem}
\end{samepage}

\cref{theo:intro} is an immediate consequence of \cref{cor:gd:convergence:random} below (applied with $\rho \with 0$ in the notation of \cref{cor:gd:convergence:random}).
 \cref{cor:gd:convergence:random}, in turn, is a direct consequence of \cref{theo:gd:convergence} (see \cref{subsection:gd:random} below for details). Roughly speaking, we prove \cref{theo:intro:random:init}, \cref{theo:intro}, and \cref{theo:gd:convergence}, respectively, 
 (i) by showing that for every number $ \width \in \N \cap [ N, \infty ) $ of neurons on the hidden layer there exists a natural number $\mdim \in \N \cap [1, \fd ) $ such that a suitable subset of the set of global minima of the risk function $\cL_{ \infty } \colon \R^{ \fd } \to \R$ in \cref{theo:intro:def:risk} is a twice continuously differentiable $k$-dimensional submanifold of the ANN parameter space $\R^{ \fd } = \R^{ 3 \width + 1 }$ (cf.~\cref{lem:manifold:help} and \cref{cor:minima:manifold:1} in \cref{section:risk:global:min} below),
 (ii) by proving that the ranks of the Hessian matrices of the risk function on this suitable set of global minima of the risk function $\cL_{ \infty } \colon \R^{ \fd } \to \R$ in \cref{theo:intro:def:risk} are equal to $\fd - k $,
and, thereafter, (iii) by applying the machinery in Fehrman et al.~\cite{FehrmanGessJentzen2020} to establish convergence of the GD optimization method.

The remainder of this article is organized as follows. In \cref{section:diff} we establish several regularity properties for the Hessian matrix of the risk function of the considered supervised learning problem. In \cref{section:risk:global:min} we employ the findings from \cref{section:diff} to establish that a suitable subset of the set of global minima of the risk function constitutes a $C^{ \infty }$-submanifold of the ANN parameter space $\R^{ \fd } = \R^{ 3 \width + 1 }$ on which the Hessian matrix of the risk function has maximal rank.
 In \cref{section:flow:convergence} we engage the findings from \cref{section:risk:global:min} to establish that the risk of certain solutions of GF differential equations converges exponentially quick to zero. Finally, in \cref{section:gd:convergence} we establish that the risk of certain GD processes converges exponentially quick to zero and, thereby, we also prove \cref{theo:intro:random:init,theo:intro} above.

\section{Second order differentiability properties of the risk function} \label{section:diff}

In this section we establish in \cref{lem:twice:diff} in \cref{subsection:risk:diff} below an explicit representation result for the Hessian matrix of the risk function of the considered supervised learning problem.
In particular, in \cref{lem:twice:diff} we identify a suitable open subset of the ANN parameter space with full Lebesgue measure on which the risk function is twice continuously differentiable (see \cref{setting:eq:defv} below for details).
This is nontrivial due to the fact that the ReLU activation function 
$ \R \ni x \mapsto \max \cu{x,0} \in \R $ is not everywhere differentiable.
Results related to \cref{lem:twice:diff} have been shown in \cite[Lemma 3.8]{CheriditoJentzenRossmannek2021}.

 \cref{cor:second:derivatives:globalmin} in \cref{subsection:risk:diff} specializes \cref{lem:twice:diff} to the specific situation where the ANN parameter represents a global minima of the risk function. In \cref{lem:twice:diff:lipschitz} in \cref{subsection:risk:diff} we employ \cref{lem:twice:diff} to conclude under the assumption that the target function is locally Lip\-schitz continuous that the second derivative of the risk function is locally Lipschitz continuous. In \cref{lem:hessian:max:upperbound}, \cref{lem:realization:upper:bound}, and \cref{cor:hessian:upperbound:max} in \cref{subsection:hess:bound} below we use \cref{lem:twice:diff} to derive suitable upper bounds for the absolute values of the second order partial derivatives of the risk function. \cref{lem:twice:diff:lipschitz}, \cref{cor:second:derivatives:globalmin}, and \cref{cor:hessian:upperbound:max} are all employed in \cref{section:risk:global:min} below. 

Our proof of \cref{lem:twice:diff} employs the well-known Leibniz integral rule type result in \cref{lem:leibniz:rule} in \cref{subsection:risk:diff}, the known representation and regularity results for the first derivative of the risk function in \cref{prop:loss:approximate:gradient} in \cref{subsection:setting} below and \cref{prop:loss:differentiable} in \cref{subsection:risk:diff}, the elementary continuity result in \cref{lem:indicator:terms:cont} in \cref{subsection:risk:diff}, 
the elementary and well-known differentiability results for certain parameter integrals in \cref{lem:interchange} and \cref{cor:interchange} in \cref{subsection:lipschitz:integrals} below, and the elementary continuity result for certain para\-meter integrals involving indicator functions in \cref{lem:integral:interval:continuous} in \cref{subsection:lipschitz:integrals} and \cref{cor:derivative:integral:cont} in \cref{subsection:intervals:lipschitz} below.
 \cref{prop:loss:differentiable} is a direct consequence of Proposition 2.11 in \cite{JentzenRiekertFlow} and \cref{prop:loss:approximate:gradient} follows directly from, e.g., item (iv) in Proposition 2.2 in \cite{JentzenRiekertFlow}. Our proof of \cref{lem:twice:diff:lipschitz} also uses the local Lipschitz continuity results for certain parameter integrals involving indicator functions in \cref{cor:derivative:integral:lip} in \cref{subsection:intervals:lipschitz}.
 Our proofs of \cref{cor:derivative:integral:cont,cor:derivative:integral:lip}, in turn, employ the elementary Lipschitz continuity result for certain parameter integrals involving indicator functions in \cref{lem:integral:interval:lipschitz} in \cref{subsection:lipschitz:integrals} as well as the local Lipschitz continuity results for active neuron regions in \cref{lem:active:intervals} and \cref{cor:active:intervals} in \cref{subsection:intervals:lipschitz}.

\subsection{Mathematical description of artificial neural networks (ANNs)} \label{subsection:setting}
\begin{setting} \label{setting:snn}
Let $ \width, \fd \in \N$, $a \in \R$, $b \in (a, \infty)$, $f \in C ( [a , b] , \R)$ satisfy $\fd = 3 \width + 1$,
let $\fw  = (( \w{\theta} _ {1}, \ldots, \w{\theta}_{\width}  ))_{ \theta \in \R^{\fd}} \colon \R^{\fd} \to \R^{\width}$,
$\fb =  (( \b{\theta} _ 1, \ldots, \b{\theta}_{\width} ) )_{ \theta \in \R^{\fd}} \colon \R^{\fd} \to \R^{\width}$,
$\fv = (( \v{\theta} _ 1, \ldots, \v{\theta}_{\width} ) )_{ \theta \in \R^{\fd}} \colon \R^{\fd} \to \R^{\width}$, 
$\fc = (\c{\theta})_{\theta \in \R^{\fd }} \colon \R^{\fd} \to \R$, and
$\fq = ( ( \q{\theta}_1, \ldots, \q{\theta}_{\width}) ) \colon \R^\fd \to (- \infty , \infty] ^\width$
 satisfy for all $\theta  = ( \theta_1 ,  \ldots, \theta_{\fd}) \in \R^{\fd}$, $j \in \{1, 2, \ldots, \width \}$ that $\w{\theta}_{ j} = \theta_{ j}$, $\b{\theta}_j = \theta_{\width + j}$, 
$\v{\theta}_j = \theta_{2 \width + j}$, 
$\c{\theta} = \theta_{\fd}$, and
\begin{equation} \label{setting:eq:defqi}
    \q{\theta}_j = \begin{cases}
    - \nicefrac{\b{\theta}_j}{\w{\theta}_j} & \colon \w{\theta}_j \not= 0 \\
    \infty & \colon \w{\theta}_j = 0,
    \end{cases}
\end{equation}
let $\dens \colon [a,b] \to (0, \infty)$ be Lipschitz continuous,
let $\Rect \colon \R \to \R$,
$\scrN = (\realization{\theta})_{\theta \in \R^{\fd } } \colon \R^{\fd } \to C(\R , \R)$, and $ \cL \colon \R^{\fd  } \to \R$
satisfy for all $\theta \in \R^{\fd}$, $x \in \R$ that $\Rect ( x ) = \max \{ x , 0 \}$, $\realization{\theta} (x) = \c{\theta} + \smallsum_{j=1}^\width \v{\theta}_j \br{ \Rect \rbr{ \w{\theta}_{j} x + \b{\theta}_j } }$, and
\begin{equation} \label{setting:eq:def:risk}
\cL (\theta) = \int_{a}^b (\realization{\theta} (y) - f ( y ) )^2 \dens ( y ) \, \d y ,
\end{equation}
let $\chi_r \in C^1 ( \R , \R )$, $r \in \N $, satisfy for all $x \in \R$ that $\sup_{r \in \N} \sup_{y \in [- \abs{x}, \abs{x} ] }  \abs{ ( \chi_r ) ' ( y ) } < \infty$ and
\begin{equation}
    \limsup\nolimits_{r \to \infty}  \rbr*{ \abs { \chi_r ( x ) - \Rect ( x ) } + \abs { (\chi_r)' ( x ) - \indicator{(0, \infty)} ( x ) } } = 0,
\end{equation}
let $\fL_r \colon \R^\fd \to \R$, $r \in \N$,
satisfy for all $r \in \N $, $\theta \in \R^{\fd}$ that
\begin{equation}
   \fL_r ( \theta ) = \int_a^b \rbr[\big]{ f ( y ) - \c{\theta} - \smallsum_{j=1}^\width \v{\theta}_j \br[\big]{ \chi_r ( \w{\theta}_{j} y + \b{\theta}_j ) } } ^2 \dens ( y ) \, \d y,
\end{equation}
let $I_j^\theta \subseteq \R$, $\theta \in \R^{\fd }$, $j \in \{1, 2, \ldots, \width \}$, satisfy for all 
$\theta \in \R^{\fd}$, $j \in \{1, 2, \ldots, \width \}$ that $I_j^\theta = \{ x \in [a , b ] \colon \w{\theta}_{j} x + \b{\theta}_j  > 0 \}$,
let $\cG = ( \cG_1 , \ldots, \cG_\fd ) \colon \R^\fd \to \R^\fd$ satisfy for all
$\theta \in \{ \vartheta \in \R^\fd \colon ( ( \nabla \fL_r ) ( \vartheta ) ) _{r \in \N} \, \allowbreak \text{is convergent} \}$
that $\cG ( \theta ) = \lim_{r \to \infty} (\nabla \fL_r) ( \theta )$, and let $\fV \subseteq \R^\fd$ satisfy
\begin{equation} \label{setting:eq:defv}
    \fV = \cu[\big]{\theta \in \R^\fd \colon \rbr[\big]{ {\textstyle\prod_{j=1}^\width} {\textstyle\prod_{v \in \{a, b\} } } (\w{\theta}_j v + \b{\theta}_j )  \not= 0 } }.
\end{equation}
\end{setting}

\begin{prop} \label{prop:loss:approximate:gradient}
Assume \cref{setting:snn}. Then it holds for all $\theta \in \R^{\fd}$, $i \in \{1, 2, \ldots, \width \}$ that
\begin{equation} \label{eq:loss:gradient}
\begin{split}
        \cG_{ i} ( \theta) &= 2 \v{\theta}_i \int_{I_i^\theta} x  ( \realization{\theta} (x) - f ( x ) ) \dens ( x ) \, \d x , \\
        \cG_{\width + i} ( \theta) &= 2 \v{\theta}_i \int_{I_i^\theta} (\realization{\theta} (x) - f ( x ) ) \dens ( x ) \,  \d x , \\
        \cG_{2 \width + i} ( \theta) &= 2 \int_{a}^b \br[\big]{\Rect \rbr{ \w{\theta}_{i} x + \b{\theta}_i } } ( \realization{\theta}(x) - f ( x ) ) \dens ( x ) \, \d x , \\
        \text{and} \qquad \cG_{\fd} ( \theta) &= 2 \int_{a}^b (\realization{\theta} (x) - f ( x ) ) \dens ( x ) \, \d x .
        \end{split}
\end{equation}
\end{prop}
\begin{cproof}{prop:loss:approximate:gradient}
\Nobs that, e.g., \cite[Item (iv) in Proposition 2.2]{JentzenRiekertFlow} establishes \cref{eq:loss:gradient}.
\end{cproof}

\subsection{Regularity properties for parametric integrals 
of Lipschitz continuous functions} \label{subsection:lipschitz:integrals}

\begin{lemma} \label{lem:interchange}
Let $\fu \in \R$, $\fv \in (\fu , \infty)$,
let $\phi \colon \R \times [\fu , \fv] \to \R$ be locally bounded and measurable,
let $\mu \colon \cB( [ \fu , \fv ] ) \to [0, \infty ]$ be a finite measure,
let $\Phi  \colon \R \to \R$ satisfy for all $x \in \R$ that
\begin{equation}
    \Phi (x) = \int_\fu ^\fv \phi (x,s) \, \mu ( \d s ) ,
\end{equation}
let $x \in \R$, 
 $\delta, c \in (0, \infty)$ satisfy for all $h \in (- \delta, \delta)$, $s \in [\fu , \fv]$ that $\abs{\phi ( x + h , s ) - \phi ( x , s ) } \leq c \abs{h}$,
let $E \subseteq [\fu , \fv]$ be measurable, assume $\mu ( [\fu , \fv] \backslash E  ) = 0 $, and assume for all $s \in E$ that $\R \ni v \mapsto \phi  ( v , s ) \in \R$ is differentiable at $x$. Then
\begin{enumerate} [label=(\roman*)]
    \item it holds that $\Phi $ is differentiable at $x$ and
    \item it holds that
    \begin{equation}
        \Phi '(x) = \int_E \rbr[\big]{\tfrac{\partial}{\partial x} \phi  }  ( x , s ) \, \mu ( \d s ).
    \end{equation}
\end{enumerate}
\end{lemma}
\begin{cproof2}{lem:interchange}
 \Nobs that the assumption that $\mu ( [\fu , \fv] \backslash E ) = 0 $ shows for all $h \in \R \backslash \{ 0 \}$ that
\begin{equation} \label{lem:interchange:eq1}
\begin{split}
     h^{-1} [ \Phi (x+h) - \Phi (x)] &= \int_\fu^\fv h^{-1} [ \phi (x+h , s ) - \phi (x,s)] \, \mu ( \d s ) \\
     &= \int_E h^{-1} [ \phi (x+h , s ) - \phi (x,s)] \, \mu ( \d s ) .
\end{split}
\end{equation}
Next \nobs that the assumption that for all $s \in E$ it holds that $\R \ni v \mapsto \phi  ( v , s) \in \R$ is differentiable at $x$ ensures that for all $s \in E$ it holds that
\begin{equation} \label{lem:interchange:eq2}
    \lim\nolimits_{\R \backslash \cu{0} \ni h \to 0} \rbr*{ h^{-1} \br{ \phi ( x+h, s) - \phi (x, s)}} = \rbr[\big]{\tfrac{\partial}{\partial x} \phi  }(x, s ).
\end{equation}
Moreover, \nobs that the assumption that for all $h \in (- \delta, \delta)$, $s \in [\fu , \fv]$ it holds that $\abs{\phi ( x + h , s ) - \phi ( x , s ) } \leq c \abs{h}$ implies that for all $h \in ( -\delta, \delta) \backslash \{0 \}$, $s \in [ \fu , \fv]$ we have that $\abs{ h^{-1} \br{ \phi ( x+h, s) - \phi (x, s)} } \leq c$. Combining this with \cref{lem:interchange:eq1}, \cref{lem:interchange:eq2}, and the dominated convergence theorem demonstrates that
\begin{equation}
    \begin{split}
       & \lim\nolimits_{ \R \backslash \cu{0} \ni h \to 0 } \rbr*{ h^{-1}  \br{\Phi (x+h) - \Phi (x)}} \\
        &= \int_E \br*{ \lim\nolimits_{ \R \backslash \cu{0} \ni h \to 0 }  \rbr*{ h^{-1} \br{ \phi ( x+h, s) - \phi (x, s)}  }} \, \mu ( \d s ) 
        = \int_E \rbr[\big]{\tfrac{\partial}{\partial x} \phi  }(x, s ) \, \mu ( \d s ) .
    \end{split}
\end{equation}
\end{cproof2}

\cfclear
\begin{cor} \label{cor:interchange}
Let $n \in \N$, $j \in \{1, 2, \ldots, n \}$, $\fu \in \R$, $\fv \in (\fu , \infty)$,
let $\phi \colon \R^n \times [\fu , \fv] \to \R$ be locally bounded and measurable, 
let $\mu \colon \cB( [ \fu , \fv ] ) \to [0, \infty ]$ be a finite measure,
let $\Phi  \colon \R^n \to \R$ satisfy for all $x \in \R^n$ that
\begin{equation}
    \Phi (x) = \int_\fu ^\fv \phi  (x , s ) \, \mu (  \d s ) ,
\end{equation}
let $x = (x_1, \ldots, x_n) \in \R^n$,
 $\delta, c \in (0, \infty)$ satisfy for all $s \in [\fu , \fv]$, $h \in (- \delta , \delta )$ that 
 \begin{equation}
     \abs{\phi ( x_1, \ldots, x_{j-1}, x_j + h , x_{j+1}, \ldots, x_n , s ) - \phi ( x , s ) } \leq c \abs{h} ,
 \end{equation}
let $E \subseteq [\fu , \fv ]$ be measurable, assume $\mu ( [\fu , \fv] \backslash E ) = 0 $, and assume for all $s \in E$ that $\R \ni v \mapsto \phi  ( x_1, \ldots, x_{j-1}, v, x_{j+1}, \ldots, x_n, s) \in \R$ is differentiable at $x_j$. Then
\begin{enumerate} [label=(\roman*)]
    \item \label{cor:interchange:item1} it holds that $\R \ni v \mapsto \Phi (x_1, \ldots, x_{j-1}, v , x_{j+1}, \ldots, x_n) \in \R$ is differentiable at $x_j$ and
    \item \label{cor:interchange:item2} it holds that
    \begin{equation}
        \rbr[\big]{\tfrac{\partial}{\partial x_j} \Phi } ( x_1, \ldots, x_n) = \int_E \rbr[\big]{\tfrac{\partial}{\partial x_j} \phi  }(x_1, \ldots, x_n, s ) \, \mu ( \d s ).
    \end{equation}
\end{enumerate}
\end{cor}
\begin{cproof}{cor:interchange}
\Nobs that \cref{lem:interchange} establishes items \ref{cor:interchange:item1} and \ref{cor:interchange:item2}.
\end{cproof}

\begin{definition} \label{def:norm}
We denote by $\norm{ \cdot } \colon \rbr*{  \bigcup_{n \in \N} \R^n  } \to \R $ and $\spro{\cdot, \cdot} \colon \rbr*{ \bigcup_{n \in \N} (\R^n \times \R^n )  } \to  \R$ the functions which satisfy for all
$n \in \N$, $x=(x_1, \ldots, x_n)$, $y = (y_1, \ldots, y_n) \in \R^n $ that
$\norm{ x } = \br[\big]{ \sum_{i=1}^n \abs*{ x_i } ^2 } ^{1/2}$ and $\spro{x,y} = \sum_{i=1}^n x_i y_i$.
\end{definition}

\cfclear 
\begin{lemma} \label{lem:integral:interval:continuous}
Let $n \in \N$, $\fu \in \R$, $\fv \in (\fu , \infty)$, $x \in \R^n$, $c , \varepsilon \in (0, \infty)$, 
$\phi \in C( \R^n \times [\fu , \fv] , \R )$, 
let $\mu \colon \cB( [ \fu , \fv ] ) \to [0, \infty ]$ be a finite measure,
let $I^y \in \cB ([\fu, \fv])$, $y \in \R^n$, satisfy for all $y , z \in \{ v \in \R^n \colon \norm{x-v} \leq \varepsilon \}$ that $\mu ( I^y \Delta I^z ) \leq c \norm{ y - z }$,
and let $\Phi \colon \R^n \to \R$ satisfy for all $y \in \R^n$ that
\begin{equation} \cfadd{def:norm} 
    \Phi ( y ) = \int_{I^y} \phi (y , s ) \, \mu ( \d s )
\end{equation}
\cfload. Then 
it holds that $\{ v \in \R^n \colon \norm{x - v} \leq \varepsilon \} \ni y \mapsto \Phi ( y ) \in \R$ is continuous.
\end{lemma}
\begin{cproof} {lem:integral:interval:continuous}
Throughout this proof let $y \in \{ v \in \R^n \colon \norm{x - v} \leq \varepsilon \}$
and let $z = (z_k )_{k \in \N} \colon \N \to \{ v \in \R^n \colon \norm{x - v} \leq \varepsilon \}$ satisfy $\limsup_{k \to \infty} \norm{z_k - y } = 0$. \Nobs that for all $k \in \N$ it holds that
\begin{equation} \label{lem:integral:interval:continuous:eq1}
\begin{split}
    \abs{\Phi (y) - \Phi (z_k )} &\leq \int_{I^y \cap I^{z_k } } \abs{\phi (y,s) - \phi (z_k , s ) } \, \mu ( \d s ) + \int_{I^y \backslash I^{z_k } } \abs{\phi (y , s)} \, \mu ( \d s ) \\
    & \quad + \int_{I^{z_k } \backslash I^y} \abs{\phi (z_k , s ) } \, \mu (  \d s ) .
\end{split}
\end{equation}
Next \nobs that the assumption that $\phi$ is continuous and the dominated convergence theorem demonstrate that
\begin{equation} \label{lem:integral:interval:continuous:eq2}
  \limsup_{k \to \infty} \br*{ \int_{I^y \cap I^{z_k } } \abs{\phi (y,s) - \phi (z_k ,s) } \, \mu ( \d s ) } = 0.
\end{equation}
Moreover, \nobs that the fact that for all $k \in \N$ it holds that $\mu ( I^y \Delta I^{z_k } ) \leq c \norm{y - z_k }$ and the assumption that $\phi$ is continuous prove that for all $k \in \N$ we have that
\begin{equation}
    \limsup_{k \to \infty} \br*{ \int_{I^y \backslash I^{ z_k } } \abs{\phi (y , s)} \, \mu ( \d s ) + \int_{I^{z_k} \backslash I^y} \abs{\phi (z_k , s ) } \, \mu (  \d s )} = 0.
\end{equation}
Combining this with \cref{lem:integral:interval:continuous:eq1,lem:integral:interval:continuous:eq2} establishes that $\limsup_{k \to \infty} \abs{\Phi ( y ) - \Phi ( z_k ) } = 0$.
\end{cproof}

\cfclear
\begin{lemma} \label{lem:integral:interval:lipschitz}
Let $n \in \N$, $\fu \in \R$, $\fv \in (\fu , \infty)$, $x \in \R^n$, 
$c, \varepsilon \in (0, \infty)$, let $\phi \colon \R^n \times [\fu , \fv] \to \R$ be locally Lipschitz continuous, 
let $\mu \colon \cB( [ \fu , \fv ] ) \to [0, \infty ]$ be a finite measure,
let $I^y \in \cB ([\fu, \fv])$, $y \in \R^n$, satisfy for all $y , z \in \{ v \in \R^n \colon \norm{x-v} \leq \varepsilon \}$ that $\mu ( I^y \Delta I^z ) \leq c \norm{ y - z }$,
and let $\Phi \colon \R^n \to \R$ satisfy for all $y \in \R^n$ that
\begin{equation} \label{lem:integral:interval:lip:defphi} \cfadd{def:norm} 
    \Phi ( y ) = \int_{I^y} \phi (y , s ) \, \mu ( \d s )
\end{equation}
\cfload. Then there exists $\fC \in \R$ such that for all $y , z \in \{ v \in \R^n \colon \norm{x - v} \leq \varepsilon \}$ 
it holds that $\abs{\Phi (y) - \Phi (z) } \leq \fC \norm{y - z}$.
\end{lemma}
\begin{cproof}{lem:integral:interval:lipschitz}
\Nobs that the assumption that $\phi $ is locally Lipschitz continuous ensures that there exists $\fC \in \R$ which satisfies for all $y , z \in \{ v \in \R^n \colon \norm{x-v} \leq \varepsilon \}$, $s \in [\fu , \fv]$ with $y \not= z$ that
\begin{equation} \label{lem:integral:interval:lipschitz:eq1}
  \tfrac{  \abs{\phi (y , s) - \phi (z , s) }} {\norm{y - z}} + \abs{\phi (y,s)} + \abs{\phi (z,s)} \leq \fC. 
\end{equation}
Furthermore, \nobs that \cref{lem:integral:interval:lip:defphi} ensures for all $y , z \in \R^n$ that
\begin{equation} \label{lem:integral:interval:lipschitz:eq2}
    \abs{\Phi (y) - \Phi (z)} \leq \int_{I^y \cap I^z } \abs{\phi (y,s) - \phi (z,s) } \, \mu ( \d s ) + \int_{I^y \backslash I^z} \abs{\phi (y , s)} \, \mu ( \d s ) + \int_{I^z \backslash I^y} \abs{\phi (z,s)} \, \mu (  \d s ) .
\end{equation}
In addition, \nobs that \cref{lem:integral:interval:lipschitz:eq1} shows for all $y , z \in \{ v \in \R^n \colon \norm{x-v} \leq \varepsilon \}$ that
\begin{equation} \label{lem:integral:interval:lipschitz:eq3}
    \int_{I^y \cap I^z } \abs{\phi (y , s) - \phi (z , s) } \, \mu ( \d s ) \leq \fC  \norm{y - z} \mu ( [ \fu , \fv ] ).
\end{equation}
Moreover, \nobs that \cref{lem:integral:interval:lipschitz:eq1} and the assumption that for all $y , z \in \{ v \in \R^n \colon \norm{x-v} \leq \varepsilon \}$ it holds that $\mu ( I^y \Delta I^z ) \leq c \norm{y - z}$ prove that for all $y , z \in \{ v \in \R^n \colon \norm{x-v} \leq \varepsilon \}$ we have that
\begin{equation} \label{lem:integral:interval:lipschitz:eq4}
    \int_{I^y \backslash I^z} \abs{\phi (y,s)} \, \mu( \d s ) + \int_{I^z \backslash I^y} \abs{\phi (z,s)} \, \mu ( \d s ) \leq  c \fC \norm{y - z}.
\end{equation}
Combining this with
\cref{lem:integral:interval:lipschitz:eq2} and \cref{lem:integral:interval:lipschitz:eq3} establishes for all $y , z\in \{ v \in \R^n \colon \norm{x - v} \leq \varepsilon \}$ that
\begin{equation}
    \abs{\Phi (y) - \Phi (z)} \leq \fC (  c + \mu ( [ \fu , \fv ] ) ) \norm{y - z}.
\end{equation}
\end{cproof}

\subsection{Local Lipschitz continuity for active neuron regions} \label{subsection:intervals:lipschitz}

\cfclear
\begin{lemma} \label{lem:active:intervals} 
Let $a \in \R$, $b \in (a, \infty)$, $u = (u_1, u_2 ) \in \R^2 \backslash \cu{0}$,
let $\dens \colon [a,b] \to \R$ be bounded and measurable, 
and let $I^v \subseteq \R$, $v \in \R^2$, satisfy for all $v = (v_1 , v_2 ) \in \R^2$ that 
$I^v = \{ x \in [a, b ] \colon v_1 x + v_2 > 0 \}$.
Then there exist $c, \varepsilon \in (0, \infty)$ such that for all $ v , w \in \R^2$ with $\max \{ \norm{u - v} , \norm{u - w } \} \leq \varepsilon$ it holds that 
\begin{equation} \label{lem:active:intervals:claim} \cfadd{def:norm}
    \abs*{ \int _{ I^v \Delta I ^w } \dens ( x ) \, \d x } \leq c \norm{v - w}
\end{equation}
\cfload.
\end{lemma}
\begin{cproof} {lem:active:intervals}
Throughout this proof let $ M \in \R$ satisfy 
$M = \sup_{x \in [a,b] } \abs{ \dens ( x ) } $.
In the following we distinguish between the case $u_1 = 0 $ and the case $u_1 \not= 0$.

We first prove \cref{lem:active:intervals:claim} in the case \begin{equation} \label{lem:active:intervals:eq:case1}
   u_1 = 0 .
\end{equation} 
\Nobs that \cref{lem:active:intervals:eq:case1} and the assumption that $u =(u_1, u_2) \in \R^{2} \backslash \cu{ 0 }$ imply that $u_2 \not= 0$.
Moreover, \nobs that \cref{lem:active:intervals:eq:case1} shows for all $v = (v_1, v_2) \in \R^2$,
$x \in I^u \Delta I^v$ that
\begin{equation} \label{proof:active:intervals:eq1}
    \abs*{(u_1 x + u_2) - (v_1 x + v_2 )} 
    =\abs*{u_1 x + u_2 }  + \abs*{v_1 x + v_2 } 
    \geq \abs*{ u_1 x + u_2 } = \abs{ u_2 }.
\end{equation}
In addition, \nobs that for all $v = (v_1, v_2) \in \R^{2}$, $x \in [a,b]$ we have that
\begin{equation} \label{proof:active:intervals:eq2}
   \abs*{(u_1 x + u_2) - (v_1 x + v_2 )} \leq \abs{u_1 - v_1} \abs{x} + \abs{u_2 - v_2}
         \leq ( 1 + \max \cu{ \abs{a} , \abs{b} } ) \norm{u - v }.
\end{equation}
Combining this with \cref{proof:active:intervals:eq1} demonstrates for all $v \in \R^2 $ with $\norm{ u - v } < \frac{\abs{u_2 }}{1 + \max \cu{\abs{a}, \abs{b} }}$ that $I^u \Delta I^v = \emptyset$
and, therefore, $I^u = I^v$. 
Hence, we obtain
for all $ v , w \in \R^2$ with $\max \{ \norm{u - v} , \norm{u - w } \} \leq \frac{\abs{u_2 }}{2 + \max \cu{\abs{a}, \abs{b} }}$ that $I^v = I^u = I^w $ and, therefore, $ \int _{ I^v \Delta I ^w } \dens ( x ) \, \d x = 0$.
This establishes \cref{lem:active:intervals:claim} in the case $u_1=0$.

In the next step we prove \cref{lem:active:intervals:claim} in the case 
$u_1 \not= 0$.
\Nobs that for all $v = (v_1, v_2)$, $w = (w_1, w_2 ) \in \R^{2}$, $\fs \in \cu{ - 1 , 1 }$ with $\min \cu{\fs v_1, \fs w_1 } > 0$ it holds that
\begin{equation}
    \begin{split}
        I^v \backslash I^w &= \cu*{ y \in [ a , b ] \colon v_1 y + v_2 > 0 \geq w_1 y + w_2} = \cu*{ y \in [ a , b ] \colon - \tfrac{\fs v_2}{v_1} < \fs y \leq - \tfrac{s w_2}{w_1} } \\
        & \subseteq \cu*{ y \in \R \colon - \tfrac{\fs v_2}{v_1} < \fs y \leq - \tfrac{s w_2}{w_1} }.
    \end{split}
\end{equation}
Hence, we obtain for all $v = (v_1, v_2)$, $w = (w_1, w_2 ) \in \R^{2}$, $\fs \in \cu{ - 1 , 1 }$ with $\min \cu{\fs v_1, \fs w_1 } > 0$ that 
\begin{equation} \label{lem:active:intervals:eq:help3}
  \int_{ I^v \backslash I^w  } 1 \, \d x \leq \abs*{ \rbr*{ - \tfrac{s w_2}{w_1}} - \rbr*{ - \tfrac{\fs v_2}{v_1} } } = \abs*{\tfrac{v_2}{v_1} - \tfrac{w_2}{w_1} }.
\end{equation}
Furthermore,
\nobs that the fact that for all $y \in \R$ it holds that $y \geq - \abs{y}$ implies that for all $v=(v_1, v_2) \in \R^2$ with $\norm{u-v} < \abs{u_1}$ it holds that
\begin{equation} \label{lem:active:intervals:eq:help1}
    u_1 v_1 = (u_1)^2 + ( v_1 - u_1 ) u_1 \geq \abs{u_1}^2 - \abs{u_1 - v_1} \abs{u_1} \geq \abs{u_1}^2 - \norm{u-v} \abs{u_1}  > 0.
\end{equation}
This ensures that for all $v = (v_1, v_2) $, $w = (w_1, w_2 ) \in \R^{2}$ with $\max \cu{\norm{u - v } , \norm{u - w } } < \abs{u_1}$ there exists $\fs \in \cu{-1 , 1 }$ such that $\min \cu{\fs v_1, \fs w_1 } > 0$.
Combining this with \cref{lem:active:intervals:eq:help3} demonstrates for all $v = (v_1, v_2) $, $w = (w_1, w_2 ) \in \R^{2}$ with $\max \cu{\norm{u - v } , \norm{u - w } } \le \frac{\abs{u_1}}{ 2 }$ that
\begin{equation}
\begin{split}
  & \abs*{ \int_{I^v \Delta I^w } \dens ( x ) \, \d x }
  \le M \br*{ \int_{I^v \Delta I^w } 1 \, \d x }
  \le 2 M \abs*{\tfrac{v_2}{v_1} - \tfrac{w_2}{w_1} } = 2 M \abs*{ \frac{v_2 ( w_1 - v_1 ) - v_1 ( w_2 - v_2 ) }{ v_1 w_1 } } \\
  & \leq 2 M \br*{ \abs*{ \frac{v_2 ( w_1 -  v_1 ) }{ v_1 w_1 } } + \abs*{ \frac{ v_1 ( w_2 - v_2 ) }{ v_1 w_1 } } } 
  \leq 2 M  \br*{  \frac{\abs{v_2} \norm{ v -  w } }{ \abs{ v_1 w_1 } } + \frac{ \abs{ v_1 } \norm{ v - w} }{  \abs{ v_1 w_1 } } } \\
& \le \frac{4 M \norm{v} \norm{v - w }}{ \abs{v_1 w_1 } }
\le \br*{ \frac{16 M \norm{v} }{\abs{u_1} ^2 } } \norm{v - w }
\le \br*{ \frac{32 M \norm{u} }{\abs{u_1} ^2 } } \norm{v - w } .
  \end{split}
\end{equation}
This establishes \cref{lem:active:intervals:claim} in the case $u_1 \not= 0$.
\end{cproof}

\cfclear 
\begin{cor} \label{cor:active:intervals}
Assume \cref{setting:snn} and let $\theta \in \fV$. Then there exist
$c, \varepsilon \in (0, \infty)$ such that for all $\vartheta_1, \vartheta_2 \in \R^\fd$ with $\max \cu{ \norm{\vartheta_1 - \theta } , \norm{\vartheta_2 - \theta } } \le \varepsilon$ it holds that
\begin{equation} \cfadd{def:norm}
    \int_{\cup_{i, j = 1 }^\width  ( ( I_i^{\vartheta_1} \cap I_j^{\vartheta_1 } ) \Delta ( I_i^{\vartheta_2 } \cap I_j^{\vartheta_2 } ) ) } \dens ( x ) \, \d x \le \int_{\cup_{i = 1}^\width ( I_i^{\vartheta_1 } \Delta I_i^{\vartheta_2 } ) } \dens ( x ) \, \d x \le c \norm{\vartheta_1 - \vartheta_2 }
\end{equation}
\cfload.
\end{cor}
\begin{cproof}{cor:active:intervals}
\Nobs that \cref{setting:eq:defv} ensures that
$ \min _{k \in  \{1, 2, \ldots, \width \} } \rbr{ \abs{\w{\theta}_k} + \abs{\b{\theta}_k}  } > 0$. Combining this with \cref{lem:active:intervals} shows that there exist $c, \varepsilon \in (0, \infty)$ such that for all $k \in \{1, 2, \ldots, \width \} $,
$\vartheta_1, \vartheta_2 \in \R^\fd$ with $\max \{ \norm{\theta - \vartheta_1}, \norm{\theta - \vartheta_2} \} \leq \varepsilon$ we have that 
\begin{equation} \label{cor:active:intervals:eq1}
    \int_{I_k^{\vartheta_1} \Delta I_k^{\vartheta_2}} \dens ( x ) \, \d x \leq c \norm{\vartheta_1 - \vartheta_2 } .
\end{equation}
Next \nobs that the fact that for all sets $A, \bbA, B, \bbB$ it holds that
\begin{equation}
(A \cap \bbA) \backslash (B \cap \bbB) \subseteq (A \backslash B ) \cup ( \bbA \backslash \bbB) \subseteq (A \Delta B ) \cup ( \bbA \Delta \bbB )
\end{equation}
implies that for all sets $A, \bbA, B, \bbB$ we have that
\begin{equation}
    (A \cap \bbA ) \Delta ( B \cap \bbB ) \subseteq (A \Delta B ) \cup ( \bbA \Delta \bbB ).
\end{equation}
Hence, we obtain for all $\vartheta_1, \vartheta_2 \in \R^\fd$, $i, j \in \cu{1, 2, \ldots, \width }$ that $(I_i^{\vartheta_1} \cap I_j^{\vartheta_1} ) \Delta (I_i^{\vartheta_2} \cap I_j^{\vartheta_2} ) \subseteq (I_i^{\vartheta_1} \Delta I_i^{\vartheta_2} ) \cup (I_j^{\vartheta_2} \Delta I_j^{\vartheta_2} ) $. Combining this with \cref{cor:active:intervals:eq1} proves for all $\vartheta_1, \vartheta_2 \in \R^\fd$ with $\max \{ \norm{\theta - \vartheta_1}, \norm{\theta - \vartheta_2} \} \leq \varepsilon$ that
\begin{equation}
    \begin{split}
          \int_{\cup_{i, j = 1 }^\width  ( ( I_i^{\vartheta_1} \cap I_j^{\vartheta_1 } ) \Delta ( I_i^{\vartheta_2 } \cap I_j^{\vartheta_2 } ) ) } \dens ( x ) \, \d x & \le \int_{\cup_{k = 1}^\width ( I_k^{\vartheta_1 } \Delta I_k^{\vartheta_2 } ) } \dens ( x ) \, \d x  \\
          & \le \sum_{k=1}^\width \br*{ \int_{I_k^{\vartheta_1} \Delta I_k^{\vartheta_2}} \dens ( x ) \, \d x } \le c \width \norm{\vartheta_1 - \vartheta_2 } .
    \end{split}
\end{equation}
\end{cproof}

\begin{cor} \label{cor:derivative:integral:cont}
Assume \cref{setting:snn} and let $i , j \in \{1, 2, \ldots, \width\}$,
 $\phi \in C( \R^\fd \times[a ,b ] , \R )$.
Then
\begin{enumerate} [label = (\roman*)]
    \item \label{cor:derivative:integral:cont:item1} it holds that 
    \begin{equation}
        \fV \ni \theta \mapsto \int_{I_i^\theta} \phi (\theta , x ) \dens ( x ) \, \d x \in \R
    \end{equation}
    is continuous and
    \item \label{cor:derivative:integral:cont:item2} it holds that 
    \begin{equation}
        \fV \ni \theta \mapsto \int_{I_i^\theta \cap I_j^\theta } \phi (\theta , x ) \dens ( x ) \, \d x \in \R
    \end{equation}
    is continuous.
\end{enumerate}
\end{cor}

\cfclear
\begin{cproof} {cor:derivative:integral:cont}
Throughout this proof let $\theta \in \fV$. 
\Nobs that \cref{cor:active:intervals} and \cref{lem:integral:interval:continuous} (applied with $n \with \fd$, $\fu \with a$, $\fv \with b$, $x \with \theta$, $\mu \with (\cB ( [ a , b ] ) \ni A \mapsto \int_A \dens ( x ) \, \d x \in [0, \infty ] )$ in the notation of \cref{lem:integral:interval:continuous}) assure that there exists $\varepsilon \in (0, \infty)$ such that 
\begin{equation}
    \{ \psi \in \R^\fd \colon \norm{\theta - \psi } \leq \varepsilon \} \ni \vartheta \mapsto \int_{I_i^\vartheta} \phi (\vartheta , x ) \dens ( x ) \, \d x \in \R
\end{equation}
and
\begin{equation}
         \{ \psi \in \R^\fd \colon \norm{\theta - \psi } \leq \varepsilon \} \ni \vartheta \mapsto \int_{I_i^\vartheta \cap I_j^\vartheta} \phi (\vartheta , x ) \dens ( x ) \, \d x \in \R
\end{equation}
are continuous.
This shows \cref{cor:derivative:integral:lip:item1,cor:derivative:integral:lip:item2}.
\end{cproof}

\begin{cor} \label{cor:derivative:integral:lip}
Assume \cref{setting:snn},
let $i , j \in \{1, 2, \ldots, \width\}$,
and let $\phi \colon \R^\fd \times[a ,b ] \to \R$ be locally Lipschitz continuous.
Then
\begin{enumerate} [label = (\roman*)]
    \item \label{cor:derivative:integral:lip:item1} it holds that 
    \begin{equation}
        \fV \ni \theta \mapsto \int_{I_i^\theta} \phi (\theta , x ) \dens ( x ) \, \d x \in \R
    \end{equation}
    is locally Lipschitz continuous and
    \item \label{cor:derivative:integral:lip:item2} it holds that 
    \begin{equation}
        \fV \ni \theta \mapsto \int_{I_i^\theta \cap I_j^\theta } \phi (\theta , x ) \dens ( x ) \, \d x \in \R
    \end{equation}
    is locally Lipschitz continuous.
\end{enumerate}
\end{cor}

\cfclear
\begin{cproof} {cor:derivative:integral:lip}
Throughout this proof let $\theta \in \fV$. \Nobs that \cref{cor:active:intervals} and \cref{lem:integral:interval:lipschitz} (applied with $n \with \fd$, $\fu \with a$, $\fv \with b$, $x \with \theta$, $\mu \with (\cB ( [ a , b ] ) \ni A \mapsto \int_A \dens ( x ) \, \d x \in [0, \infty ] )$ in the notation of \cref{lem:integral:interval:lipschitz}) demonstrate that there exist $\varepsilon , \fC \in (0, \infty)$ such that for all $\vartheta_1, \vartheta_2 \in \R^\fd$ with $\max \{ \norm{\theta - \vartheta_1}, \norm{\theta - \vartheta_2} \} \leq \varepsilon$ it holds that
\begin{equation}
    \abs*{\int_{I_i^{\vartheta_1}} \phi ( \vartheta_1 , x ) \dens ( x ) \, \d x - \int_{I_i^{\vartheta_2}} \phi ( \vartheta_2 , x ) \dens ( x ) \, \d x } \leq \fC \norm{\vartheta_1 - \vartheta_2 }
\end{equation}
and
\begin{equation}
        \abs*{\int_{I_i^{\vartheta_1} \cap I_j^{\vartheta_1} } \phi ( \vartheta_1 , x ) \dens ( x ) \, \d x - \int_{I_i^{\vartheta_2} \cap I_j^{\vartheta_2} } \phi ( \vartheta_2 , x ) \dens ( x ) \, \d x } \leq \fC \norm{\vartheta_1 - \vartheta_2 }.
\end{equation}
This establishes \cref{cor:derivative:integral:lip:item1,cor:derivative:integral:lip:item2}.
\end{cproof}

\subsection{Explicit representations for the Hessian matrix of the risk function} \label{subsection:risk:diff}

\cfclear
\begin{prop} \label{prop:loss:differentiable}
 Assume \cref{setting:snn} and let $\theta \in \fV$.
Then
\begin{enumerate} [label=(\roman*)]
    \item \label{prop:loss:diff:item1} it holds that $\cL$ is differentiable at $\theta$ and
    \item \label{prop:loss:diff:item2} it holds that $( \nabla \cL )( \theta ) = \cG ( \theta ) $. 
\end{enumerate} 
\end{prop}
\begin{cproof}{prop:loss:differentiable}
\Nobs that the assumption that $\theta \in \fV$ implies that for all $i \in \{1, 2, \ldots, \width\}$ it holds that $\abs{\w{\theta}_i} + \abs{\b{\theta}_i } > 0$. Hence, we obtain that
\begin{equation} 
  \cL ( \theta )\rbr[\big]{ \smallsum_{i = 1}^\width \abs{\v{\theta}_i } \indicator{ \{ 0 \} } \rbr[\big]{\abs{\w{\theta}_i } + \abs{ \b{\theta}_{i} } } }  = 0.
\end{equation}
Combining this with \cite[Proposition 2.11]{JentzenRiekertFlow} establishes \cref{prop:loss:diff:item1,prop:loss:diff:item2}.
\end{cproof}

\begin{lemma} \label{lem:indicator:terms:cont}
Assume \cref{setting:snn},
let $i \in \cu{1, 2, \ldots, \width}$,
$r, s \in \N_0$,
let $\psi \colon \R \to \R$ satisfy for all $x \in \R \backslash \cu{0}$ that $\psi ( x ) = x^{-1}$,
and let $\bfc \colon ( - \infty , \infty ] \to \R$ satisfy for all $x \in  ( - \infty , \infty ] $ that $\bfc ( x ) = \max \{ \min \{x , b \} , a \}$.
Then 
\begin{enumerate} [ label = (\roman*)]
    \item \label{lem:indicator:terms:cont:item1} it holds for all continuous $\phi \colon \fV \times [a, b] \to \R$ that
    \begin{equation}
        \fV \ni \theta \mapsto \br[\big]{ \psi ( [ \w{\theta}_i]^r \abs{\w{\theta}_i}^s ) } \br[\big]{ \phi ( \theta , \bfc ( \q{\theta}_i ) ) } \indicator{[a,b]} ( \q{\theta}_i ) \in \R 
    \end{equation}
    is continuous and
        \item \label{lem:indicator:terms:cont:item2} it holds for all locally Lipschitz continuous $\phi \colon \fV \times [a, b] \to \R$ that
    \begin{equation}
        \fV \ni \theta \mapsto \br[\big]{ \psi ( [ \w{\theta}_i]^r \abs{\w{\theta}_i}^s ) } \br[\big]{ \phi ( \theta , \bfc ( \q{\theta}_i ) ) } \indicator{[a,b]} ( \q{\theta}_i ) \in \R 
    \end{equation}
    is locally Lipschitz continuous.
\end{enumerate}
\end{lemma}
\begin{cproof}{lem:indicator:terms:cont}
\Nobs that \cref{setting:eq:defv} shows for all $\theta \in \fV$ that $\abs{\w{\theta}_i} + \abs{\b{\theta}_i } > 0$.
Hence, we obtain for all $\theta \in \fV$ with $\w{\theta}_i = 0$ that $\b{\theta}_i \not= 0$.
This implies that for all $\theta \in \fV$ with $\w{\theta}_i = 0$ there exists $\varepsilon \in (0, \infty)$ such that for all $\vartheta \in \cu{\psi \in \R^\fd \colon \norm{\psi - \theta } < \varepsilon }$ it holds that $\q{\vartheta}_i \notin [a,b]$.
Combining this with \cref{setting:eq:defqi} and the fact that for all $\theta \in \fV$ it holds that $\q{\theta}_i \notin \cu{a,b}$ establishes \cref{lem:indicator:terms:cont:item1,lem:indicator:terms:cont:item2}.
\end{cproof}

\cfclear
\begin{lemma} \label{lem:leibniz:rule}
Let $a \in \R$,
$b \in (a, \infty)$,
let $U \subseteq \R$ be open,
let $\phi = (\phi_x(t))_{(x,t) \in [a,b] \times U } \in C([a,b] \times U , \R)$ satisfy for all $x \in [a,b]$ that $\phi_x \in C^1 ( U , \R)$,
assume that $[a,b] \times U \ni (x,t) \mapsto ( \phi_x ) ' (t) \in \R$ is continuous,
let $\psi_0, \psi_1 \in C^1( U , [a,b])$,
and let $\Phi \colon U \to \R$ satisfy for all $t \in U$ that
\begin{equation} \label{lem:leibniz:eq:defphi}
    \Phi(t) = \int_{\psi_0(t)}^{\psi_1 ( t ) } \phi_x ( t ) \, \d x.
\end{equation}
Then
\begin{enumerate} [label = (\roman*)]
    \item it holds that $\Phi \in C^1 ( U , \R)$ and
    \item it holds for all $t \in U$ that
    \begin{equation}
        \Phi' ( t ) = \br*{ \phi_{\psi_1(t)} ( t )} \br*{ ( \psi_1 ) ' ( t ) } - \br*{ \phi_{\psi_0 ( t ) } ( t )} \br*{ ( \psi_0 ) ' ( t ) } + \int_{\psi_0 ( t ) }^{\psi_1 ( t ) } ( \phi_x ) ' ( t ) \, \d x .
    \end{equation}
\end{enumerate}
\end{lemma}
\begin{cproof}{lem:leibniz:rule}
Throughout this proof let $\Psi \colon [a,b] \times U \to \R$ satisfy for all $x \in [a,b]$, $t \in U$ that
\begin{equation} \label{lem:leibniz:eq:defpsi}
    \Psi ( x , t ) = \int_a^x \phi_y(t) \, \d y.
\end{equation}
\Nobs that \cref{lem:leibniz:eq:defphi,lem:leibniz:eq:defpsi} imply for all $t \in U$ that
\begin{equation} \label{lem:leibniz:eq:phipsi}
    \Phi ( t ) = \int_a^{\psi_1 ( t ) } \phi_x ( t ) \, \d x - \int_a^{\psi_0 ( t ) } \phi_x ( t ) \, \d x = \Psi ( \psi_1 ( t ) , t ) - \Psi ( \psi_0 ( t ) , t ).
\end{equation}
Next \nobs that the fundamental theorem of calculus ensures for all $x \in [a,b]$, $t \in U$ that $\frac{\partial}{\partial x} \Psi ( x , t ) = \phi_x ( t ) $.
In addition, \nobs that \cref{lem:interchange} assures for all $x \in [a,b]$, $t \in U$ that $\frac{\partial}{\partial t} \Psi ( x , t ) = \int_a^x ( \phi_y ) ' ( t ) \, \d y$.
Furthermore, \nobs that the assumption that $[a,b] \times U \ni (x,t) \mapsto \phi_x (t) \in \R$ is continuous,
the assumption that $[a,b] \times U \ni (x,t) \mapsto ( \phi_x ) ' ( t ) \in \R$ is continuous,
and the dominated convergence theorem demonstrate that $[a,b] \times U \ni ( x , t ) \mapsto \frac{\partial}{\partial x} \Psi ( x , t ) \in \R$ and $[a,b] \times U \ni ( x , t ) \mapsto \frac{\partial}{\partial t} \Psi ( x , t ) \in \R$
are continuous.
Hence, we obtain that $\Psi \in C^1 ( [a, b ] \times U , \R)$. Combining this with \cref{lem:leibniz:eq:phipsi} and the chain rule shows for all $t \in U$ that $\Phi \in C^1 ( U , \R)$ and
\begin{equation}
    \begin{split}
        \Phi' ( t ) &= ( \psi_1 ) ' ( t ) \rbr[\big]{\tfrac{\partial}{\partial x} \Psi } ( \psi_1  ( t ) , t ) + \rbr[\big]{\tfrac{\partial}{\partial t} \Psi } ( \psi_1  ( t ) , t ) \\
        & \quad - ( \psi_0 ) ' ( t ) \rbr[\big]{\tfrac{\partial}{\partial x} \Psi } ( \psi_0  ( t ) , t ) - \rbr[\big]{\tfrac{\partial}{\partial t} \Psi } ( \psi_0  ( t ) , t ) \\
        &= \br*{ ( \psi_1 ) ' ( t ) } \br*{ \phi_{\psi_1 ( t ) } ( t ) } + \int_a^{\psi_1 ( t )} ( \phi_x ) ' ( t ) \, \d x  - \br*{ ( \psi_0 )  ' ( t ) } \br*{ \phi_{\psi_0 ( t ) } ( t ) } - \int_a^{\psi_0 ( t )} ( \phi_x ) ' ( t ) \, \d x \\
        &= \br*{ ( \psi_1 ) ' ( t ) } \br*{ \phi_{\psi_1 ( t ) } ( t ) } - \br*{ ( \psi_0 ) ' ( t ) } \br*{ \phi_{\psi_0 ( t ) } ( t )  } + \int_{\psi_0 ( t ) }^{\psi_1 ( t )} ( \phi_x ) ' ( t ) \, \d x .
    \end{split}
\end{equation}
\end{cproof}

\cfclear
\begin{lemma} \label{lem:twice:diff}
Assume \cref{setting:snn},
let $\psi \colon \R \to \R$ satisfy for all $x \in \R \backslash \{ 0 \}$ that $\psi ( x ) = x^{-1}$,
and let $\bfc \colon ( - \infty , \infty ] \to \R$ satisfy for all $x \in  ( - \infty , \infty ] $ that $\bfc ( x ) = \max \{ \min \{x , b \} , a \}$. Then 
\begin{enumerate} [ label = (\roman*)]
    \item \label{lem:twice:diff:item0} it holds that $\fV \subseteq \R^\fd$ is open,
    \item \label{lem:twice:diff:item1} it holds that $\cL |_\fV \in C^2( \fV , \R)$, and
    \item \label{lem:twice:diff:item2} it holds for all $\theta = (\theta_1, \ldots, \theta_\fd ) \in \fV$, $i,j \in \{1, 2, \ldots, \width \}$ that
    \begin{equation} \label{lem:twice:diff:eq4}
         \rbr[\big]{\tfrac{\partial^2}{ \partial \theta_j \partial \theta_{\fd } } \cL } ( \theta) =   2 \v{\theta}_j \tint_{I_j^\theta} x \dens ( x ) \, \d x ,
        \end{equation}
    \begin{equation} \label{lem:twice:diff:eq5}
           \rbr[\big]{\tfrac{\partial^2}{ \partial \theta_{\width + j} \partial \theta_{\fd }} \cL } ( \theta) =   2 \v{\theta}_j \tint_{I_j^\theta}  \dens ( x ) \, \d x , 
        \end{equation}
  \begin{equation} \label{lem:twice:diff:eq2}
          \rbr[\big]{\tfrac{\partial^2}{ \partial \theta_{2 \width + j} \partial \theta_{\fd } } \cL } ( \theta) = 2 \tint_a^b \br[\big]{\Rect ( \w{\theta}_j x+ \b{\theta}_j) } \dens ( x ) \, \d x , 
        \end{equation}
      \begin{equation} \label{lem:twice:diff:eq1}
          \rbr[\big]{\tfrac{\partial^2}{\partial \theta_{\fd }^2} \cL } ( \theta) = 2 \tint_a^b \dens ( x ) \, \d x, 
    \end{equation}
    \begin{multline} \label{lem:twice:diff:eq6}
           \rbr[\big]{\tfrac{\partial^2}{ \partial \theta_j \partial \theta_{2\width + i} } \cL } ( \theta) 
           = 2 \v{\theta}_j \tint_{I_j^\theta} x \br[\big]{ \Rect( \w{\theta}_i x + \b{\theta}_i ) } \dens ( x ) \, \d x \\
           + 2 \indicator{ \{ i \} } ( j ) \tint_{I_i^\theta} x ( \realization{\theta} ( x ) - f ( x ) ) \dens ( x ) \, \d x, 
        \end{multline}
    \begin{multline} \label{lem:twice:diff:eq7}
        \rbr[\big]{\tfrac{\partial^2}{\partial \theta_{\width + j}\partial \theta_{2 \width + i} } \cL } ( \theta) 
        = 2 \v{\theta}_j \tint_{I_j^\theta}  \br[\big]{\Rect( \w{\theta}_i x + \b{\theta}_i ) } \dens ( x ) \, \d x \\
        + 2 \indicator{ \{ i \} } ( j ) \tint_{I_i^\theta}  ( \realization{\theta} ( x ) - f ( x ) ) \dens ( x ) \, \d x,
\end{multline}
     \begin{equation} \label{lem:twice:diff:eq3}
            \rbr[\big]{\tfrac{\partial^2}{\partial \theta_{2 \width + j } \partial \theta_{2 \width + i } } \cL } ( \theta) = 2 \tint_a^b \br[\big]{ \Rect ( \w{\theta}_i x+ \b{\theta}_i) } \br[\big]{ \Rect (  \w{\theta}_j x+ \b{\theta}_j) } \dens ( x ) \, \d x, 
        \end{equation}
\begin{multline} \label{lem:twice:diff:eq8}
        \rbr[\big]{\tfrac{\partial^2}{ \partial \theta_j \partial \theta_i } \cL } ( \theta) = 2 \v{\theta}_i \v{\theta}_j \tint_{I_i^\theta \cap I_j^\theta} x^2 \dens ( x ) \, \d x \\
         -  2 \v{\theta}_i \b{\theta}_i \indicator{\{ i \}} ( j ) \indicator{[a , b]}( \q{\theta}_i ) [\psi ( \w{\theta}_i \abs{\w{\theta}_i } )] [ \bfc ( \q{\theta}_i ) ] (\realization{\theta} ( \bfc (\q{\theta}_i) ) - f ( \bfc ( \q{\theta}_i ) ) ) \dens ( \bfc ( \q{\theta}_i ) ), 
            \end{multline}
    \begin{multline} \label{lem:twice:diff:eq9}
\rbr[\big]{\tfrac{\partial^2}{ \partial \theta_j \partial \theta_{\width + i } } \cL } ( \theta) = 2 \v{\theta}_i \v{\theta}_j \tint_{I_i^\theta \cap I_j^\theta} x \dens ( x ) \, \d x \\
 + 2 \v{\theta}_i \indicator{\{ i \}} ( j ) \indicator{[a , b]}( \q{\theta}_i ) [ \psi ( \abs{ \w{\theta}_i } ) ] [ \bfc( \q{\theta}_i ) ] (\realization{\theta} ( \bfc ( \q{\theta}_i ) ) - f ( \bfc ( \q{\theta}_i ) ) ) \dens (\bfc ( \q{\theta}_i ) ), 
         \end{multline}
         and
    \begin{multline} \label{lem:twice:diff:eq10}
\rbr[\big]{\tfrac{\partial^2}{ \partial \theta_{\width + j } \partial \theta_{\width + i } } \cL } ( \theta) = 2 \v{\theta}_i \v{\theta}_j \tint_{I_i^\theta \cap I_j^\theta }  \dens ( x ) \, \d x \\
 +  2 \v{\theta}_i \indicator{\{ i \}} ( j ) \indicator{[a , b]}( \q{\theta}_i ) [ \psi ( \abs{ \w{\theta}_i } ) ] (\realization{\theta} ( \bfc ( \q{\theta}_i ) ) - f ( \bfc ( \q{\theta}_i ) ) ) \dens (\bfc ( \q{\theta}_i ) ).
    \end{multline}
\end{enumerate} 
\end{lemma}

\cfclear
\begin{cproof} {lem:twice:diff}
\Nobs that \cref{setting:eq:defv} establishes \cref{lem:twice:diff:item0}. Next \nobs that \cref{prop:loss:differentiable} ensures that $\fV \ni \theta \mapsto \cL ( \theta ) \in \R$ is differentiable and satisfies $\nabla (\cL |_\fV) = \cG_\fV$. In addition, \nobs that \cref{eq:loss:gradient} and \cref{cor:derivative:integral:cont} prove that $\cG |_\fV$ is continuous. Hence, we obtain that $\cL | _\fV \in C^1 ( \fV , \R)$ and \begin{equation} \label{lem:twice:diff:eq:c1}
    \nabla (\cL |_\fV) = \cG | _\fV.
\end{equation}
Combining this with \cref{eq:loss:gradient}, \cref{cor:interchange}, and the product rule establishes \cref{lem:twice:diff:eq1,lem:twice:diff:eq2,lem:twice:diff:eq3,lem:twice:diff:eq4,lem:twice:diff:eq5,lem:twice:diff:eq6,lem:twice:diff:eq7}. 
In the next step we prove \cref{lem:twice:diff:eq8,lem:twice:diff:eq9,lem:twice:diff:eq10}
and for this let $\theta = (\theta_1, \ldots, \theta_\fd ) \in \fV$, $i, j \in \cu{1, 2, \ldots, \width }$. In our proof of \cref{lem:twice:diff:eq8,lem:twice:diff:eq9,lem:twice:diff:eq10} we distinguish between the case $( i \not= j ) $, 
the case $( ( i = j ) \wedge ( \max \{ \w{\theta}_i a + \b{\theta}_i , \w{\theta}_i b + \b{\theta}_i \} < 0 ) )$,
the case $( ( i = j ) \wedge ( \min \{ \w{\theta}_i a + \b{\theta}_i , \w{\theta}_i b + \b{\theta}_i \} > 0 ) )$,
the case $( ( i = j ) \wedge ( \w{\theta}_i a + \b{\theta}_i < 0 < \w{\theta}_i b + \b{\theta}_i ) )$,
and the case $( ( i = j ) \wedge ( \w{\theta}_i a + \b{\theta}_i > 0 > \w{\theta}_i b + \b{\theta}_i ) ) $.
We first establish \cref{lem:twice:diff:eq8,lem:twice:diff:eq9,lem:twice:diff:eq10} in the case $( i \not= j ) $.
\Nobs that for all $k \in \cu{0,1}$ and almost all $x \in [a,b]$ it holds that
\begin{equation}
    \tfrac{\partial}{\partial \theta_{k \width + j} } \realization{\theta} ( x ) =   \tfrac{\partial}{\partial \theta_{k \width + j} }  \rbr[\big]{ \theta_{2 \width + j }  [\Rect ( \theta_j x + \theta_{\width + j } ) ] } =\v{\theta}_j x^{ 1 - k } \indicator{I_j^\theta} ( x ).
\end{equation}
Combining this with \cref{eq:loss:gradient},
\cref{lem:twice:diff:eq:c1},
and \cref{cor:interchange} (applied for every $k, \ell \in \cu{0, 1 }$ with $n \with \fd$, $j \with k \width + j$, $\phi \with ( \R^\fd \times [a,b] \ni (\vartheta , x ) \mapsto x ^{1 - \ell } ( \realization{\vartheta} ( x ) - f ( x ) ) \dens ( x ) \indicator{I_i^\vartheta} ( x ) \in \R ) $ in the notation of \cref{cor:interchange}) demonstrates for all $k, \ell \in \cu{0,1}$ that
\begin{equation}
\begin{split}
     &\rbr[\big]{\tfrac{\partial^2}{ \partial \theta_{k \width + j} \partial \theta_{\ell \width + i } } \cL } ( \theta) 
     =  \rbr[\big]{\tfrac{\partial}{ \partial \theta_{k \width + j } } \cG_{\ell \width + i } } ( \theta) \\
     &= \tfrac{\partial}{\partial \theta_{k \width + j } } \rbr*{ 2 \v{\theta}_i \int_a^b x ^{ 1 - \ell} ( \realization{\theta} ( x ) - f ( x ) ) \dens ( x ) \indicator{I_i^\theta} ( x ) \, \d x } 
     = 2 \v{\theta}_i \v{\theta}_j \int_{I_i^\theta \cap I_j^\theta} x^{2 - k - \ell} \dens ( x ) \, \d x.
\end{split}
\end{equation}
This establishes \cref{lem:twice:diff:eq8,lem:twice:diff:eq9,lem:twice:diff:eq10} in the case $(i \not= j )$.

We next prove \cref{lem:twice:diff:eq8,lem:twice:diff:eq9,lem:twice:diff:eq10} in the case \begin{equation} \label{lem:diff:eq:case1}
    ( i = j ) \wedge (\max \{ \w{\theta}_i a + \b{\theta}_i , \w{\theta}_i b + \b{\theta}_i \} < 0 ). 
\end{equation}
\Nobs that \cref{lem:diff:eq:case1} implies that there exists $\delta \in (0, \infty)$ such that for all $h \in \R^\fd$ with $\norm{h} < \delta$ it holds that $\q{\theta+h}_i \notin [a , b]$ and $I_i^{\theta + h } = \emptyset$ \cfadd{def:norm}\cfload. Combining this with \cref{eq:loss:gradient} and
\cref{lem:twice:diff:eq:c1} ensures that $\rbr[\big]{\tfrac{\partial^2}{ \partial \theta_i ^2} \cL } ( \theta) =  \rbr[\big]{\tfrac{\partial^2}{\partial \theta_i \partial \theta_{\width + i}} \cL } ( \theta) = \rbr[\big]{\tfrac{\partial^2}{ \partial \theta_{\width + i} ^2} \cL } ( \theta) = 0$, as desired. 

In the next step we prove \cref{lem:twice:diff:eq8,lem:twice:diff:eq9,lem:twice:diff:eq10} in the case 
\begin{equation} \label{lem:diff:eq:case2}
    ( i = j ) \wedge (\min \{ \w{\theta}_i a + \b{\theta}_i , \w{\theta}_i b + \b{\theta}_i \} > 0 ).
\end{equation}
\Nobs that \cref{lem:diff:eq:case2} implies that there exists $\delta \in (0, \infty)$ such that for all $h \in \R^\fd$ with $\norm{h} < \delta$ it holds that $\q{\theta+h}_i \notin [a , b]$ and $I_i^{\theta + h } = [a , b]$.
Combining \cref{eq:loss:gradient},
\cref{lem:twice:diff:eq:c1},
and \cref{cor:interchange} hence shows that $\rbr[\big]{\tfrac{\partial^2}{ \partial \theta_i ^2} \cL } ( \theta) = 2 (\v{\theta}_i)^2 \int_a^b x^2 \dens ( x ) \, \d x $,
$\rbr[\big]{\tfrac{\partial^2}{ \partial \theta_i \partial \theta_{\width + i}} \cL } ( \theta) = 2 (\v{\theta}_i)^2 \int_a^b x \dens ( x ) \, \d x $, and  $\rbr[\big]{\tfrac{\partial^2}{ \partial \theta_{\width + i} ^2} \cL } ( \theta) = 2 (\v{\theta}_i)^2 \int_a^b \dens ( x ) \, \d x $, as claimed.

In the remaining cases we employ \cref{lem:leibniz:rule} since the interval $I_i^\theta$ depends on $\w{\theta}_i$ and $\b{\theta}_i$ in these cases.
We first consider the case 
\begin{equation} \label{lem:diff:eq:case3}
    (i=j) \wedge ( \w{\theta}_i a + \b{\theta}_i < 0 < \w{\theta}_i b + \b{\theta}_i ) .
\end{equation}
\Nobs that \cref{lem:diff:eq:case3} ensures that there exists an open neighborhood $U \subseteq \R^\fd$ of $\theta$ which satisfies for all $\vartheta \in U$ that 
$\w{\vartheta }_i>0$, $\q{\vartheta }_i \in (a , b)$, and $I_i^{\vartheta } = ( \q{\vartheta }_i, b ]$. 
Furthermore, \nobs that $U \ni \vartheta \mapsto \q{\vartheta}_i = - \frac{\b{\vartheta}_i}{\w{\vartheta}_i } \in \R$ is continuously differentiable and satisfies $\frac{\partial}{\partial \theta_i} \q{\theta}_i = \frac{\b{\theta}_i}{( \w{\theta}_i ) ^2  } = - \frac{\q{\theta}_i}{\w{\theta}_i}$ and $\frac{\partial}{\partial \theta_{\width + i }} \q{\theta}_i = - \frac{1}{\w{\theta}_i }$. 
Combining \cref{lem:leibniz:rule}, \cref{eq:loss:gradient}, and \cref{lem:twice:diff:eq:c1} hence shows that
\begin{equation}
    \begin{split}
        \rbr[\big]{\tfrac{\partial^2}{ \partial \theta_i ^2} \cL } ( \theta) &= 2 (\v{\theta}_i )^2 \tint_{I_i^\theta} x^2 \dens ( x ) \, \d x - \br*{ \frac{2 \v{\theta}_i \b{\theta}_i}{(\w{\theta}_i)^2} } \q{\theta}_i (\realization{\theta} ( \q{\theta}_i) - f ( \q{\theta}_i) ) \dens ( \q{\theta}_i ), \\
\rbr[\big]{\tfrac{\partial^2}{ \partial \theta_i \partial \theta_{\width + i}} \cL } ( \theta) &= 2( \v{\theta}_i) ^2 \tint_{I_i^\theta} x \dens ( x ) \, \d x + \br*{ \frac{2 \v{\theta}_i }{ \w{\theta}_i} } \q{\theta}_i (\realization{\theta} ( \q{\theta}_i) - f ( \q{\theta}_i ) ) \dens ( \q{\theta}_i ), \\ 
\text{and} \qquad \rbr[\big]{\tfrac{\partial^2}{ \partial \theta_{\width + i} ^2} \cL } ( \theta) &= 2 (\v{\theta}_i )^2 \tint_{I_i^\theta} \dens ( x ) \, \d x + \br*{ \frac{2 \v{\theta}_i}{ \w{\theta}_i} } (\realization{\theta} ( \q{\theta}_i) - f ( \q{\theta}_i ) ) \dens ( \q{\theta}_i ).
    \end{split}
\end{equation}
This establishes \cref{lem:twice:diff:eq8,lem:twice:diff:eq9,lem:twice:diff:eq10} in the case $((i=j) \wedge ( \w{\theta}_i a + \b{\theta}_i < 0 < \w{\theta}_i b + \b{\theta}_i ) )$.
It remains to consider the case 
\begin{equation} \label{lem:diff:eq:case4}
     ( i = j ) \wedge ( \w{\theta}_i a + \b{\theta}_i > 0 > \w{\theta}_i b + \b{\theta}_i ) 
\end{equation}
\Nobs that \cref{lem:diff:eq:case4} assures that $\w{\theta}_i < 0$, $\q{\theta}_i \in (a , b)$, and $I_i^\theta = [a, \q{\theta}_i)$.
Combining \cref{lem:leibniz:rule}, \cref{eq:loss:gradient}, and \cref{lem:twice:diff:eq:c1} therefore demonstrates that
\begin{equation}
    \begin{split}
        \rbr[\big]{\tfrac{\partial^2}{ \partial \theta_i ^2} \cL } ( \theta) &= 2 (\v{\theta}_i )^2 \tint_{I_i^\theta} x^2 \dens ( x ) \, \d x + \br*{ \frac{2 \v{\theta}_i \b{\theta}_i}{(\w{\theta}_i)^2} } \q{\theta}_i (\realization{\theta} ( \q{\theta}_i) - f ( \q{\theta}_i) ) \dens ( \q{\theta}_i ), \\
\rbr[\big]{\tfrac{\partial^2}{ \partial \theta_i \partial \theta_{\width + i}} \cL } ( \theta) &= 2( \v{\theta}_i) ^2 \tint_{I_i^\theta} x \dens ( x ) \, \d x - \br*{ \frac{2 \v{\theta}_i }{ \w{\theta}_i} } \q{\theta}_i (\realization{\theta} ( \q{\theta}_i) - f ( \q{\theta}_i ) ) \dens ( \q{\theta}_i ), \\ 
\text{and} \qquad \rbr[\big]{\tfrac{\partial^2}{ \partial \theta_{\width + i} ^2} \cL } ( \theta) &= 2 (\v{\theta}_i )^2 \tint_{I_i^\theta} \dens ( x ) \, \d x - \br*{ \frac{2 \v{\theta}_i}{ \w{\theta}_i} } (\realization{\theta} ( \q{\theta}_i) - f ( \q{\theta}_i ) ) \dens ( \q{\theta}_i ).
    \end{split}
\end{equation}
This establishes \cref{lem:twice:diff:eq8,lem:twice:diff:eq9,lem:twice:diff:eq10} in the case $((i=j) \wedge ( \w{\theta}_i a + \b{\theta}_i > 0 > \w{\theta}_i b + \b{\theta}_i )$.

Finally, \nobs that \cref{cor:derivative:integral:cont} and \cref{lem:indicator:terms:cont:item1} in \cref{lem:indicator:terms:cont} imply that the partial derivatives in \cref{lem:twice:diff:eq1,lem:twice:diff:eq2,lem:twice:diff:eq3,lem:twice:diff:eq4,lem:twice:diff:eq5,lem:twice:diff:eq6,lem:twice:diff:eq7,lem:twice:diff:eq8,lem:twice:diff:eq9,lem:twice:diff:eq10} are continuous on $\fV$.
\end{cproof}

\begin{lemma} \label{lem:twice:diff:lipschitz}
Assume \cref{setting:snn} and assume that $f$ is Lipschitz continuous. Then
\begin{enumerate} [label = (\roman*)]
\item \label{lem:twice:diff:lip:item0} it holds that $\fV \subseteq \R^\fd$ is open,
    \item \label{lem:twice:diff:lip:item1} it holds that $\cL |_\fV \in C^2(\fV , \R)$, and
    \item \label{lem:twice:diff:lip:item2} it holds that $\fV \ni \theta \mapsto (\Hs  \cL) ( \theta ) \in \R^{\fd \times \fd }$ is locally Lipschitz continuous.
\end{enumerate}
\end{lemma}
\begin{cproof}{lem:twice:diff:lipschitz}
\Nobs that \cref{lem:twice:diff} establishes \cref{lem:twice:diff:lip:item0,lem:twice:diff:lip:item1}.
Moreover, \nobs that \cref{lem:twice:diff}, \cref{cor:derivative:integral:lip},
\cref{lem:indicator:terms:cont:item2} in \cref{lem:indicator:terms:cont},
the assumption that $f$ is Lipschitz continuous,
and the assumption that $\dens$ is Lipschitz continuous establish \cref{lem:twice:diff:lip:item2}.
\end{cproof}

\cfclear 
\begin{cor} \label{cor:second:derivatives:globalmin} 
Assume \cref{setting:snn}, let $\theta \in \fV$, $i,j \in \{1, 2, \ldots, \width\}$,
and assume for all $x \in [a , b]$ that $ \realization{\theta} ( x ) = f ( x )$. Then
\begin{equation} \label{cor:second:derivatives:globalmin:eq}
    \begin{split}
          \rbr[\big]{\tfrac{\partial^2}{\partial \theta_i \partial \theta_j} \cL } ( \theta) &= 2 \v{\theta}_i \v{\theta}_j \int_{I_i^\theta \cap I_j^\theta} x^2 \dens ( x ) \, \d x, \\
          \rbr[\big]{\tfrac{\partial^2}{\partial \theta_i \partial \theta_{\width + j}} \cL } ( \theta) &= 2 \v{\theta}_i \v{\theta}_j \int_{I_i^\theta \cap I_j^\theta} x \dens ( x ) \, \d x, \\
        \text{and} \qquad  \rbr[\big]{\tfrac{\partial^2}{\partial \theta_{\width + i} \partial \theta_{\width + j}} \cL } ( \theta) &= 2 \v{\theta}_i \v{\theta}_j \int_{I_i^\theta \cap I_j^\theta} \dens ( x ) \, \d x.
    \end{split}
\end{equation}
\end{cor}
\begin{cproof}{cor:second:derivatives:globalmin}
\Nobs that the assumption that for all $x \in [a , b]$ it holds that $ \realization{\theta} ( x ) = f ( x )$ and \cref{lem:twice:diff} establish \cref{cor:second:derivatives:globalmin:eq}.
\end{cproof}

\subsection{Upper bounds for the entries of the Hessian matrix of the risk function} \label{subsection:hess:bound}

\begin{lemma} \label{lem:hessian:max:upperbound}
Assume \cref{setting:snn},
let $\tbound \in [1, \infty)$, $A \in \R$ satisfy $A = \max \{ 1 , \abs{a}, \abs{b} , b - a \}$,
and let $\theta \in \fV$ satisfy $\max_{i \in \{1, 2, \ldots, \fd \} } \abs{\theta_i} \leq \tbound$ and $\min_{j \in \{1, 2, \ldots, \width \} } \rbr[\big]{ ( \w{\theta}_j - \frac{1}{2} ) \indicator{[a,b]} ( \q{\theta}_j ) }  \geq 0$. Then 
\begin{equation} \label{lem:hessian:upperbound:eq}
\begin{split}
   & \max\nolimits_{i,j \in \{1, 2, \ldots, \fd \} } \abs[\big]{ \rbr[\big]{ \tfrac{\partial^2}{\partial \theta_i \partial \theta_j} \cL } ( \theta )} \\
   &\leq \rbr[\big]{ 8 A^3 \tbound^2 + 8 A ^2 \tbound ^2 \br[\big]{  \sup\nolimits_{x \in [a,b]} \abs{\realization{\theta} ( x ) - f ( x ) } } } \rbr[\big]{ \sup\nolimits_{x \in [a,b]} \dens ( x ) }.
    \end{split}
\end{equation}
\end{lemma}
\begin{cproof}{lem:hessian:max:upperbound}
Throughout this proof let $\psi \colon \R \to \R$ satisfy for all $x \in \R \backslash \cu{0} $ that $\psi ( x ) = x^{-1}$ and let $\bfc \colon ( - \infty , \infty ] \to \R$ satisfy for all $x \in ( - \infty , \infty ]$ that $\bfc ( x ) = \max \{ \min \{x , b \} , a \}$.
\Nobs that \cref{lem:twice:diff} implies for all $i,j \in \{1, 2, \ldots, \width \}$ that 
\begin{equation}
    \abs[\big]{ \rbr[\big]{\tfrac{\partial^2}{\partial \theta_{\fd }^2} \cL } ( \theta)} = 2 \abs[\big]{\tint_a^b \dens ( x ) \, \d x} \leq 2 A \rbr[\big]{ \sup\nolimits_{x \in [a,b]} \dens ( x ) },
\end{equation}
    \begin{equation}
    \begin{split}
         \abs[\big]{ \rbr[\big]{\tfrac{\partial^2}{\partial \theta_{2 \width + j} \partial \theta_{\fd }} \cL } ( \theta) } &= 2 \abs[\big]{ \tint_a^b \br{\Rect ( \w{\theta}_j x+ \b{\theta}_j ) } \dens ( x ) \, \d x } \leq 2 \tint_a^b \abs{\Rect( \w{\theta}_j x + \b{\theta}_j ) } \dens ( x ) \, \d x \\
         &\leq 2 A ( \abs{\w{\theta}_j } + \abs{\b{\theta}_j } ) \tint_a^b \dens ( x ) \, \d x \leq 4 A^2 \tbound \rbr[\big]{ \sup\nolimits_{x \in [a,b]} \dens ( x ) }, 
         \end{split}
    \end{equation}
    \begin{equation}
        \begin{split} 
              &\abs[\big]{ \rbr[\big]{\tfrac{\partial^2}{\partial \theta_{2 \width + i} \partial \theta_{2 \width + j}} \cL } ( \theta) } = 2 \abs[\big]{ \tint_a^b \br{ \Rect ( \w{\theta}_i x+ \b{\theta}_i ) } \br{ \Rect (  \w{\theta}_j x+ \b{\theta}_j ) } \dens ( x ) \, \d x } \\
              & \leq  2 \tint_a^b \abs{ \Rect ( \w{\theta}_i x+ \b{\theta}_i) \Rect (  \w{\theta}_j x+ \b{\theta}_j) } \dens ( x ) \, \d x \\
              &\leq 2 A^2 ( \abs{\w{\theta}_i } + \abs{\b{\theta}_i } ) ( \abs{\w{\theta}_j } + \abs{\b{\theta}_j } ) \tint_a^b \dens ( x ) \, \d x \leq 8 A^3 \tbound^2 \rbr[\big]{ \sup\nolimits_{x \in [a,b]} \dens ( x ) }, 
        \end{split}
    \end{equation}
    \begin{equation} 
       \abs[\big]{  \rbr[\big]{\tfrac{\partial^2}{\partial \theta_{\fd } \partial \theta_j} \cL } ( \theta) }=   2 \abs{ \v{\theta}_j} \abs[\big]{ \tint_{I_j^\theta} x \dens ( x ) \, \d x } \leq 2A^2 \tbound \rbr[\big]{ \sup\nolimits_{x \in [a,b]} \dens ( x ) } , 
       \end{equation}
       \begin{equation} 
           \abs[\big]{\rbr[\big]{\tfrac{\partial^2}{\partial \theta_{\fd } \partial \theta_{\width + j}} \cL } ( \theta) } =   2 \abs{\v{\theta}_j}  \tint_{I_j^\theta} \dens ( x ) \, \d x \leq 2 A \tbound \rbr[\big]{ \sup\nolimits_{x \in [a,b]} \dens ( x ) } ,
          \end{equation}
          \begin{equation}
              \begin{split} 
         \abs[\big]{  \rbr[\big]{\tfrac{\partial^2}{\partial \theta_{2\width + i} \partial \theta_j} \cL } ( \theta) } &\leq 2 \abs{ \v{\theta}_j } \abs[\big]{ \tint_{I_j^\theta} x \br{ \Rect( \w{\theta}_i x + \b{\theta}_i ) } \dens ( x ) \, \d x } + 2 \abs[\big]{\tint_{I_i^\theta} x ( \realization{\theta} ( x ) - f ( x ) ) \dens ( x ) \, \d x } \\
         &\leq \rbr[\big]{ 4 A^3 \tbound ^2 + 2 A^2 \br[\big]{ \sup\nolimits_{x \in [a,b]} \abs{\realization{\theta} ( x ) - f ( x ) } } } \rbr[\big]{ \sup\nolimits_{x \in [a,b]} \dens ( x ) }, \\
         \end{split}
         \end{equation}
         and
         \begin{equation} 
         \begin{split} 
        \abs[\big]{ \rbr[\big]{\tfrac{\partial^2}{\partial \theta_{2 \width + i} \partial \theta_{\width + j}} \cL } ( \theta) } &\leq 2 \abs{ \v{\theta}_j } \abs[\big]{ \tint_{I_j^\theta} \br{ \Rect( \w{\theta}_i x + \b{\theta}_i ) } \dens ( x ) \, \d x } + 2 \abs[\big]{\tint_{I_i^\theta}  ( \realization{\theta} ( x ) - f ( x ) ) \dens ( x ) \, \d x } \\
        &\leq  \rbr[\big]{ 4 A^2 \tbound ^2 + 2 A \br[\big]{ \sup\nolimits_{x \in [a,b]} \abs{\realization{\theta} ( x ) - f ( x ) } } } \rbr[\big]{ \sup\nolimits_{x \in [a,b]} \dens ( x ) }.
    \end{split}
\end{equation}
In addition, \nobs that \cref{lem:twice:diff} and the fact that for all $i \in \{1, 2, \ldots, \width\}$ with $\q{\theta}_i \in [a,b]$ it holds that $\w{\theta}_i \geq \frac{1}{2}$ show that for all $i,j \in \{1, 2, \ldots, \width \}$ it holds that 
\begin{equation}
    \begin{split}
         \abs[\big]{ \rbr[\big]{\tfrac{\partial^2}{ \partial \theta_i \partial \theta_j} \cL } ( \theta) } 
         &\leq 2 \abs{\v{\theta}_i \v{\theta}_j} \abs[\big]{ \tint_{I_i^\theta \cap I_j^\theta} x^2 \dens ( x ) \, \d x } \\
         & \quad + \indicator{[a,b]} ( \q{\theta}_i ) \abs[\big]{2 \v{\theta}_i \b{\theta}_i [ \psi ( \abs{\w{\theta}_i}^2 ) ] [ \bfc ( \q{\theta}_i ) ] (\realization{\theta} ( \bfc ( \q{\theta}_i ) ) - f ( \bfc (  \q{\theta}_i ) ) ) \dens ( \bfc ( \q{\theta}_i ) ) } \\
         &\leq \rbr[\big]{ 2 A^3 \tbound^2 + 8 A \tbound ^2 \br[\big]{ \sup\nolimits_{x \in [a,b]} \abs{\realization{\theta} ( x ) - f ( x ) } } } \rbr[\big]{ \sup\nolimits_{x \in [a,b]} \dens ( x ) },
    \end{split}
\end{equation}
\begin{equation}
    \begin{split}
         \abs[\big]{ \rbr[\big]{\tfrac{\partial^2}{ \partial \theta_i \partial \theta_{\width + j}} \cL } ( \theta) } 
         &\leq 2 \abs{\v{\theta}_i \v{\theta}_j} \abs[\big]{ \tint_{I_i^\theta \cap I_j^\theta} x \dens ( x ) \, \d x } \\
         & \quad + \indicator{[a,b]} ( \q{\theta}_i ) \abs[\big]{2 \v{\theta}_i  [ \psi ( \w{\theta}_i ) ] [ \bfc (\q{\theta}_i ) ] (\realization{\theta} ( \bfc (\q{\theta}_i ) ) - f ( \bfc ( \q{\theta}_i) ) ) \dens ( \bfc ( \q{\theta}_i ) ) } \\
         &\leq  \rbr[\big]{ 2 A^2 \tbound^2 + 4 A \tbound \br[\big]{ \sup\nolimits_{x \in [a,b]} \abs{\realization{\theta} ( x ) - f ( x ) } } } \rbr[\big]{ \sup\nolimits_{x \in [a,b]} \dens ( x ) },
    \end{split}
\end{equation}
and
\begin{equation}
    \begin{split}
         \abs[\big]{ \rbr[\big]{\tfrac{\partial^2}{ \partial \theta_{\width + i} \partial \theta_{\width + j}} \cL } ( \theta) } 
         &\leq 2 \abs{\v{\theta}_i \v{\theta}_j} \abs[\big]{ \tint_{I_i^\theta \cap I_j^\theta} \dens ( x ) \, \d x } \\
         & \quad + \indicator{[a,b]} ( \q{\theta}_i ) \abs[\big]{2 \v{\theta}_i [ \psi ( \w{\theta}_i ) ] (\realization{\theta} ( \bfc ( \q{\theta}_i ) ) - f ( \bfc ( \q{\theta}_i ) ) ) \dens ( \bfc ( \q{\theta}_i ) ) } \\
         &\leq \rbr[\big]{ 2 A \tbound^2 + 4 \tbound \br[\big]{ \sup\nolimits_{x \in [a,b]} \abs{\realization{\theta} ( x ) - f ( x ) } } } \rbr[\big]{ \sup\nolimits_{x \in [a,b]} \dens ( x ) }.
    \end{split}
\end{equation}
Combining this with the fact that $\cu{ A, \tbound } \subseteq [1, \infty)$ establishes \cref{lem:hessian:upperbound:eq}.
\end{cproof}

\begin{lemma} \label{lem:realization:upper:bound}
Assume \cref{setting:snn} and let $\theta \in \R^\fd$,
$A \in \R$ satisfy $A = \max \{ 1 , \abs{a}, \abs{b} \}$. Then 
\begin{equation}
\begin{split} 
    \sup\nolimits_{x \in [a,b]} \abs{\realization{\theta} ( x )  }
    &\leq \abs{\c{\theta}} + A \br*{ \smallsum_{i=1}^\width \abs{\v{\theta}_i} ( \abs{\w{\theta}_i } + \abs{\b{\theta}_i } ) } \\
    &\leq \br*{ \max\nolimits_{i \in \{1, 2, \ldots, \fd \} } \abs{\theta_i} } + 2 A \width \br*{ \max\nolimits_{i \in \{1, 2, \ldots, \fd \} } \abs{\theta_i}^2 } .
    \end{split}
\end{equation}
\end{lemma}
\begin{cproof}{lem:realization:upper:bound}
\Nobs that for all $i \in \{1, 2, \ldots, \width\}$, $x \in [a,b]$ it holds that
\begin{equation}
    \abs{\v{\theta}_i \Rect ( \w{\theta}_i x + \b{\theta}_i) } \leq \abs{\v{\theta}_i} ( \abs{\w{\theta}_i x } + \abs{\b{\theta}_i } ) \leq   \abs{\v{\theta}_i}  ( \abs{\w{\theta}_i } + \abs{\b{\theta}_i } ) A.
\end{equation}
This and the triangle inequality demonstrate that for all $x \in [a,b]$ it holds that
\begin{equation}
\begin{split}
      \abs{\realization{\theta} ( x )  }
      &\leq \abs{\c{\theta}} + \smallsum_{i=1}^\width   \abs{\v{\theta}_i \Rect ( \w{\theta}_i x + \b{\theta}_i) }
      \leq \abs{\c{\theta}} + A \br*{ \smallsum_{i=1}^\width \abs{\v{\theta}_i} ( \abs{\w{\theta}_i } + \abs{\b{\theta}_i } ) } \\
      &\leq \br*{ \max\nolimits_{i \in \{1, 2, \ldots, \fd \} } \abs{\theta_i} } + 2 A \width \br*{ \max\nolimits_{i \in \{1, 2, \ldots, \fd \} } \abs{\theta_i}^2 } .
\end{split}
\end{equation}
\end{cproof}

\begin{cor} \label{cor:hessian:upperbound:max}
Assume \cref{setting:snn},
let $\tbound \in [1, \infty)$, $A \in \R$ satisfy $A = \max \{ 1 , \abs{a}, \abs{b} , b - a \}$,
and let $\theta \in \fV$ satisfy $\max_{i \in \{1, 2, \ldots, \fd \} } \abs{\theta_i} \leq \tbound$ and $\min_{j \in \{1, 2, \ldots, \width \} } \rbr[\big]{ ( \w{\theta}_j - \frac{1}{2} ) \indicator{[a,b]} ( \q{\theta}_j) }  \geq 0$.
Then 
\begin{equation} \label{cor:hessian:upperbound:eq}
\begin{split}
    & \max\nolimits_{i , j \in \{1, 2, \ldots, \fd \} } \abs[\big]{ \rbr[\big]{ \tfrac{\partial^2}{\partial \theta_i \partial \theta_j} \cL } ( \theta )} \\
    &\leq \br[\big]{ 8 A^3 \tbound^2 + 8 A^2 \tbound ^2 \rbr[\big]{\tbound + 2 A \width \tbound^2 + \sup\nolimits_{x \in [a,b]} \abs{ f ( x ) } } } \rbr[\big]{ \sup\nolimits_{x \in [a,b]} \dens ( x ) } \\
    & = \br[\big]{ 8 A^3 \tbound ^2 + 8 A^2 \tbound ^3 + 16 A^3 \width \tbound ^4 + 8 A^2 \tbound ^2 \rbr[\big]{ \sup\nolimits_{x \in [a,b]} \abs{ f ( x ) } } } \rbr[\big]{ \sup\nolimits_{x \in [a,b]} \dens ( x ) } .
     \end{split}
\end{equation}
\end{cor}
\begin{cproof}{cor:hessian:upperbound:max}
\Nobs that \cref{lem:realization:upper:bound} and the triangle inequality prove that for all $x \in [a,b]$ it holds that
\begin{equation}
    \abs{\realization{\theta} ( x ) - f ( x ) } \leq \tbound + 2 A \width \tbound ^2 + \abs{f(x)} \leq \tbound + 2 A \width \tbound^2 + \sup\nolimits_{y \in [a,b]} \abs{ f ( y ) } .
\end{equation}
This and \cref{lem:hessian:max:upperbound} establish \cref{cor:hessian:upperbound:eq}.
\end{cproof}

\section{Regularity properties for the set of 
global minima of the risk function} \label{section:risk:global:min}

In this section we establish in \cref{cor:minima:manifold:1} 
in \cref{subsection:global:min:risk} below under the assumption that the target function is piecewise affine linear that there exists a natural number $k \in \cu{ 1, 2, \ldots, \fd }$ such that a suitable subset of the set of global minima of the considered risk function constitutes a k-dimensional $C^{ \infty }$-submanifold of the ANN parameter space on which the Hessian matrix of the risk function has the maximal rank $\fd - k$. 

Our proof of \cref{cor:minima:manifold:1} employs \cref{prop:minima:manifold:1} 
in \cref{subsection:global:min:risk} as well as the elementary and well-known eigenvalue estimate in \cref{lem:matrix:norm:estimate} in \cref{subsection:global:min:risk}. 
In \cref{prop:minima:manifold:1} we establish under the assumption 
that the target function is piecewise affine linear with varying slopes 
in consecutive subintervals that a suitable subset of the set of global minima 
of the risk function represents an $( \width + 1 )$-dimensional $C^{ \infty }$-submanifold 
of the ANN parameter space on which the Hessian matrix of the risk function 
has the maximal rank $\fd - ( \width + 1 ) = ( 3 \width + 1 ) - ( \width + 1 ) = \width $ 
where $\width \in \N$ represents the number of neurons on the hidden layer (see \cref{setting:snn} for details).

Our proof of \cref{prop:minima:manifold:1} uses \cref{lem:manifold:help} in \cref{subsection:parameter:submanifold}, \cref{prop:hessian:minor:det} in \cref{subsection:hessian:determinant}, and the elementary and well-known properties for tangent spaces of submanifolds in \cref{lem:rank:hessian:upperbound} in \cref{subsection:global:min:risk}. The notion of tanget spaces is recalled in \cref{def:tangent:space} in \cref{subsection:global:min:risk}. Our proof of \cref{prop:hessian:minor:det}, in turn, is based on an application of the auxiliary result in \cref{lem:determinant:positive} in \cref{subsection:hessian:determinant} and in \cref{lem:determinant:positive} and \cref{prop:hessian:minor:det} we show that certain matrices involving appropriate subintegrals of the unnormalized density function have a strictly positive determinant. 

In \cref{lem:manifold:help} in \cref{subsection:parameter:submanifold} we verify that a suitable subset of the ANN parameter space is a non-empty $( \width + 1 )$-dimensional $C^{ \infty }$-submanifold of the ANN parameter space $\R^{ \fd }$. Our proof of \cref{lem:manifold:help} is based on an application of the regular level set theorem which we recall in \cref{lem:regular:value} below. In the scientific literature \cref{lem:regular:value} is sometimes also referred to as submersion level set theorem, regular value theorem, or preimage theorem.
 \cref{lem:regular:value} is, e.g., proved as Theorem 9.9 in Tu~\cite{Tu2011}. 
 Only for the sake of completeness we include in this section the detailed proofs 
 for \cref{lem:rank:hessian:upperbound} and \cref{lem:matrix:norm:estimate}. 
 In the scientific literature \cref{lem:matrix:norm:estimate} is, e.g., proved 
 in Golub \& Van Loan~\cite[Section 2.3.2]{GolubLoan2013}.

\subsection{Submanifolds of the ANN parameter space} \label{subsection:parameter:submanifold}

\begin{prop} \label{lem:regular:value}
Let $\fd, n \in \N$,
let $U \subseteq \R^\fd$ be open,
let $g \in C^\infty ( U, \R^n)$, and assume for all $x \in g^{-1} ( \{ 0 \} )$ that $\rank ( g'(x)) = n$. Then it holds that $g^{-1}(\{ 0 \} ) \subseteq U$ is a $(\fd - n)$-dimensional $C^\infty$-submanifold of $\R^\fd$.
\end{prop}

\begin{lemma} \label{lem:manifold:help}
Assume \cref{setting:snn},
let $\fx_0, \fx_1, \ldots, \fx_\width, \alpha_1, \alpha_2, \ldots, \alpha_\width, \tbound, \fy  \in \R$ satisfy $a = \fx_0 < \fx_1 < \cdots  < \fx_\width = b$ and
\begin{equation}
    \tbound \geq 1 + \abs{\fy} + ( 1 + 2 \max\nolimits_{j \in \{1, 2, \ldots, \width \} } \abs{\alpha_j} ) ( 1 + \abs{a} + \abs{b} ),
\end{equation}
and let $\cM \subseteq \R^{\fd }$ be given by
\begin{multline} 
    \cM = \bigl\{  \theta \in (- \tbound , \tbound )^{\fd } \colon \bigl( \br[\big]{ \min \{ \w{\theta}_1 a + \b{\theta}_1, \w{\theta}_1 b + \b{\theta}_1 , \v{\theta}_1 \} > 0 } ,  \,
    [ \v{\theta}_1 ( \w{\theta}_1 a + \b{\theta}_1 ) + \c{\theta} = \fy ] , \\
    [\w{\theta}_1 \v{\theta}_1 = \alpha_1 ], \, \br[\big]{\forall \, j \in \N \cap (1, \width ] \colon \w{\theta}_j > \nicefrac{1}{2}, \,  \q{\theta}_j = \fx_{j-1}, \, \w{\theta}_j \v{\theta}_j = \alpha_j - \alpha_{j - 1} }\bigr) \bigr\}.
\end{multline}
Then
\begin{enumerate} [label = (\roman*)]
    \item \label{lem:manifold:help:item1} it holds that $\cM \not= \emptyset$ and
    \item \label{lem:manifold:help:item2} it holds that $\cM$ is a $(\width + 1)$-dimensional $C^\infty$-submanifold of $\R^\fd$.
\end{enumerate}
\end{lemma}
\begin{cproof} {lem:manifold:help}
Throughout this proof let $U \subseteq \R^{\fd }$ satisfy
\begin{equation}  \label{lem:manifold:help:eq:defu}
    U = \cu[\big]{ \theta \in (-\tbound , \tbound )^{\fd } \colon \rbr[\big]{ \br[\big]{ \min \{ \w{\theta}_1 a + \b{\theta}_1, \w{\theta}_1 b + \b{\theta}_1 , \v{\theta}_1 \} > 0 },  \, \br[\big]{  \forall \, j \in \N \cap (1, \width ] \colon \w{\theta}_j > \nicefrac{1}{2}} }},
\end{equation}
let $g =(g_1, \ldots, g_{2 \width}) \colon U \to \R^{2 \width}$ satisfy for all $\theta \in U$, $j \in \{1, 2, \ldots, \width\}$ that
\begin{equation} \label{eq:manifold:help:defg1}
    g_j ( \theta ) = \begin{cases}
    \w{\theta}_1 \v{\theta}_1 - \alpha_1 & \colon j=1 \\
    \w{\theta}_j \v{\theta}_j - ( \alpha_j - \alpha_{j-1} )& \colon j>1
    \end{cases}
\end{equation}
and
\begin{equation} \label{eq:manifold:help:defg2}
     g_{\width + j} ( \theta ) = \begin{cases}
    \v{\theta}_1 ( \w{\theta}_1 a + \b{\theta}_1) + \c{\theta} - \fy & \colon j=1 \\
    \q{\theta}_j - \fx_{j-1} & \colon j>1,
    \end{cases}
\end{equation}
and let $\vartheta \in \R^\fd$ satisfy
\begin{multline} \label{eq:manifold:help:nonempty}
   \bigl([\w{\vartheta}_1 = \alpha_1] , \, [\forall \, i \in \N \cap (1, \width] \colon \w{\vartheta}_i = 1], \, [\b{\vartheta}_1 = \abs{\alpha_1} (\abs{a} + \abs{b}) + 1] , \,  [\forall \, i \in \N \cap (1, \width ] \colon \b{\vartheta}_i =- \fx_{i-1} ] , \\
    [ \v{\vartheta}_1 = 1] , \, [\forall \, i \in \N \cap (1, \width ] \colon \v{\vartheta}_i = \alpha_i - \alpha_{i-1} ], \,
    [ \c{\vartheta} = \fy - \v{\vartheta}_1 ( \w{\vartheta}_1 a + \b{\vartheta}_1) ] \bigr) .
\end{multline}
 \Nobs that \cref{eq:manifold:help:nonempty} ensures that $\v{\vartheta}_1 > 0$, $\w{\vartheta}_1 \v{\vartheta}_1 = \alpha_1$, and $\v{\vartheta}_1 ( \w{\vartheta}_1 a + \b{\vartheta}_1 ) + \c{\vartheta} = \fy$.
Moreover, \nobs that $\min\{\w{\vartheta}_1 a + \b{\vartheta}_1 , \w{\vartheta}_1 b + \b{\vartheta }_1 \} = \min \{ \alpha_1 a , \alpha_1 b \} + \abs{\alpha_1} (\abs{a} + \abs{b}) + 1 \geq 1 > 0$.
In addition, \nobs that for all $j \in \N \cap (1, \width]$ we have that $\w{\vartheta}_j = 1 > \nicefrac{1}{2}$, 
$\q{\vartheta}_j = - \nicefrac{\b{\vartheta}_j}{\w{\vartheta}_j} = \fx_{j-1}$, and $\w{\vartheta}_j \v{\vartheta}_j = \alpha_j - \alpha_{j-1}$.
Furthermore, \nobs that for all $i \in \N \cap ( 1, \width]$ it holds that $\abs{\w{\vartheta}_i } = 1 < \tbound$, $\abs{\v{\vartheta}_i } \leq 2 \max_{j \in \{1, 2, \ldots, \width \} } \abs{\alpha_j } < \tbound$, and $\abs{\b{\vartheta}_i} \leq 1 + \abs{a} + \abs{b} < \tbound$.
Moreover, \nobs that $\abs{\w{\vartheta}_1 } = \abs{\alpha_1} < \tbound$, $\abs{\b{\vartheta}_1 } \leq ( 1 +  \max\nolimits_{j \in \{1, 2, \ldots, \width \} } \abs{\alpha_j} ) ( 1 + \abs{a} + \abs{b} ) < \fD$, $\abs{\v{\vartheta}_1} = 1 < \fD$, and 
\begin{equation}
\begin{split}
    \abs{\c{\vartheta}} &\leq \abs{\fy} + \abs{\v{\vartheta}_1 \w{\vartheta}_1 a } + \abs{\v{\vartheta}_1 \b{\vartheta}_1 } 
     = \abs{\fy} + \abs{\alpha_1} \abs{a} + \abs{\alpha_1} ( \abs{a} + \abs{b} ) + 1  \\
    & \leq \abs{\fy} +  \rbr[\big]{ 1 + 2 \max\nolimits_{j \in \{1, 2, \ldots, \width \} } \abs{\alpha_j} } ( 1 + \abs{a} + \abs{b} ) < \tbound.    
\end{split}
\end{equation}
This implies that $\vartheta \in (-\tbound , \tbound ) ^\fd$.
Hence, we obtain that $\vartheta \in \cM$. This establishes \cref{lem:manifold:help:item1}.
In the next step we prove \cref{lem:manifold:help:item2} through an application of the regular value theorem in \cref{lem:regular:value}.
\Nobs that \cref{lem:manifold:help:eq:defu} assures that $U \subseteq \R^\fd$ is open.
In addition, \nobs that the fact that for all $\theta \in U$, $j \in \N \cap (1, \width ]$ it holds that $\w{\theta}_j > 0$ ensures that $g \in C^\infty ( U , \R^{2 \width } )$.
Moreover, \nobs that
\begin{multline}
        g^{-1} ( \{ 0 \} ) = \bigl\{ \theta \in U \colon \bigl( [\w{\theta}_1 \v{\theta}_1 = \alpha_1 ], \, [ \v{\theta}_1 ( \w{\theta}_1 a + \b{\theta}_1 ) + \c{\theta} = \fy ] , \\
         \br[\big]{\forall \, j \in \N \cap (1, \width ] \colon \q{\theta}_j = \fx_{j-1}, \, \w{\theta}_j \v{\theta}_j = \alpha_j - \alpha_{j - 1} } \bigr) \bigr\} .
\end{multline}
This implies that
\begin{multline}
    g^{-1} ( \{ 0 \} ) =  \bigl\{  \theta \in (- \tbound , \tbound )^{\fd } \colon \bigl( \br[\big]{ \min \{ \w{\theta}_1 a + \b{\theta}_1, \w{\theta}_1 b + \b{\theta}_1 , \v{\theta}_1 \} > 0 } , \\
      [\forall \, j \in \N \cap (1, \width ] \colon \w{\theta}_j > \nicefrac{1}{2} ], \, [\w{\theta}_1 \v{\theta}_1 = \alpha_1 ], 
   \, [ \v{\theta}_1 ( \w{\theta}_1 a + \b{\theta}_1 ) + \c{\theta} = \fy ] , \\
    \br[\big]{\forall \, j \in \N \cap (1, \width ] \colon \q{\theta}_j = \fx_{j-1}, \, \w{\theta}_j \v{\theta}_j = \alpha_j - \alpha_{j - 1} }\bigr) \bigr\}
    = \cM.
\end{multline}
Next \nobs that \cref{eq:manifold:help:defg1}, \cref{eq:manifold:help:defg2},
and the fact that for all $\theta \in U$, $j \in \N \cap [1, \width]$ it holds that $\w{\theta}_j = \theta_j$, $\b{\theta}_j = \theta_{\width + j }$, and $\v{\theta}_j = \theta_{2 \width + j}$ ensure that for all $\theta \in U$,
$j \in \N \cap (1, \width]$, $\ell \in \N \cap [1, 2 \width]$ it holds that
\begin{equation}
    \rbr[\big]{\tfrac{\partial}{\partial \theta_{2 \width + j} } g_\ell } ( \theta ) = \begin{cases} \w{\theta}_j \not= 0 & \colon  \ell =  j  \\
    0 & \colon \ell \not= j
    \end{cases}
\end{equation}
and
\begin{equation}
    \rbr[\big]{\tfrac{\partial}{\partial \theta_{\width + j } } g_\ell } ( \theta ) = \begin{cases} - (\w{\theta}_j)^{-1} \not= 0 & \colon \ell = \width + j \\
    0 & \colon \ell \not= \width + j .
    \end{cases}
\end{equation}
In addition, \nobs that \cref{eq:manifold:help:defg1} and \cref{eq:manifold:help:defg2} show that for all $\theta \in U$,
$\ell \in \N \cap [1, 2 \width]$ it holds that
\begin{equation}
    \rbr[\big]{\tfrac{\partial}{\partial \theta_{ 1} } g_\ell } ( \theta ) = \begin{cases} \v{\theta}_1 \not= 0 & \colon \ell = 1 \\
    \v{\theta}_1 a & \colon \ell = \width + 1 \\
    0 & \colon \ell \notin  \{ 1 , \width + 1 \}
    \end{cases}
\end{equation}
and
\begin{equation}
    \rbr[\big]{\tfrac{\partial}{\partial \theta_{ \width + 1 } } g_\ell } ( \theta ) = \begin{cases} \v{\theta}_1 \not= 0 & \colon \ell = \width + 1  \\
    0 & \colon  \ell \not= \width + 1 .
    \end{cases}
\end{equation}
This demonstrates that for all $\theta \in U$ it holds that the $( ( 2 \width ) \times ( 2 \width ) )$-matrix with entries $\rbr[\big]{\tfrac{\partial}{\partial \theta_{i } } g_\ell } ( \theta ) \in \R$, $(i , \ell ) \in ( \{1 \} \cup \{\width + j \colon j \in \N \cap [ 1, \width ]\} \cup \{ 2 \width + j \colon j \in \N \cap (1, \width ] \} ) \times \{1,2, \ldots, 2 \width \}$, is invertible. 
Hence, we obtain for all $\theta \in U$ that $\rank ( g' ( \theta ) ) = 2 \width $. Combining this with \cref{lem:regular:value} establishes \cref{lem:manifold:help:item2}.
\end{cproof}

\subsection{Determinants of submatrices of the Hessian matrix of the risk function} \label{subsection:hessian:determinant}

\begin{lemma} \label{lem:determinant:positive}
Let $a \in \R$, $b \in (a, \infty)$, let $\dens \colon [a,b] \to (0, \infty)$ be bounded and measurable, 
let $\cQ_N \subseteq \R^{N + 1 }$, $N \in \N$, satisfy for all $N \in \N$ that $\cQ_N = \cu{ \fx = (\fx_1, \ldots, \fx_{N+1} ) \in \R^{N+1} \colon a \le \fx_1 < \fx_2 < \cdots < \fx_{N+1} \le b }$,
and let $A^{N,\fx} = (A^{N,\fx}_{i,j} )_{(i,j) \in \cu{1, 2, \ldots, 2N} ^2 } \in \R^{(2N) \times (2N) }$, $\fx \in \cQ_N$,
$N \in \N$,
satisfy for all $N \in \N$, $\fx=(\fx_1, \ldots, \fx_{N+1} ) \in \cQ_N$, $i,j \in \cu{1, 2, \ldots, N} $ that
\begin{multline}
A^{N,\fx}_{i,j} = \tint_{\fx_{\max \{i,j \}} }^{\fx_{N+1}} x^2 \dens ( x ) \, \d x, \qquad A^{N,\fx}_{N + i, j} = A^{N,\fx}_{i, N + j} = \tint_{  \fx_{\max \{i,j \} } }^{\fx_{N+1}} x \dens ( x ) \, \d x, \\
\text{and} \qquad A^{N,\fx}_{N + i, N + j} = \tint_{  \fx_{\max \{i,j \} } }^{\fx_{N+1}} \dens ( x ) \, \d x.    
\end{multline}
Then it holds for all $N \in \N$, $\fx \in \cQ_N$ that
\begin{equation} \label{lem:det:positive:eq:claim}
\det (A^{N, \fx }) = \prod_{i=1}^N \rbr*{ \br*{ \tint _{\fx_i}^{\fx_{i+1} } x^2 \dens ( x ) \, \d x } \br*{ \tint _{\fx_i}^{\fx_{i+1} }  \dens ( x ) \, \d x } - \br*{ \tint _{\fx_i}^{\fx_{i+1} } x \dens ( x ) \, \d x } ^2 } > 0 .
\end{equation}
\end{lemma}
\begin{cproof}{lem:determinant:positive}
Throughout this proof let $E_i^{N, \fx} \in \R$, $i \in \cu{1, 2, \ldots, N}$, $\fx \in \cQ_N$, $N \in \N$, satisfy for all $N \in \N$, $\fx \in \cQ_N$, $i \in \cu{1, 2, \ldots, N}$ that
\begin{equation}
    E_i^{N, \fx} =  \br*{ \tint _{\fx_i}^{\fx_{i+1} } x^2 \dens ( x ) \, \d x } \br*{ \tint _{\fx_i}^{\fx_{i+1} }  \dens ( x ) \, \d x } - \br*{ \tint _{\fx_i}^{\fx_{i+1} } x \dens ( x ) \, \d x } ^2 .
\end{equation}
\Nobs that the Cauchy-Schwarz inequality and the fact that for all $x \in [a,b]$ it holds that $\dens ( x ) > 0 $ ensure that for all $N \in \N$, $\fx \in \cQ_N$, $i \in \cu{1, 2, \ldots, N}$ it holds that
\begin{equation}
\begin{split}
    \abs*{ \tint_{\fx_i}^{\fx_{i+1}} x \dens ( x ) \, \d x } 
    &= 
    \abs*{ \tint_{\fx_i}^{\fx_{i+1}} \br[\big]{ x \sqrt{ \dens ( x )} } \br[\big]{ \sqrt{\dens ( x ) } } \, \d x }  \\
    &< \br*{ \tint _{\fx_i}^{\fx_{i+1} } x^2 \dens ( x ) \, \d x } ^{\nicefrac{1}{2}} \br*{ \tint _{\fx_i}^{\fx_{i+1} }  \dens ( x ) \, \d x }^{\nicefrac{1}{2}}.
\end{split}
\end{equation}
Hence, we obtain for all $N \in \N$, $\fx \in \cQ_N$, $i \in \cu{1, 2, \ldots, N}$
that $E_i^{N, \fx} > 0$. Next we claim that for all $N \in \N$, $x \in \cQ_N$ it holds that
\begin{equation} \label{lem:det:positive:eq:inductclaim}
    \det ( A ^{N , \fx } ) = \textstyle\prod_{i=1}^N E_i^{N, \fx} > 0.
\end{equation}
We now prove \cref{lem:det:positive:eq:inductclaim} by induction on $N \in \N$. 
For the base case $N=1$ \nobs that for all $\fx = (\fx_1, \fx_2 ) \in \cQ_1$ it holds that
\begin{equation}
       \det ( A ^{1 , \fx } ) = \det \!  \begin{pmatrix} \tint_{x_1}^{x_2} x^2 \dens ( x ) \, \d x & \tint_{x_1}^{x_2} x \dens ( x ) \, \d x \\
    \tint_{\fx_1}^{\fx_2} x \dens ( x ) \, \d x & \tint_{\fx_1}^{\fx_2} \dens ( x ) \, \d x 
    \end{pmatrix}
    = E_1^{1, x} > 0.
\end{equation}
This establishes \cref{lem:det:positive:eq:inductclaim} in the base case $N=1$. For the induction step let $N \in \N \cap [2, \infty )$ and assume for all $\fx \in \cQ_{N - 1 }$ that
\begin{equation} \label{lem:det:positive:eq:hypo}
    \det ( A^{N - 1, \fx} ) = \textstyle\prod_{i=1}^{N - 1 } E_i^{N - 1, \fx } > 0.
\end{equation}
Next let $\fx = (\fx_1, \ldots, \fx_{N + 1 } ) \in \cQ_{N}$
and let $B = (B_{i,j} )_{ ( i , j ) \in \cu{1, 2, \ldots, 2N } ^2 } \in \R^{(2N ) \times (2N)}$ satisfy for all $i,j \in \cu{1, 2, \ldots, 2N} $ that
\begin{equation} \label{lem:det:positive:eq:defb}
    B_{i,j} = \begin{cases}
            A_{i,j}^{N, \fx } & \colon i \notin \cu{1, N+1 } \\[1ex] 
            A_{1,j}^{N, \fx } - A_{2,j}^{N, \fx } & \colon i=1 \\[1ex] 
            A_{N+1, j}^{N, \fx } - A_{N+2, j}^{N, \fx } & \colon i=N+1.
    \end{cases}
\end{equation}
\Nobs that
$B$ is the matrix that is obtained from $A^{N, \fx}$ by subtracting the 2nd row from the 1st row and the $(N+2)$-th row from the $(N+1)$-th row.
In particular, \nobs that \cref{lem:det:positive:eq:defb} implies that $\det ( B ) = \det ( A^{N, \fx} ) $.
Next \nobs that the fact that for all $j \in \N \cap (1 , N]$ it holds that $A^{N, \fx}_{1,j} = A^{N, \fx}_{2,j}$, $A^{N, \fx}_{1, N + j} = A^{N , \fx}_{2, N + j}$, $A^{N, \fx}_{N+1, j} = A^{N, \fx}_{N+2, j}$, and $A^{N, \fx}_{N+1, N+j} = A^{N, \fx}_{N+2, N+j}$ demonstrates that for all $i, j \in \N \cap ( 1 , N]$ we have that
\begin{equation}
\begin{split}
    B_{1,1} &= A_{1,1}^{N, \fx} - A_{2,1}^{N, \fx} = \tint_{\fx_1}^{\fx_{N+1}} x^2 \dens ( x ) \, \d x - \tint_{\fx_2}^{\fx_{N+1}} x^2 \dens ( x ) \, \d x
    =\tint_{\fx_1}^{\fx_2} x^2 \dens ( x ) \, \d x, \\
    B_{N+1, 1} &= B_{1, N+1} =  \tint_{\fx_1}^{\fx_{N+1}} x \dens ( x ) \, \d x - \tint_{\fx_2}^{\fx_{N+1}} x \dens ( x ) \, \d x
    = \tint_{\fx_1}^{\fx_2} x \dens ( x ) \, \d x , \\
    B_{N+1, N+1} &= \tint_{\fx_1}^{\fx_{N+1}} \dens ( x ) \, \d x - \tint_{\fx_2}^{\fx_{N+1}}  \dens ( x ) \, \d x
    = \tint_{\fx_1}^{\fx_2} \dens ( x ) \, \d x, \\
    B_{1, j} &= B_{N+1, j} = B_{1, N + j} = B_{N+1, N+j} = 0, \quad B_{i,j} = A^{N, \fx}_{i,j}, \quad B_{N+i, j} = A^{N, \fx}_{N+i, j}, \\
    B_{i, N+j} &= A^{N, \fx}_{i, N+j}, \qandq B_{N+i, N+j} = A^{N, \fx }_{N+i, N+j}.
    \end{split}
\end{equation}
Hence, we obtain that
\begin{equation}
\begin{split}
    \det ( B ) &= ( B_{1,1} B_{N+1, N+1} - B_{N+1, 1} B_{1, N+1} ) \det \rbr*{ (B_{i,j})_{( i , j) \in (\{1, 2, \dots, 2 N\} \backslash \{ 1, N+1 \})^2 } } \\
    &= E_1^{N, \fx} \det \rbr*{ (B_{i,j})_{(i, j) \in ( \{1, 2, \dots, 2 N\} \backslash \{ 1, N+1 \} )^2 } }.
    \end{split}
\end{equation}
In addition, \nobs that \cref{lem:det:positive:eq:hypo} proves that
\begin{equation}
\begin{split}
    \det \rbr*{ (B_{i,j})_{( i , j ) \in ( \{1, \ldots, 2 N\} \backslash \{ 1, N+1 \} ) ^2 } } 
    &= \det (A^{N-1, (\fx_2, \fx_3, \ldots, \fx_{N+1} )} ) \\
    &= \textstyle\prod_{i=1}^{N-1} E_i^{N-1 , (\fx_2, \fx_3, \ldots, \fx_{N+1} ) } = \textstyle\prod_{i = 2}^{N} E_i^{N, \fx} > 0.
\end{split}
\end{equation}
Hence, we obtain that $\det ( A^{N , \fx } ) = \det ( B ) = \prod_{i=1}^N E_i^{N, \fx }$.
Induction thus proves \cref{lem:det:positive:eq:inductclaim}. Furthermore, \nobs that \cref{lem:det:positive:eq:inductclaim} establishes \cref{lem:det:positive:eq:claim}.
\end{cproof}

\begin{prop} \label{prop:hessian:minor:det}
Let $N \in \N$, $v_1, v_2, \ldots, v_N \in \R \backslash \{0 \}$, $\fx_0, \fx_1, \ldots, \fx_{N} \in \R$ satisfy $\fx_0 < \fx_1 < \cdots < \fx_{N }$,
let $I_j \subseteq \R$, $j \in \{1, 2, \ldots, N\}$, satisfy for all $j \in \{1, 2, \ldots, N\}$ that $I_j = [\fx_{j - 1}, \fx_{N}]$,
let $\dens \colon [\fx_0, \fx_{N}] \to (0, \infty)$ be bounded and measurable,
and let $A = (A_{i,j})_{( i , j ) \in \{1, 2, \ldots, 2N \} ^2} \in \R^{ ( 2 N ) \times ( 2 N ) }$ satisfy for all $i,j \in \{1, 2, \ldots, N\}$ that 
\begin{multline} \label{prop:hessian:minor:det:eq1}
    A_{i,j} =2 v_i v_j \tint_{I_i \cap I_j} x^2 \dens ( x ) \, \d x, \qquad A_{N + i, j} = A_{i, N + j} =2 v_i v_j \tint_{I_i \cap I_j} x \dens ( x ) \, \d x, \\  \text{and} \qquad A_{N + i, N + j} =2 v_i v_j \tint_{I_i \cap I_j } \dens ( x ) \, \d x.
\end{multline}
Then $ \det ( A ) > 0$.
\end{prop}
\begin{cproof}{prop:hessian:minor:det}
Throughout this proof let $B = (B_{i,j})_{( i , j ) \in \{1, 2, \ldots, 2N \} ^2 } \in \R^{( 2 N )  \times ( 2 N ) }$ satisfy for all $i,j \in \{1, 2, \ldots, N\}$ that $B_{i,j} = \int_{I_i \cap I_j} x^2 \dens ( x ) \, \d x$, $B_{N + i, j} = B_{i, N + j} = \int_{I_i \cap I_j} x \dens ( x ) \, \d x$, and $B_{N + i, N + j} =  \int_{I_i \cap I_j } \dens ( x ) \, \d x$.
\Nobs that for all $i,j \in \{1, 2, \ldots, N\}$ it holds that 
\begin{multline} \label{prop:hessian:minor:det:eq2}
B_{i,j} = \tint_{\fx_{\max \{i - 1 , j - 1 \}} }^{\fx_{N}} x^2 \dens ( x ) \, \d x, \qquad B_{N + i, j} = B_{i, N + j} = \tint_{  \fx_{\max \{i - 1 , j - 1 \} } }^{\fx_{N}} x \dens ( x ) \, \d x, \\
\text{and} \qquad B_{N + i, N + j} = \tint_{  \fx_{\max \{i - 1 , j - 1 \} } }^{\fx_{N}} \dens ( x ) \, \d x.    
\end{multline}
Furthermore, \nobs that \cref{prop:hessian:minor:det:eq1}
and the fact that the determinant is linear in each row and each column show that
\begin{equation} \label{prop:hessian:minor:det:eq3}
    \det ( A ) = 4^N \rbr[\big] { \textstyle\prod_{i=1}^N \abs{v_i}^4 } \det  (B ).
\end{equation}
In addition, \nobs that \cref{prop:hessian:minor:det:eq2} and \cref{lem:determinant:positive} (applied with $a \with \fx_0$, $b \with \fx_{N}$, $\dens \with \dens$, $N \with N$, $\fx \with (\fx_0, \fx_1, \ldots, \fx_{N} ) $ in the notation of \cref{lem:determinant:positive}) demonstrate that
$\det ( B ) > 0$. Combining this with \cref{prop:hessian:minor:det:eq3} ensures that $\det ( A ) > 0 $.
\end{cproof}

\subsection{Regularity properties for the set of global minima of the risk function} \label{subsection:global:min:risk}

\begin{definition} [Tangent space] \label{def:tangent:space}
Let $\fd \in \N$,
let $\cM \subseteq \R^\fd$ be a set,
 and let $x \in \cM$.
Then we denote by $\cT_\cM^x \subseteq \R^\fd$ the set given by
\begin{equation}
    \cT_\cM^x = \cu*{v \in \R^\fd \colon \br*{\exists \, \gamma \in C^1 ( \R , \R^\fd ) \colon \rbr*{ [ \gamma ( \R ) \subseteq \cM ] , \, [ \gamma ( 0 ) = x ] , \, [ \gamma'(0) = v ] } } }.
\end{equation}
\end{definition}

\cfclear
\begin{lemma} \label{lem:rank:hessian:upperbound}
 Let $\fd , \mdim \in \N$, 
 let $U \subseteq \R^\fd$ be open, let $f \in C^2 ( U , \R) $ have locally Lipschitz continuous derivatives, let $\cM \subseteq U$ satisfy $\cM = \{ x \in U \colon f(x) = \inf_{y \in U} f(y) \}$, assume that $\cM $ is a $\mdim$-dimensional $C^2$-submanifold of $\R^\fd$,
 and let $x \in \cM $. Then
\begin{enumerate} [label = (\roman*)] \cfadd{def:tangent:space}
\item \label{lem:rank:hessian:upperbound:item1} it holds for all $v \in \cT_\cM^x$ that $\rbr[\big]{ (\Hs  f)(x) } v = 0$,
    \item \label{lem:rank:hessian:upperbound:item2} it holds that $\rank ( (\Hs  f)(x) ) \leq \fd - \mdim$, and
      \item \label{lem:rank:hessian:upperbound:item3}  it holds for all $v \in ( \cT_\cM^x ) ^\perp$ that $ \rbr[\big]{(\Hs  f)(x) } v \in ( \cT_\cM^x)^\perp$
\end{enumerate}
\cfload.
\end{lemma}
\begin{cproof} {lem:rank:hessian:upperbound}
 \Nobs that the assumption that $\cM = \cu{y \in U \colon f(y) = \inf_{z \in U} f(z) }$ ensures for all $y \in \cM$ that $(\nabla f ) ( y ) = 0$. 
 This implies for all $\gamma \in C^1 ( \R , \R^\fd ) $, $t \in \R$ with $\gamma ( \R ) \subseteq \cM$ that $(\nabla f ) ( \gamma ( t ) ) = 0$. Hence, we obtain for all
 $\gamma \in C^1 ( \R , \R^\fd ) $, $t \in \R $ with $\gamma ( \R ) \subseteq \cM$ that
 \begin{equation}
     0 = \tfrac{\d}{\d t} \rbr[\big]{ (\nabla f)(\gamma ( t ) ) } = \rbr[\big]{ (\Hs  f)(\gamma(t)) } \gamma' ( t ) .
 \end{equation}
This shows for all $\gamma \in C^1 ( \R , \R^\fd ) $ with $\gamma ( \R ) \subseteq \cM$ and $\gamma ( 0 ) = x $ that $( ( \Hs f ) ( x ) ) \gamma ' ( 0 ) = 0$.
This establishes \cref{lem:rank:hessian:upperbound:item1}.

Next \nobs that the assumption that $\cM $ is a $\mdim$-dimensional $C^2$-submanifold of $\R^\fd$
proves that $\dim ( \cT_\cM^x) = \mdim$. Combining this with \cref{lem:rank:hessian:upperbound:item1} establishes \cref{lem:rank:hessian:upperbound:item2}.

Moreover, \nobs that \cref{lem:rank:hessian:upperbound:item1} and the fact that $(\Hs  f)(x)$ is symmetric demonstrate for all $v \in \cT_\cM^x$, $w \in (\cT_\cM^x)^\perp $ that
\begin{equation}
    \spro*{v , \rbr[\big]{ ( \Hs f ) ( x ) } w } = \spro*{ \rbr[\big]{ ( \Hs f ) ( x ) } v , w } = \spro{0 , w } = 0.
\end{equation}
This establishes \cref{lem:rank:hessian:upperbound:item3}.
\end{cproof}

\begin{prop} \label{prop:minima:manifold:1}
Assume \cref{setting:snn},
let $\fx_0, \fx_1, \ldots, \fx_\width, \alpha_1, \alpha_2, \ldots, \alpha_\width \in \R$ satisfy $a = \fx_0 < \fx_1 < \cdots  < \fx_\width = b$,
assume for all $i \in \{1, 2, \ldots, \width \}$, $x \in [\fx_{i-1}, \fx_i ]$ that $f(x) = f( \fx_{i-1} ) + \alpha_i ( x - \fx_{i-1} )$, 
assume $ \prod_{i=1}^{\width - 1 } ( \alpha_{i+1} - \alpha_i )  \not= 0$,
and let $\tbound \in \R$ satisfy
\begin{equation}
    \tbound = 1 + \abs{f(a)} + ( 1 + 2 \max\nolimits_{j \in \{1, 2, \ldots, \width \} } \abs{\alpha_j} ) ( 1 +  \abs{a} + \abs{b} ).
\end{equation}
Then there exists an open $U \subseteq (- \tbound , \tbound )^{\fd }$ such that
\begin{enumerate} [label=(\roman*)]
\item \label{prop:minima:manifold:1:item0} it holds that $U \subseteq \fV$,
\item \label{prop:minima:manifold:1:item1} it holds that $\cL |_U \in C^2( U , \R)$,
\item \label{prop:minima:manifold:1:item2} it holds that $U \ni \theta \mapsto (\Hs  \cL) ( \theta ) \in \R^{\fd \times \fd}$ is locally Lipschitz continuous,
\item \label{prop:minima:manifold:1:itemhess} it holds for all $\theta = (\theta_1, \ldots, \theta_\fd) \in U$ that 
\begin{equation}
    \max\nolimits_{i,j \in \{1, 2, \ldots, \fd \} } \abs[\big]{ \rbr[\big]{ \tfrac{\partial^2}{\partial \theta_i \partial \theta_j} \cL } ( \theta )} \leq \rbr[\big]{24 \tbound ^5 + 16 \width \tbound ^7 } \rbr[\big]{ \sup\nolimits_{x \in [a,b]} \dens ( x ) } ,
\end{equation}
    \item \label{prop:minima:manifold:1:item3} it holds that $\{ \vartheta \in U \colon \cL ( \vartheta ) = 0\} \not= \emptyset$,
    \item \label{prop:minima:manifold:1:item4} it holds that $\{ \vartheta \in U \colon \cL ( \vartheta ) = 0\}$ is a $(\width + 1)$-dimensional $C^\infty$-submanifold of $\R^\fd$, and
    \item \label{prop:minima:manifold:1:item5}  it holds for all $\theta \in \{ \vartheta \in U \colon \cL ( \vartheta ) = 0\}$ that $ \rank ( (\Hs  \cL ) ( \theta ) ) = 2 \width = \fd  - (\width + 1)$.
\end{enumerate}
\end{prop}
\cfclear
\begin{cproof}{prop:minima:manifold:1}
Throughout this proof let $U \subseteq \R^{\fd }$ satisfy
\begin{multline} \label{prop:minima:1:eq:defu}
    U = \bigl\{ \theta \in ( - \tbound , \tbound )^{\fd } \colon \bigl( \br[\big]{ \min \{ \w{\theta}_1 a + \b{\theta}_1, \w{\theta}_1 b + \b{\theta}_1 , \v{\theta}_1 \} > 0 }, \, \br[\big]{  \forall \, j \in \N \cap (1, \width ] \colon \w{\theta}_j > \nicefrac{1}{2}} , \\
    \br[\big]{\forall \, j \in \N \cap (1, \width] \colon \q{\theta}_j \in (a,b)}, \br[\big]{\forall \, j \in \N \cap (1, \width ) \colon \q{\theta}_j < \q{\theta}_{j+1} }  \bigr) \bigr\}
\end{multline}
and let $\cM \subseteq \R^{\fd }$ be given by
\begin{multline}  \label{prop:minima:1:eq:defm}
    \cM = \bigl\{  \theta \in (- \tbound , \tbound )^{\fd } \colon \bigl( \br[\big]{ \min \{ \w{\theta}_1 a + \b{\theta}_1, \w{\theta}_1 b + \b{\theta}_1 , \v{\theta}_1 \} > 0 } ,  \,  
    [ \v{\theta}_1 ( \w{\theta}_1 a + \b{\theta}_1 ) + \c{\theta} = f( a ) ] , \\
    [\w{\theta}_1 \v{\theta}_1 = \alpha_1 ], \,
    \br[\big]{\forall \, j \in \N \cap (1, \width ] \colon \w{\theta}_j > \nicefrac{1}{2} , \,  \q{\theta}_j = \fx_{j-1}, \, \w{\theta}_j \v{\theta}_j = \alpha_j - \alpha_{j - 1} } \bigr) \bigr\}.
\end{multline}
\Nobs that \cref{prop:minima:1:eq:defu} ensures that $U$ is open. Furthermore, \nobs that \cref{setting:eq:defv} and \cref{prop:minima:1:eq:defu} assure that $U \subseteq \fV$. This proves \cref{prop:minima:manifold:1:item0}. In addition,
\nobs that \cref{prop:minima:manifold:1:item0}, \cref{lem:twice:diff:lipschitz}, and the fact that $U$ is open establish \cref{prop:minima:manifold:1:item1,prop:minima:manifold:1:item2}.

Next \nobs that \cref{cor:hessian:upperbound:max}, the fact that for all $\theta \in U$,
$j \in \{1, 2, \ldots, \width\}$ with $\q{\theta}_j \in [a,b]$ it holds that $\w{\theta}_j > \frac{1}{2}$,
and the fact that $\tbound \geq \max \{ \abs{a}, \abs{b}, b - a , \sup\nolimits_{x \in [a,b]} \abs{ f ( x ) }, 1 \} \geq 1$ prove that for all $\theta \in U \subseteq (- \tbound , \tbound ) ^\fd$ we have that
\begin{equation}
    \begin{split}
        \max_{i,j \in \{1, 2, \ldots, \fd \} } \abs[\big]{ \rbr[\big]{ \tfrac{\partial^2}{\partial \theta_i \partial \theta_j} \cL } ( \theta )} 
        & \leq \rbr[\big]{ 16 \tbound ^5 + 16 \width \tbound ^7 + 8 \tbound^4 \rbr[\big]{ \sup\nolimits_{x \in [a,b]} \abs{ f ( x ) } } } \rbr[\big]{ \sup\nolimits_{x \in [a,b]} \dens ( x ) } \\
        & \leq \rbr[\big]{24 \tbound ^5 + 16 \width \tbound ^7 } \rbr[\big]{ \sup\nolimits_{x \in [a,b]} \dens ( x ) }.
    \end{split}
\end{equation}
This establishes \cref{prop:minima:manifold:1:itemhess}.

Next \nobs that \cref{prop:minima:1:eq:defm} and \cref{lem:manifold:help} imply that $\cM$ is a non-empty $(\width +1)$-dimensional $C^\infty$-submanifold of $\R^{\fd }$.
Furthermore, \nobs that \cref{prop:minima:1:eq:defu}, \cref{prop:minima:1:eq:defm}, and the fact that $a < \fx_1 < \fx_2 < \cdots < \fx_\width = b$ show that $\cM \subseteq U$.
In the next step we intend to prove that for all $\theta \in \cM$ it holds that $\cL ( \theta ) = 0$. \Nobs that \cref{prop:minima:1:eq:defu} and the fact that for all $\theta \in U$, $x \in [a,b]$ it holds that 
\begin{equation} \label{prop:minima:1:eq:positive}
\w{\theta}_1 x + \b{\theta}_1 =  \br*{ \tfrac{b - x}{b - a} } (\w{\theta}_1 a + \b{\theta}_1 ) + \br*{ \tfrac{x - a}{b - a} } (\w{\theta}_1 b + \b{\theta}_1 )  > 0    
\end{equation}
ensure that for all $\theta \in U$, $x \in [a , b]$ it holds that
\begin{equation} \label{prop:minima:1:eq:realization}
    \begin{split}
         \realization{\theta} ( x ) &= \c{\theta} + \v{\theta}_1 \max \{\w{\theta}_1 x + \b{\theta}_1, 0 \} + \smallsum_{j=2}^\width \v{\theta}_j \max \{\w{\theta}_j x + \b{\theta}_j, 0 \} \\
        &= \c{\theta} + \v{\theta}_1 ( \w{\theta}_1 x + \b{\theta}_1 ) + \smallsum_{j=2}^\width \v{\theta}_j \w{\theta}_j \max \{ x - \q{\theta}_{j}, 0 \} .
    \end{split}
\end{equation}
Combining this with \cref{prop:minima:1:eq:defm} demonstrates that for all $\theta \in \cM$, $x \in [a,b]$ we have that 
\begin{equation} \label{prop:minima:1:eq:help1}
    \begin{split}
        \realization{\theta} ( x )
        &= \v{\theta}_1 \w{\theta}_1 x + \v{\theta}_1\b{\theta}_1 + \c{\theta} + \smallsum_{j=2}^\width \v{\theta}_j \w{\theta}_j \max \{ x - \fx_{j-1}, 0 \} \\
        &= \v{\theta}_1 \w{\theta}_1 x + f ( a ) - \v{\theta}_1 \w{\theta}_1 a + \smallsum_{j=2}^\width \v{\theta}_j \w{\theta}_j \max \{ x - \fx_{j-1}, 0 \} \\
        &= f( a ) + \alpha_1 (x - a) +  \smallsum_{j=2}^\width ( \alpha_j - \alpha_{j-1} ) \max \{ x - \fx_{j-1}, 0 \}.
    \end{split}
\end{equation}
In addition, \nobs that the assumption that for all $i \in \{1, 2, \ldots \width \}$, $x \in [\fx_{i-1}, \fx_i ]$ it holds that $f(x) = f( \fx_{i-1} ) + \alpha_i ( x - \fx_{i-1} )$ proves that for all $j \in \{0, 1, \ldots, \width - 1 \}$, $x \in [\fx_j, \fx_{j+1}]$ it holds that
\begin{equation} 
    \begin{split}
      f(x)  &= f( \fx_0 ) + \br[\big]{\smallsum_{k=1}^j [ f(\fx_k) - f(\fx_{k-1} ) ] } + [f(x) - f(\fx_j ) ] \\
    &= f(a) + \br[\big]{ \smallsum_{k=1}^j \alpha_k ( \fx_k - \fx_{k-1}) } + \alpha_{j+1} (x-\fx_j) \\
    &= f ( a ) + \alpha_{j+1} x + \br[\big]{ \smallsum_{k=1}^j \alpha_k ( \fx_k - \fx_{k-1} ) } - \alpha_{j+1} \fx_j \\
     &= f(a) + \alpha_{j+1} x + \br[\big]{ \smallsum_{k=1}^{j} \alpha_{k} \fx_{k} } - \br[\big]{ \smallsum_{k=1}^{j} \alpha_k \fx_{k-1} } - \alpha_{j+1} \fx_j \\
        &= f(a) + \alpha_{j+1} x - \rbr[\big]{ \br[\big]{ \smallsum_{k=1}^{j+1} \alpha_k \fx_{k-1} } - \br[\big]{ \smallsum_{k=1}^{j} \alpha_{k} \fx_{k} } } \\
         &= f(a) + \alpha_{j+1} x - \rbr[\big]{ \alpha_1 \fx_0 + \br[\big]{ \smallsum_{k=2}^{j+1} \alpha_k \fx_{k-1} } - \br[\big]{ \smallsum_{k=2}^{j+1} \alpha_{k-1} \fx_{k-1} } } \\
          & = f( a ) + \rbr[\big]{ \alpha_1 x + \br[\big]{ \smallsum_{k=2}^{j+1} (\alpha_k - \alpha_{k-1}) x } } - \rbr[\big]{ \alpha_1 \fx_0 + \br[\big]{ \smallsum_{k=2}^{j+1} (\alpha_k - \alpha_{k-1}) \fx_{k-1} } }\\
          &= f( a ) + \alpha_1 (x - a) + \smallsum_{k = 2}^{j+1} ( \alpha_k - \alpha_{k - 1} ) ( x - \fx_{k - 1} ) \\
          &= f( a ) + \alpha_1 (x - a) +  \smallsum_{k = 2}^{\width} ( \alpha_k - \alpha_{k - 1} ) \max \{ x - \fx_{k - 1}, 0 \}.
    \end{split}
\end{equation}
This implies that for all $x \in [a,b]$ we have that
\begin{equation} \label{prop:minima:1:eq:help2}
     f(x) = f( a ) + \alpha_1 (x - a) +  \smallsum_{j=2}^\width ( \alpha_j - \alpha_{j-1} ) \max \{ x - \fx_{j-1}, 0 \}.
\end{equation}
Combining this with \cref{prop:minima:1:eq:help1} demonstrates that for all $\theta \in \cM$, $x \in [a , b]$ it holds that $\realization{\theta} ( x ) = f ( x )$. Hence, we obtain that for all $\theta \in \cM$ it holds that $\cL ( \theta ) = 0$.
Next we intend to prove that for all $\theta \in U$ with $\cL ( \theta ) = 0$ it holds that $\theta \in \cM$. \Nobs that \cref{setting:eq:def:risk} and the fact that for all $\theta \in \R^\fd$ it holds that $[a , b] \ni x \mapsto \realization{\theta} ( x ) - f ( x ) \in \R$ is continuous show that for all $\theta \in \{ \vartheta \in U \colon \cL ( \vartheta ) = 0\}  \subseteq \R^\fd$, $x \in [a,b]$ we have that 
\begin{equation} \label{prop:minima:1:eq:help4}
    \realization{\theta} ( x )= f ( x ). 
\end{equation}
Combining this with \cref{prop:minima:1:eq:defu},
\cref{prop:minima:1:eq:defm},
\cref{prop:minima:1:eq:realization},
and the fact that $\cM \subseteq U$ demonstrates that for all $\theta \in \{ \vartheta \in U \colon \cL ( \vartheta ) = 0\} $,
$x \in [\fx _ 0 , \fx_1 +  \min \{0 , ( \q{\theta}_{ \min \{ 2 , \width \} } - \fx_1 ) \indicator{(1, \infty )} ( \width ) \} ]$ it holds that
\begin{equation}
  f(a) + \alpha_1 ( x - a) = f ( x ) =   \realization{\theta} ( x ) = \v{\theta}_1 ( \w{\theta}_1 x + \b{\theta}_1 ) + \c{\theta} = \v{\theta}_1 \w{\theta}_1 ( x - a ) + \v{\theta}_1 ( \w{\theta}_1 a + \b{\theta}_1 ) + \c{\theta} .
\end{equation}
The fact that for all $\theta \in U$ it holds that $\fx_1 +  \min \{0 , ( \q{\theta}_{ \min \{ 2 , \width \} } - \fx_1 ) \indicator{(1, \infty )} ( \width ) \} > \fx_0 $ hence ensures that for all $\theta \in \{ \vartheta \in U \colon \cL ( \vartheta ) = 0\} $ we have that 
\begin{equation} \label{prop:minima:1:eq:help3}
    \w{\theta}_1 \v{\theta}_1 = \alpha_1 \qquad \text{and} \qquad \v{\theta}_1 ( \w{\theta}_1 a + \b{\theta}_1 ) + \c{\theta} = f ( a ).
\end{equation}
Next \nobs that the fact that for all $\theta \in U$ it holds that $(a,b) \backslash \{ \q{\theta}_1 , \q{\theta}_2, \ldots, \q{\theta}_\width \}$ is an open set shows that there exists $\varepsilon = (\varepsilon_{\theta , x } )_{(\theta, x ) \in U \times \R } \colon U \times \R \to (0, \infty)$ which satisfies
for all $\theta \in U$, $x \in (a,b) \backslash \{ \q{\theta}_1 , \q{\theta}_2, \ldots, \q{\theta}_\width \}$ that $(x - \varepsilon_{\theta , x }, x + \varepsilon_{\theta , x } ) \subseteq (a,b) \backslash \{ \q{\theta}_1 , \q{\theta}_2, \ldots, \q{\theta}_\width \}$. Combining this with \cref{prop:minima:1:eq:defu} and
\cref{prop:minima:1:eq:realization} demonstrates for all $\theta \in U$, $x \in (a,b) \backslash \{ \q{\theta}_1 , \q{\theta}_2, \ldots, \q{\theta}_\width \}$ that $(x - \varepsilon_{\theta , x } , x + \varepsilon_{\theta , x } ) \ni y \mapsto \realization{\theta} ( y ) \in \R$ is affine linear. 
This,
\cref{prop:minima:1:eq:help2}, \cref{prop:minima:1:eq:help4},
and the fact that for all $i \in \N \cap [1, \width )$ it holds that $\alpha_{i+1} \not= \alpha_i$ prove that for all $\theta \in \{ \vartheta \in U \colon \cL ( \vartheta ) = 0\} $,
$i \in \N \cap [1 , \width )$ it holds that $\fx_i \in \{\q{\theta}_1, \q{\theta}_2, \ldots, \q{\theta}_\width \}$. 
Combining this with the fact that
for all $\theta \in U$ it holds that $\q{\theta}_1 \notin [a,b]$,
the fact that for all $\theta \in U$, $j \in \N \cap (1, \width]$ it holds that $\q{\theta}_j \in (a,b)$, the fact that for all $\theta \in U$, $j \in \N \cap (1, \width)$ it holds that $\q{\theta}_j < \q{\theta}_{j+1}$, 
and the fact that $a < \fx_1 < \fx_2 < \cdots < \fx_\width = b$ shows that for all
$\theta \in \{ \vartheta \in U \colon \cL ( \vartheta ) = 0\} $, $j \in \N \cap (1, \width ]$ we have that $ \q{\theta}_j = \fx_{j-1}$. 
This, 
\cref{prop:minima:1:eq:positive},
\cref{prop:minima:1:eq:help2}, \cref{prop:minima:1:eq:help4},
and \cref{prop:minima:1:eq:help3} assure that for all $\theta \in \{ \vartheta \in U \colon \cL ( \vartheta ) = 0\} $,
$x \in [a , b]$ it holds that
\begin{equation}
\begin{split}
     & f( a ) + \alpha_1 (x - a) +  \smallsum_{j=2}^\width ( \alpha_j - \alpha_{j-1} ) \max \{ x - \fx_{j-1}, 0 \} = f(x) \\
     &= \realization{\theta} ( x ) = \c{\theta} + \smallsum_{j = 1}^\width \v{\theta}_j \max \{ \w{\theta}_j x + \b{\theta}_j , 0 \} \\
     &= \c{\theta} + \v{\theta}_1 \max \{\w{\theta}_1 x + \b{\theta}_1 , 0 \} + \smallsum_{j=2}^\width \v{\theta}_j \w{\theta}_j \max \{  x + (\w{\theta}_j)^{-1} \b{\theta}_j , 0 \} \\
     &= \c{\theta} + \v{\theta}_1 ( \w{\theta}_1 x + \b{\theta}_1) + \smallsum_{j=2}^\width \v{\theta}_j \w{\theta}_j \max \{  x - \q{\theta}_j , 0 \} \\
     &= \c{\theta} + \v{\theta}_1  \w{\theta}_1 (x - a ) +  \v{\theta}_1  \w{\theta}_1 a + \v{\theta}_1 \b{\theta}_1 + \smallsum_{j=2}^\width \v{\theta}_j \w{\theta}_j \max \{  x - \fx_{j-1} , 0 \} \\
     &= ( \c{\theta} +  \v{\theta}_1  \w{\theta}_1 a + \v{\theta}_1 \b{\theta}_1 ) 
    + \alpha_1 ( x - a ) + \smallsum_{j=2}^\width \v{\theta}_j \w{\theta}_j \max \{  x - \fx_{j-1} , 0 \} \\
     &= f ( a ) + \alpha_1 ( x - a ) + \smallsum_{j=2}^\width \v{\theta}_j \w{\theta}_j \max \{ x - \fx_{j-1}, 0 \} .  
\end{split}
\end{equation}
Hence, we obtain for all $\theta \in \{ \vartheta \in U \colon \cL ( \vartheta ) = 0\} $,
$j \in \N \cap (1, \width ]$ that $\v{\theta}_j \w{\theta}_j = \alpha_j - \alpha_{j-1}$. Combining this with \cref{prop:minima:1:eq:help3} proves that for all $\theta \in \{ \vartheta \in U \colon \cL ( \vartheta ) = 0\} $ it holds that $\theta \in \cM$.
Hence, we obtain that $\cM = \{ \vartheta \in U \colon \cL ( \vartheta ) = 0\}$. This and the fact that $\cM$ is a non-empty $(\width + 1 )$-dimensional $C^\infty$-submanifold of $\R^\fd$ establish \cref{prop:minima:manifold:1:item3,prop:minima:manifold:1:item4}.

In the next step \nobs that \cref{prop:minima:1:eq:positive} ensures that for all $\theta \in \cM$ it holds that $I_1^\theta = [a,b]$. 
In addition, \nobs that \cref{prop:minima:1:eq:defm} shows that for all $\theta \in \cM$, $j \in \N \cap (1, \width ]$ it holds that $I_j^\theta = (\fx_{j-1}, b ]$.
Furthermore, \nobs that \cref{prop:minima:1:eq:defm} and the fact that for all $j \in \N \cap (1, \width ]$ it holds that $\alpha_j - \alpha_{j-1} \not= 0$ demonstrate that for all $\theta \in \cM$, $i \in \N \cap [1, \width]$ it holds that $\v{\theta}_i \not= 0$.
This, \cref{cor:second:derivatives:globalmin}, and \cref{prop:hessian:minor:det} assure that for all $\theta \in \cM$ it holds that $\det \rbr[\big]{ \rbr[\big]{ \rbr[\big]{ \frac{\partial^2 }{\partial \theta_i \partial \theta_j } \cL} ( \theta ) }_{(i,j) \in \{1, 2, \ldots, 2 \width \} ^2 } } \not= 0$.
Hence, we obtain for all $\theta \in \cM$ that
\begin{equation} \label{prop:minima:1:eq:help5}
\rank ( (\Hs  \cL ) ( \theta ) ) \ge 2 \width.    
\end{equation}
Moreover, \nobs that the fact that $\cM = \{ \vartheta \in U \colon \cL ( \vartheta ) = 0\} $ is a $(\width + 1 )$-dimensional $C^\infty$-submanifold of $\R^\fd$ and \cref{lem:rank:hessian:upperbound} imply that for all $\theta \in \cM$ we have that $\rank ( (\Hs  \cL ) ( \theta ) ) \leq \fd - ( \width + 1 ) = 2 \width$.
This and \cref{prop:minima:1:eq:help5} establish \cref{prop:minima:manifold:1:item5}.
\end{cproof}

\begin{definition} \label{def:minimal:eigenvalue}
Let $n \in \N$ and let $A \in \R^{n \times n} \backslash \{ 0 \} $ be symmetric. Then we denote by $\sigma ( A ) \in (0, \infty)$ the real number given by
\begin{equation}
    \sigma ( A ) = \min \cu*{ \ell \in (0, \infty )  \colon \br*{ \exists \, \lambda \in \cu{- \ell , \ell } , v \in \R^n \backslash \{ 0 \}\colon Av = \lambda v } }
\end{equation}
and we denote by $\Lambda (A)  \in (0, \infty)$ the real number given by
\begin{equation}
    \Lambda (A) = \max \cu*{ \ell \in (0, \infty )  \colon \br*{ \exists \, \lambda \in \cu{- \ell , \ell } , v \in \R^n \backslash \{ 0 \}\colon Av = \lambda v } }
\end{equation}
\end{definition}

\cfclear
\begin{lemma} \label{lem:matrix:norm:estimate}
Let $n \in \N$ and let $A = (a_{i,j})_{(i , j ) \in \cu{1, 2, \ldots, n } ^2 } \in \R^{n \times n} \backslash \{ 0 \} $ be symmetric. Then $\Lambda (A) \leq \br[\big]{ \sum_{i,j=1}^n \abs{a_{i,j}} ^2 }^{1/2}$ \cfadd{def:minimal:eigenvalue}\cfload.
\end{lemma}
\begin{cproof} {lem:matrix:norm:estimate}
Throughout this proof let $\lambda \in \R \backslash \cu{0}$, $v \in \R^n \backslash \cu{0}$ satisfy
\begin{equation} \label{lem:matrix:norm:estimate:eq1}
    A v = \lambda v .
\end{equation}
\Nobs that \cref{lem:matrix:norm:estimate:eq1} ensures that
\begin{equation} \label{lem:matrix:norm:estimate:eq2} \cfadd{def:norm}
    \tfrac{\norm{A v }^2 }{ \norm{v}^2 } = \tfrac{\norm{\lambda v } ^2 }{ \norm{v}^2} = \abs{\lambda}^2
\end{equation}
\cfload.
Moreover, \nobs that the Cauchy-Schwarz inequality demonstrates for all $w = (w_1, \ldots, w_n) \in \R^n$ that
\begin{equation}
    \begin{split}
        \norm{A w } ^2 &= \smallsum_{i=1}^n \abs*{ \sum_{j=1}^n a_{i,j} w_j } ^2 \le \sum_{i=1}^n \br*{ \sum_{j=1}^n \abs{ a_{i,j} w_j } } ^2 \\
        & \le \smallsum_{i=1}^n \br*{ \rbr*{ \smallsum_{j=1}^n \abs{a_{i,j} }^2 } \rbr*{ \smallsum_{j=1}^n \abs{w_j}^2 } } = \norm{w} ^2 \br*{ \smallsum_{i,j=1}^n \abs{a_{i,j}}^2 }.
    \end{split}
\end{equation}
Combining this with \cref{lem:matrix:norm:estimate:eq2} shows that $\abs{\lambda}^2 \le \sum_{i,j=1}^n \abs{a_{i,j}} ^2$.
\end{cproof}

\cfclear
\begin{cor} \label{cor:minima:manifold:1}
Assume \cref{setting:snn}, let $N \in \N \cap [1 , \width ]$, $\fx_0, \fx_1, \ldots, \fx_N, \alpha_1, \alpha_2, \ldots, \alpha_N \in \R$ satisfy $a = \fx_0 < \fx_1 < \cdots < \fx_N = b$, 
assume for all $i \in \{1, 2, \ldots, N \}$, $x \in [\fx_{i-1}, \fx_i ]$ that $f(x) = f( \fx_{i-1} ) + \alpha_i ( x - \fx_{i-1} )$,
and let $\tbound \in \R$ satisfy
\begin{equation}
    \tbound = 1 + \abs{f(a)} + ( 1 + 2 \max\nolimits_{j \in \{1, 2, \ldots, N \} } \abs{\alpha_j} ) ( 1 + \abs{a} + \abs{b} ).
\end{equation}
Then there exist $\mdim \in \N \cap [1, \fd)$ and an open $U \subseteq (- \tbound , \tbound )^\fd$ such that
\begin{enumerate} [label=(\roman*)]
\item \label{cor:minima:manifold:1:item0} it holds that $U \subseteq \fV$,
\item \label{cor:minima:manifold:1:item1} it holds that $\cL |_U \in C^2( U , \R)$,
\item \label{cor:minima:manifold:1:item2} it holds that $U \ni \theta \mapsto (\Hs  \cL) ( \theta ) \in \R^{\fd \times \fd}$ is locally Lipschitz continuous,
\item \label{cor:minima:manifold:1:itemhess} it holds for all $\theta \in U$ that \begin{equation}
    \Lambda ( ( \Hs \cL ) ( \theta ) ) \leq  (3 N + 1 ) \rbr[\big]{24 \tbound ^5 + 16 N \tbound ^7 } \rbr[\big]{ \sup\nolimits_{x \in [a,b]} \dens ( x ) },
\end{equation}
    \item \label{cor:minima:manifold:1:item3} it holds that $\{ \vartheta \in U \colon \cL ( \vartheta ) = 0\} \not= \emptyset$,
    \item \label{cor:minima:manifold:1:item4} it holds that $\{ \vartheta \in U \colon \cL ( \vartheta ) = 0\}$ is a $\mdim$-dimensional $C^\infty$-submanifold of $\R^\fd$,
    \item \label{cor:minima:manifold:1:item5}  it holds for all $\theta \in \{ \vartheta \in U \colon \cL ( \vartheta ) = 0\}$ that $ \rank ( (\Hs  \cL ) ( \theta ) ) = \fd - \mdim$,
    and
    \item \label{cor:minima:manifold:1:item6} it holds that 
        $k = \fd - 2 \br*{ \# \cu{\alpha_1, \alpha_2, \ldots, \alpha_N } }$
\end{enumerate}
\cfadd{def:minimal:eigenvalue}\cfload.
\end{cor}
\begin{cproof}{cor:minima:manifold:1}
Throughout this proof assume without loss of generality that $\prod_{i=1}^{N-1} ( \alpha_{i+1} - \alpha_i ) \not= 0$ (otherwise we can simply remove the points $\fx_i$ which satisfy $\alpha_{i+1} = \alpha_i$ and thereby reduce the number $N$),
let $P \colon \R^{\fd } \to \R^{3N + 1}$ satisfy for all $\theta \in \R^{\fd }$ that $P(\theta) = (\w{\theta}_1, \ldots, \w{\theta}_N, \b{\theta}_1, \ldots, \allowbreak \b{\theta}_N, \allowbreak \v{\theta}_1, \ldots, \v{\theta}_N, \c{\theta} )$,
and let $\scrL \colon \R^{3 N + 1} \to \R$ satisfy for all $\theta = ( \theta_1, \ldots, \theta_{3 N + 1 }) \in \R^{3N+1}$ that
\begin{equation}
    \scrL ( \theta ) = \tint_a^b \rbr[\big]{f ( x ) - \theta_{3 N + 1 } - \smallsum_{j=1}^N \theta_{2N+j} \br{ \Rect ( \theta_j x + \theta_{N+j}) } }^2 \dens ( x ) \, \d x.
\end{equation}
\Nobs that \cref{prop:minima:manifold:1} (applied with $\width \with N$, $\cL \with \scrL$ in the notation of \cref{prop:minima:manifold:1}) demonstrates that there exists an open
$V \subseteq (- \tbound , \tbound )^{3N+1}$ which satisfies that
\begin{enumerate} [label=(\Roman*)]
\item \label{help:defv:item0} it holds that
\begin{equation}
    V \subseteq \cu[\big]{ \theta = (\theta_1, \ldots, \theta_{3 N + 1 } ) \in \R^{3 N + 1 } \colon \rbr[\big]{ \textstyle\prod_{j=1}^N  \textstyle\prod_{v \in \cu{a,b} } (\theta_j v + \theta_{N + j } )  \not= 0 } },
\end{equation}
\item \label{help:defv:item1} it holds that $\scrL |_V \in C^2( V , \R)$,
\item \label{help:defv:item3} it holds for all $\theta = (\theta_1, \ldots, \theta_{3 N + 1 } ) \in V$ that 
\begin{equation} \label{cor:minima:manifold:eq:hess:assumption}
    \max\nolimits_{i,j \in \{1, 2, \ldots, 3 N + 1 \} } \abs[\big]{ \rbr[\big]{ \tfrac{\partial^2}{\partial \theta_i \partial \theta_j} \scrL } ( \theta )} \leq \rbr[\big]{24 \tbound ^5 + 16 N \tbound ^7 } \rbr[\big]{ \sup\nolimits_{x \in [a,b]} \dens ( x ) },
\end{equation}
    \item \label{help:defv:item4} it holds that $\{ \vartheta \in V \colon \scrL ( \vartheta ) = 0\} \not= \emptyset$,
    \item \label{help:defv:item5} it holds that $\{ \vartheta \in V \colon \scrL ( \vartheta ) = 0\}$ is an $(N + 1)$-dimensional $C^\infty$-submanifold of $\R^{3 N + 1}$, and
    \item \label{help:defv:item6} it holds for all $\theta \in \{ \vartheta \in V \colon \scrL ( \vartheta ) = 0\}$ that $ \rank ( (\Hs  \scrL ) ( \theta ) ) = 2 N = (3 N + 1 )  - (N + 1)$.
\end{enumerate} 
In the following let $U \subseteq \R^{\fd }$ satisfy
\begin{equation} \label{cor:minima:manifold:eq:defu}
    U = \cu[\big]{ \theta  \in (- \tbound , \tbound )^{\fd } \cap (  P^{-1} ( V ) ) \colon
    \rbr[\big]{ \forall \, j \in \N \cap (N, \width ] \colon \max \cu[\big]{ \w{\theta}_j a + \b{\theta}_j , \w{\theta}_j b + \b{\theta}_j } < 0 } }.
\end{equation}
\Nobs that \cref{cor:minima:manifold:eq:defu} assures that $U \subseteq \R^\fd$ is open.
In addition, \nobs that \cref{setting:eq:defv}, \cref{cor:minima:manifold:eq:defu}, and \cref{help:defv:item0} imply that $U \subseteq \fV$. This establishes \cref{cor:minima:manifold:1:item0}.
Next \nobs that \cref{cor:minima:manifold:1:item0} and \cref{lem:twice:diff:lipschitz} prove \cref{cor:minima:manifold:1:item1,cor:minima:manifold:1:item2}.
Furthermore, \nobs that for all $\theta \in U$, $x \in [a , b]$, $i \in \N \cap (N, \width ]$ it holds that $\Rect ( \w{\theta}_i x + \b{\theta}_i ) = 0$.
Therefore, we obtain for all $\theta \in U$, $x \in [a,b]$ that
\begin{equation}
    \realization{\theta} (x) = \c{\theta} + \smallsum_{j=1}^\width \v{\theta}_j \br[\big]{\Rect \rbr{ \w{\theta}_{j} x + \b{\theta}_j } } =  \c{\theta} + \smallsum_{j=1}^N \v{\theta}_j \br[\big]{ \Rect \rbr{ \w{\theta}_{j} x + \b{\theta}_j } } .
\end{equation}
This implies for all $\theta \in U$ that 
\begin{equation} \label{cor:minima:manifold:eq:projrisk}
    \cL ( \theta ) = \scrL ( P ( \theta ) ).
\end{equation}
Combining this with \cref{cor:minima:manifold:eq:hess:assumption} ensures for all $\theta \in U$, $i,j \in \N \cap ( (0, N ] \cup (\width , \width + N ] \cup (2 \width , 2 \width + N ] \cup \cu{3 \width + 1 } )$ that 
\begin{equation} \label{cor:minima:manifold:eq:hessest}
\abs[\big]{ \rbr[\big]{ \tfrac{\partial^2}{\partial \theta_i \partial \theta_j} \cL } ( \theta )} \leq \rbr[\big]{24 \tbound ^5 + 16 N \tbound ^7 } \rbr[\big]{ \sup\nolimits_{x \in [a,b]} \dens ( x ) }.    
\end{equation}
Moreover, \nobs that \cref{cor:minima:manifold:eq:projrisk} shows that for all $\theta \in U$, $i \in \{ 1, 2, \ldots, \fd\} \backslash  ( (0, N ] \cup (\width , \width + N ] \cup (2 \width , 2 \width + N ] \cup \cu{3 \width + 1 } )$, $j \in \{1, 2, \ldots, \fd\}$ we have that 
\begin{equation} \label{cor:minima:manifold:eq:hesszero}
    \rbr[\big]{ \tfrac{\partial^2}{\partial \theta_i \partial \theta_j} \cL } ( \theta ) = 0.
\end{equation}
Combining this with \cref{lem:matrix:norm:estimate} and \cref{cor:minima:manifold:eq:hessest} assures for all $\theta \in U$ that
\begin{equation}
    \Lambda ( ( \Hs \cL ) ( \theta ) ) \leq \sqrt{\smallsum_{i,j=1}^\width \abs[\big]{ \rbr[\big]{ \tfrac{\partial^2}{\partial \theta_i \partial \theta_j} \scrL } ( \theta ) }^2} \leq (3 N + 1 ) \rbr[\big]{24 \tbound ^5 + 16 N \tbound ^7 } \rbr[\big]{ \sup\nolimits_{x \in [a,b]} \dens ( x ) } .
\end{equation}
This establishes \cref{cor:minima:manifold:1:itemhess}.
Furthermore, \nobs that \cref{help:defv:item4,help:defv:item5},
\cref{cor:minima:manifold:eq:defu},
and \cref{cor:minima:manifold:eq:projrisk}
establish \cref{cor:minima:manifold:1:item3,cor:minima:manifold:1:item4,cor:minima:manifold:1:item6}.
In addition, 
\nobs that \cref{cor:minima:manifold:eq:projrisk}, \cref{cor:minima:manifold:eq:hesszero}, and \cref{help:defv:item6} demonstrate for all $\theta \in \cu{\vartheta \in U \colon \cL ( \vartheta ) = 0 }$ that $\rank ( ( \Hs  \cL ) ( \theta ) ) = 2N$.
 Combining this with \cref{cor:minima:manifold:1:item6} establishes \cref{cor:minima:manifold:1:item5}.
\end{cproof}

\section{Local convergence to the set of global minima for gradient flow (GF)} \label{section:flow:convergence}

In this section we employ \cref{cor:minima:manifold:1} from \cref{section:risk:global:min} 
to establish in \cref{theo:flow:convergence} in \cref{subsection:flow:convergence} below 
and \cref{cor:flow:convergence:random} in \cref{subsection:flow:random} below that 
the risk of certain solutions of GF differential equations converges 
under the assumption that the target function is piecewise 
constant exponentially quick to zero. Our proof of \cref{theo:flow:convergence} 
employs the abstract local convergence result for GF trajectories 
in \cref{prop:gradient:flow:abstract} in \cref{subsection:flow:abstract}. 
\cref{prop:gradient:flow:abstract} and its proof are strongly 
inspired by Fehrman et al.~\cite[Proposition 16]{FehrmanGessJentzen2020}.
Our proofs of \cref{prop:gradient:flow:abstract,theo:flow:convergence} also 
use the several well-known concepts and results from differential geometry 
which we recall in \cref{subsection:proj:submanifold} below. 

In particular,
\cref{lem:unique:projection} is a direct consequence of, e.g.,~\cite[Proposition 7]{FehrmanGessJentzen2020}, \cref{lem:projection:perp} is proved as, e.g.,~\cite[Lemma 10]{FehrmanGessJentzen2020}, \cref{lem:gradient:distance} is proved as, e.g.,~\cite[Lemma 11]{FehrmanGessJentzen2020}, \cref{def:local:tubular:nbhd} is a slight reformulation of, e.g.,~\cite[Definition 12]{FehrmanGessJentzen2020}, \cref{prop:tubular:nbhd} is a slight extension of, e.g.,~\cite[Proposition 13]{FehrmanGessJentzen2020},
\cref{prop:hessian:decomp} is a reformulation of \cite[Lemma 15]{FehrmanGessJentzen2020},
and \cref{lem:second:order:approx} is a slight generalization of \cite[Lemma 14]{FehrmanGessJentzen2020}.
 
 \subsection{Differential geometric preliminaries} \label{subsection:proj:submanifold}

 \cfclear
\begin{definition} \label{def:dist:set}
Let $\fd \in \N$
and let $\cM \subseteq \R^\fd$ satisfy $\cM \not= \emptyset$. Then we denote by $\dist _\cM \colon \R^\fd \to \R$ the function which satisfies for all $x \in \R^\fd$ that $\dist _\cM ( x ) = \inf_{y \in \cM } \norm{x-y}$ \cfadd{def:norm}\cfload.
\end{definition}

 \cfclear
\begin{definition}\label{def:subset:unique:projection}
Let $\fd \in \N$
and let $\cM \subseteq \R^\fd$ satisfy $\cM \not= \emptyset$.
Then we denote by $\scrP_{\cM} \subseteq \R^{\fd}$ the set given by
\begin{equation} \cfadd{def:norm}
\scrP_{\cM} = \cu*{ x \in \R^{\fd} \colon (\exists_{1} \, y \in \cM \colon \norm{x - y} = \dist _\cM ( x ) ) }
\end{equation}
and we denote by $\projm _\cM \colon \scrP_{\cM} \to \R^{\fd}$ the function which satisfies for all $x \in \scrP_{\cM}$ that $\projm _\cM ( x ) \in \cM$ and
\begin{equation} \cfadd{def:dist:set}
\norm{x - \projm _\cM(x)} = \dist_\cM ( x )
\end{equation}
\cfload.
\end{definition}

\cfclear
\begin{definition}\label{def:unique:projection}
 Let $\fd \in \N$ and let $\cM \subseteq \R^\fd$ satisfy $\cM \not= \emptyset$. Then we denote by $\tubn _\cM \subseteq \R^{\fd}$ the set given by
\begin{equation} \cfadd{def:subset:unique:projection}
\tubn _\cM = \bigcup\nolimits_{\substack{U \subseteq \R^{\fd} \; \text{is open}, \; U \subseteq \scrP_{\cM}, \\ \text{and} \; \projm _\cM|_{U} \in C^1(U, \R^{\fd})}} U
\end{equation}
\cfload.
\end{definition}

\cfclear
\begin{lemma} \label{lem:unique:projection}
Let $\fd , \mdim \in \N$, 
let $\cM \subseteq \R^\fd$ be a $\mdim$-dimensional $C^2$-submanifold of $\R^\fd$, and let $x \in \cM$. 
Then there exists an open $V \subseteq \R^\fd$ such that
\begin{enumerate} [label = (\roman*)] \cfadd{def:dist:set} \cfadd{def:norm} \cfadd{def:subset:unique:projection}
    \item it holds that $x \in V \subseteq \scrP_\cM$ and
    \item it holds that $\projm _\cM |_V  \in C^1 ( V , \R^\fd)$.
\end{enumerate}
\cfload.
\end{lemma}

\cfclear
\begin{prop} \label{prop:unique:projection}
Let $\fd , \mdim \in \N$ and let $\cM \subseteq \R^\fd$ be a non-empty $\mdim$-dimensional $C^2$-submanifold of $\R^\fd$. Then $\cM \subseteq \tubn_\cM  $
\cfadd{def:unique:projection}\cfload.
\end{prop}
\begin{cproof}{prop:unique:projection}
\Nobs that \cref{lem:unique:projection} assures that $\cM \subseteq \tubn_\cM$.
\end{cproof}

\cfclear
\begin{lemma} \label{lem:projection:perp}
Let $\fd , \mdim \in \N$,
 let $\cM \subseteq \R^\fd$ be a non-empty $\mdim$-dimensional $C^2$-submanifold of $\R^\fd$,
 and let $x \in \tubn_\cM$ \cfadd{def:unique:projection}\cfload.
Then $x - \projm _\cM ( x ) \in ( \cT_{\cM}^{\projm _\cM (x) }) ^\perp$ 
\cfadd{def:subset:unique:projection}\cfadd{def:unique:projection}\cfadd{def:tangent:space}\cfload.
\end{lemma}

\cfclear
\begin{lemma} \label{lem:gradient:distance}
Let $\fd , \mdim \in \N$
and let $\cM \subseteq \R^\fd$ be a non-empty $\mdim$-dimensional $C^2$-submanifold of $\R^\fd$. Then 
\begin{enumerate} [label = (\roman*)] \cfadd{def:norm}
\item it holds that $\tubn_\cM \backslash \cM \subseteq \R^\fd$ is open,
    \item it holds that $\tubn_\cM \backslash \cM \ni y \mapsto \dist _\cM ( y  ) \in \R$ is continuously differentiable, and
    \item it holds for all $y \in \tubn_\cM \backslash \cM$ that
\begin{equation}
    (\nabla \dist _\cM) ( y ) = \tfrac{y - \projm_\cM ( y ) }{\norm{y - \projm_\cM ( y ) } }
\end{equation}
\end{enumerate}
\cfadd{def:dist:set}\cfadd{def:subset:unique:projection}\cfadd{def:unique:projection}\cfload.
\end{lemma}

\cfclear
\begin{definition} \label{def:local:tubular:nbhd}
Let $\fd , \mdim \in \N$, 
let $\cM \subseteq \R^\fd$ be a $\mdim$-dimensional $C^2$-submanifold of $\R^\fd$,
and let $x \in \cM $, $r , s \in (0, \infty)$. Then we denote by $V^{r,s}_{\cM , x }  \subseteq \R^\fd$ the set given by
\begin{equation}  \label{eq:def:local:tubular:nbhd} \cfadd{def:norm} \cfadd{def:tangent:space}
    V^{r , s}_{\cM , x }  = \cu*{ y \in \R^\fd \colon \exists \, \fm \in \cM \colon \exists \, v \in (\cT_\cM^{\fm})^\perp \colon \br[\big]{ ( \norm{\fm - x } \le r ) , \, ( \norm{v} < s ) , \, ( y = \fm + v ) } }
\end{equation}
\cfload.
\end{definition}

\cfclear
\begin{lemma} \label{lem:local:tubular:nbhd}
	Let $\fd , \mdim \in \N$, 
	let $\cM \subseteq \R^\fd$ be a $\mdim$-dimensional $C^2$-submanifold of $\R^\fd$,
	and let $x \in \cM $, $r , s \in (0, \infty)$.
	Then
	\begin{enumerate} [label = (\roman*)]
		\item \label{lem:local:tubular:nbhd:item1} it holds that
		\begin{equation} \cfadd{def:norm} \cfadd{def:tangent:space}
		V^{r , s}_{\cM , x }  = \cu*{ y \in \R^\fd \colon \exists \, \fm \in \cM \colon  \br[\big]{ ( \norm{\fm - x } \le r ) , \, ( \norm{y - \fm } < s ) , \, ( y - \fm \in (\cT_\cM^\fm)^\perp ) } },
		\end{equation}
	\item \label{lem:local:tubular:nbhd:item2} it holds that
	\begin{equation} \cfadd{def:unique:projection} \cfadd{def:local:tubular:nbhd} \cfadd{def:subset:unique:projection}
			V^{r , s}_{\cM , x }  \supseteq \cu*{y \in \tubn_\cM \colon \br[\big]{( \norm{x - \projm_\cM ( y ) } \le r ) , \, ( \norm{y - \projm_\cM ( y ) } < s)}},
	\end{equation}
	and 
		\item \label{lem:local:tubular:nbhd:item3} it holds that $x \in (V^{r,s}_{\cM , x })^\circ$
		\cfload.
	\end{enumerate}
\end{lemma}
\begin{cproof}{lem:local:tubular:nbhd}
	\Nobs that \cref{eq:def:local:tubular:nbhd} establishes \cref{lem:local:tubular:nbhd:item1}. Next \nobs that \cref{eq:def:local:tubular:nbhd,lem:projection:perp} establish \cref{lem:local:tubular:nbhd:item2}.
	Furthermore, \nobs that \cref{lem:local:tubular:nbhd:item2} implies that
	\begin{equation} \label{lem:local:tubular:nbhd:eq1}
		V^{r , s}_{\cM , x }  \supseteq \cu*{y \in \tubn_\cM \colon \br[\big]{( \norm{x - \projm_\cM ( y ) } < r ) , \, ( \norm{y - \projm_\cM ( y ) } < s)}}.
	\end{equation}
	Furthermore, \nobs that the fact that $\tubn_\cM \ni y \mapsto \projm ( y ) \in \R^\fd$ is continuous shows that $\cu{y \in \tubn_\cM \colon \br{( \norm{x - \projm_\cM ( y ) } < r ) , \, ( \norm{y - \projm_\cM ( y ) } < s)}} \subseteq \R^\fd$ is open. Combining this with \cref{lem:local:tubular:nbhd:eq1} and the fact that $x \in \cu{y \in \tubn_\cM \colon \br{( \norm{x - \projm_\cM ( y ) } < r ) , \, ( \norm{y - \projm_\cM ( y ) } < s)}}$ establishes \cref{lem:local:tubular:nbhd:item3}.
\end{cproof}

\cfclear
\begin{prop} \label{prop:tubular:nbhd}
Let $\fd , \mdim \in \N$,
let $\cM \subseteq \R^\fd$ be a $\mdim$-dimensional $C^2$-submanifold of $\R^\fd$,
let $U \subseteq \tubn_\cM$ be open,
and let $x \in \cM \cap U$
\cfadd{def:unique:projection}\cfload.
Then there exist $R , S \in (0, \infty)$ such that 
\begin{enumerate} [label = (\roman*)] \cfadd{def:local:tubular:nbhd}  \cfadd{def:norm}
    \item \label{prop:tubular:nbhd:item1} it holds for all $r \in (0, R ]$, $s \in (0, S]$ that $\overline{V^{r , s}_{ \cM , x } }  \subseteq U$,
    \item \label{prop:tubular:nbhd:item2} it holds for all $r \in (0, R ]$, $s \in (0, S]$ that 
    \begin{equation}
        V^{r , s}_{\cM , x }  = \cu*{ y \in \R^\fd \colon \dist _\cM ( y ) = \dist _{\cu{\fm \in \cM \colon \norm{x - \fm} \le r } } ( y ) < s},
    \end{equation}
    \item \label{prop:tubular:nbhd:item3} it holds for all $r \in (0, R ]$, $s \in (0, S]$, $\fm \in \cM $, $v \in ( \cT_\cM^\fm )^\perp$ with $\norm{\fm - x } \le r$ and $\norm{v} < s$ that $\fm + v \in V^{r , s}_{\cM , x } $ and $\projm_\cM ( \fm + v ) = \fm$,
    and
    \item \label{prop:tubular:nbhd:item4} it holds for all $r \in (0, R]$, $s \in (0, S]$ that
    \begin{equation}
    V^{r , s}_{\cM , x } = \cu*{y \in \tubn_\cM \colon \br[\big]{( \norm{x - \projm_\cM ( y ) } \le r ) , \, ( \norm{y - \projm_\cM ( y ) } < s)}}
    \end{equation}
\end{enumerate}
\cfadd{def:subset:unique:projection}\cfadd{def:dist:set}\cfadd{def:tangent:space}\cfout.
\end{prop}
\begin{cproof}{prop:tubular:nbhd}
	\Nobs that \cite[Proposition 13]{FehrmanGessJentzen2020} establishes \cref{prop:tubular:nbhd:item1,prop:tubular:nbhd:item2,prop:tubular:nbhd:item3}.
	In addition, \nobs that \cref{prop:tubular:nbhd:item2,prop:tubular:nbhd:item3} and \cref{eq:def:local:tubular:nbhd} establish \cref{prop:tubular:nbhd:item4}.
\end{cproof}

 \begin{setting} \label{setting:mfld}
 Let $\fd \in \N$, $\mdim  \in \N \cap (0, \fd)$, let $U \subseteq \R^\fd$ be open, let $f \in C^2 ( U , \R) $ have locally Lipschitz continuous derivatives, let $\cM \subseteq U$ satisfy $\cM = \{ x \in U \colon f(x) = \inf_{y \in U} f(y) \}$,
 and assume that $\cM $ is a $\mdim$-dimensional $C^2$-submanifold of $\R^\fd$.
 \end{setting}

\cfclear
\begin{prop} \label{prop:hessian:decomp}
Assume \cref{setting:mfld} and let $x \in \cM $ satisfy $\rank ( (\Hs  f ) ( x ) ) = \fd - \mdim$. Then
\begin{enumerate} [label = (\roman*)] \cfadd{def:minimal:eigenvalue} \cfadd{def:norm} \cfadd{def:tangent:space}
    \item \label{prop:hessian:decomp:item2}  it holds for all $v \in ( ( \cT_\cM^x ) ^\perp ) \backslash \{0 \}$ that $\spro{ ((\Hs  f) (x) ) v, v } \geq \br{ \sigma ( ( \Hs  f ) ( x ) ) } \norm{v}^2 > 0$
    and
    \item \label{prop:hessian:decomp:item3} it holds for all $v \in ( ( \cT_\cM^x ) ^\perp ) \backslash \{0 \}$, $r \in [0, \rbr{\Lambda (( \Hs  f ) ( x )) }^{-1} ]$ that $\norm{v - r ( ( \Hs  f)(x ) ) v } \leq [ 1 - r \sigma ( ( \Hs  f ) ( x ) ) ] \norm{v}$.
\end{enumerate}
\cfload.
\end{prop}
\begin{cproof}{prop:hessian:decomp}
Throughout this proof let $\{v_1, v_2, \ldots, v_{\fd- \mdim} \} \subseteq ( ( \cT_\cM^x ) ^\perp ) \backslash \{ 0 \}$ be an orthogonal basis of $( \cT_\cM^x ) ^\perp$ with respect to which $(\Hs f ) ( x )$ is diagonal and let $\lambda_1, \lambda_2, \ldots, \lambda_{\fd -  \mdim} \in \R$ satisfy for all $i \in \{1, 2, \ldots, \fd -  \mdim \}$ that $(( \Hs f) ( x ) ) v_i = \lambda_i v_i$. \Nobs that the fact that $x$ is a local minimum of $f$ shows for all $i \in \{1, 2, \ldots, \fd -  \mdim\}$ that $\lambda_i \geq 0$. This and the assumption that $\rank ( (\Hs  f)(x ) ) = \fd -  \mdim$ imply for all $i \in \{1, 2, \ldots, \fd -  \mdim\}$ that $\lambda_i > 0$. Hence, we obtain for all $i \in \{1, 2, \ldots, \fd -  \mdim\}$ that $\lambda_i \in [\sigma ( ( \Hs f ) ( x ) ), \Lambda ( ( \Hs f ) ( x ) ) ]$.
Next let $\bfv \in ( ( \cT_\cM^x ) ^\perp ) \backslash \{ 0 \}$ and let $u_1, u_2, \ldots, u_{\fd -  \mdim} \in \R$ satisfy $\bfv = \sum_{i=1}^{\fd -  \mdim} u_i v_i$. \Nobs that
\begin{equation}
\begin{split}
    \spro*{ \rbr[\big]{  ( \Hs  f) (x ) } \bfv, \bfv } &= \smallsum_{i=1}^{\fd -  \mdim} \rbr[\big]{ \lambda_i \abs{u_i}^2 \norm{v_i}^2 }
    \geq \br[\big]{ \sigma ( ( \Hs f ) ( x ) ) } \br[\big]{\smallsum_{i=1}^{\fd -  \mdim} \abs{u_i}^2 \norm{v_i}^2 } \\
    &= \br[\big]{ \sigma ( ( \Hs f ) ( x ) ) } \norm{\bfv}^2 > 0.
    \end{split}
\end{equation}
This establishes \cref{prop:hessian:decomp:item2}.
Furthermore, \nobs that the fact that for all $i \in \{1, 2, \ldots, \fd -  \mdim\}$ it holds that $\lambda_i \in [\sigma ( ( \Hs f ) ( x ) ), \Lambda ( ( \Hs f ) ( x ) ) ]$ ensures that for all $r \in [0, (\Lambda ( \Hs f ) ( x ) )^{-1}]$ we have that
\begin{equation}
    \begin{split}
        \norm*{\bfv - r \rbr[\big]{ ( \Hs  f ) ( x ) } \bfv }^2
        &= \smallsum_{i=1}^{\fd -  \mdim}\rbr[\big]{ \abs{u_i}^2 \norm{v_i}^2 ( 1 - r \lambda_i )^2 } \\
        & \leq \smallsum_{i=1}^{\fd -  \mdim}\rbr[\big]{ \abs{u_i}^2 \norm{v_i}^2 \rbr[\big]{ 1 - r \br[\big]{\sigma ( ( \Hs f ) ( x ) ) } }^2 } \\
        &= \rbr[\big]{ 1 - r \br[\big]{ \sigma ( ( \Hs  f)(x) ) } }  ^2 \norm{\bfv}^2.
    \end{split}
\end{equation}
This establishes \cref{prop:hessian:decomp:item3}.
\end{cproof}

\cfclear
\begin{lemma} \label{lem:second:order:approx}
Assume \cref{setting:mfld} and let $x \in \cM $.
Then there exist $c, r, s \in (0, \infty)$ such that for all $y \in V^{r , s}_{\cM , x }$ it holds that $\overline{V^{r , s}_{\cM , x }} \subseteq ( \tubn_\cM \cap U ) $ and
\begin{equation} \cfadd{def:norm} \cfadd{def:dist:set}\cfadd{def:subset:unique:projection}\cfadd{def:unique:projection}\cfadd{def:local:tubular:nbhd}
    \norm*{ (\nabla f) ( y ) - \rbr[\big]{ ( \Hs  f ) ( \projm_\cM ( y ) ) } ( y - \projm_\cM ( y ) )} \leq c ( \dist _\cM ( y  ) )^2
\end{equation}
\cfout.
\end{lemma}
\begin{cproof}{lem:second:order:approx}
\Nobs that \cref{prop:tubular:nbhd} ensures that there exist $r , s \in (0, \infty)$ which satisfy $\overline{V^{r,s}_{\cM , x } } \subseteq  U$, which satisfy
    \begin{equation} \label{lem:second:order:eq1}
          V^{r , s}_{\cM , x } = \cu*{y \in \tubn_\cM \colon \br[\big]{( \norm{x - \projm_\cM ( y ) } \le r ) , \, ( \norm{y - \projm_\cM ( y ) } < s)}},
    \end{equation}
    and which satisfy for all $\fm \in \cM $, $v \in ( \cT_\cM^\fm )^\perp$ with $\norm{\fm - x } \le r$ and $\norm{v} < s$ that $\fm + v \in  V^{r , s}_{\cM , x }$ and
\begin{equation} \label{lem:second:order:eq2} \cfadd{def:tangent:space}
	\projm_\cM ( \fm + v ) = \fm
\end{equation} 
\cfload.
\Nobs that \cref{lem:second:order:eq1},
\cref{lem:second:order:eq2},
and \cref{lem:projection:perp} imply for all $y \in V^{r,s}_{\cM , x }$, $t \in [0,1]$ that $\projm_\cM ( y ) + t ( y - \projm_\cM  ( y ) ) \in V^{r,s} _{\cM , x }$.
In addition,
\nobs that the fact that $\overline{V^{r , s}_{\cM , x }} $ is compact and the assumption that $U \ni y \mapsto ( \Hs f ) ( y ) \in \R^{\fd \times \fd}$ is locally Lipschitz continuous prove that there exists $c \in (0, \infty)$ which satisfies for all $y,z \in \overline{V^{r , s}_{\cM , x }} $, $v \in \R^\fd$ that $\norm{((\Hs  f)(y) - (\Hs  f)(z)) v } \leq c \norm{y-z} \norm{v}$.
 Furthermore, \nobs that the fact that for all $y \in V^{r , s}_{\cM , x }$ it holds that $(\nabla f) ( \projm_\cM ( y ) ) = 0$ and the assumption that $f $ is twice continuously differentiable demonstrate that for all $y \in V^{r , s}_{\cM , x }$ it holds that 
\begin{equation}
\begin{split}
    (\nabla f) ( y ) &= \int_0^1 \rbr[\big]{(\Hs  f) ( \projm _\cM ( y ) + t (y - \projm _\cM ( y ) ) } (y - \projm _\cM ( y ) ) \, \d t \\
    &= \rbr[\big]{ (\Hs  f) ( \projm _\cM ( y ) ) } ( y - \projm _\cM ( y ) ) \\
    & \quad + \int_0^1 \rbr[\Big]{(\Hs  f) \rbr[\big]{ \projm _\cM ( y ) + t (y - \projm _\cM ( y ) ) } - (\Hs  f) ( \projm _\cM ( y ) ) } (y - \projm _\cM ( y ) ) \, \d t.
    \end{split}
\end{equation}
Combining this with the fact that for all $y \in V^{r , s}_{\cM , x }$, $t \in [0,1]$ it holds that
\begin{equation}
    \norm*{ \rbr[\big]{ (\Hs  f) ( \projm _\cM ( y ) + t (y - \projm _\cM ( y ) ) - (\Hs  f) ( \projm _\cM ( y ) ) } ( y - \projm _\cM ( y ) ) } \leq c t \norm{y - \projm _\cM ( y ) } ^2
\end{equation}
implies that for all $y \in V^{r , s}_{\cM , x }$ we have that
\begin{equation}
    \begin{split}
         \norm*{ (\nabla f) ( y ) - ( ( \Hs  f ) ( \projm _\cM ( y ) ) ) ( y - \projm _\cM ( y ) )} \leq c \norm{y - \projm _\cM ( y ) } ^2 \br*{ \tint_0^1 t \, \d t } = \tfrac{c}{2} ( \dist _\cM ( y  ) )^2.
    \end{split}
\end{equation}
\end{cproof}

\subsection{Abstract convergence result for GF to a submanifold of global minima} \label{subsection:flow:abstract}

\cfclear
\begin{prop} \label{prop:gradient:flow:abstract}
Assume \cref{setting:mfld}, assume for all $x \in \cM $ that $\rank ( (\Hs  f) (x) ) = \fd - \mdim$, let $\cG \colon \R^\fd \to \R^\fd$ be locally bounded and measurable, assume for all $x \in U$ that $\cG ( x ) = (\nabla f ) ( x )$, let $\Theta^\theta \in C([0, \infty), \R^\fd)$, $\theta \in \R^\fd$, satisfy for all $\theta \in \R^\fd$, $t \in [0, \infty)$ that $\Theta_t^\theta = \theta - \int_0^t \cG ( \Theta_s^\theta ) \, \d s $, and let $x \in \cM $. Then there exist $r , s \in (0, \infty)$ such that
\begin{enumerate} [label = (\roman*)]
\cfadd{def:dist:set}
\cfadd{def:minimal:eigenvalue} \cfadd{def:local:tubular:nbhd}
    \item it holds for all $\theta \in V^{r /2 , s}_{\cM , x }$, $t \in [0, \infty)$ that $\Theta_t^\theta \in V^{r , s}_{\cM , x }$,
    \item it holds that $ \inf_{y \in \cM  \cap V^{r,s}_{\cM , x }} \br*{ \sigma ((\Hs  f) ( y ) ) } > 0$, and
    \item it holds for all $\theta \in V^{r/2 , s}_{\cM , x }$, $t \in [0, \infty)$ that 
    \begin{equation}
        \dist  _{\cM  } ( \Theta_t^\theta ) \leq \exp\rbr[\big]{- \tfrac{t}{2} \br[\big]{ \inf\nolimits_{y \in \cM  \cap V^{r,s}_{\cM , x }} \br[\big]{ \sigma ((\Hs  f) ( y ) ) } } } \dist _\cM ( \theta )
    \end{equation}
\end{enumerate}
\cfload.
\end{prop}

\cfclear
\begin{cproof} {prop:gradient:flow:abstract}
\Nobs that \cref{prop:tubular:nbhd,lem:second:order:approx} prove 
that there exist $r, \varepsilon, \fc \in (0, \infty)$ 
which satisfy $\overline{V^{r , \varepsilon }_{\cM , x } } \subseteq U$,
which satisfy
    \begin{equation} \cfadd{def:unique:projection}
 V^{r , s}_{\cM , x } = \cu*{y \in \tubn_\cM \colon \br[\big]{( \norm{x - \projm_\cM ( y ) } \le r ) , \, ( \norm{y - \projm_\cM ( y ) } < s)}},
    \end{equation}
    and which satisfy for all $y \in \overline{V^{r , \varepsilon }_{\cM , x } } $ that
    \begin{equation} \label{prop:gf:eq:hess}
    \norm*{ (\nabla f) ( y ) - (\Hs  f) ( \projm _\cM ( y ) ) ( y - \projm _\cM ( y ) )} \leq \fc ( \dist _\cM ( y ) )^2
\end{equation}
    \cfload.
    In the following let $\kappa \in \R$ satisfy $\kappa = \frac{1}{2} \inf_{y \in \cM  \cap V^{r , \varepsilon }_{\cM , x }} \br[\big]{ \sigma ((\Hs  f) ( y ) ) }$. \Nobs that the fact that $\Hs  f$ is locally Lipschitz continuous and the fact that the eigenvalues are continuous functions of a matrix 
    (cf., e.g., Kato~\cite[Section 2.5.1]{Kato1966}) prove that $\kappa > 0$.
    Next \nobs that the fact that $\overline{V^{r , \varepsilon }_{\cM , x } } $ is compact, the fact that for all $y \in \tubn_\cM$ it holds that $ (\nabla f ) ( \projm _\cM (y ) ) = 0$, the fact that $\tubn_\cM \ni y \mapsto \projm _\cM ( y ) \in \R^\fd $ is continuously differentiable, and the assumption that $f \in C^2 ( U , \R)$ prove that there exists $c \in (0, \infty)$ which satisfies for all $y \in \overline{V^{r , \varepsilon }_{\cM , x } } $ that
    \begin{equation} \label{prop:gf:eq:projection} \cfadd{def:norm} \cfadd{def:subset:unique:projection}
        \norm{( \projm _\cM ) ' ( y ) [(\nabla f)(y)] } = \norm{ ( \projm _\cM ) ' ( y ) [(\nabla f)(y) - (\nabla f ) ( \projm _\cM ( y ) ) ] } \leq c \norm{y - \projm _\cM ( y ) } = c \dist  _\cM ( y )
    \end{equation}
    \cfload.
In the following let $s \in (0, \infty)$ satisfy
\begin{equation} \label{prop:gf:eq:defdelta}
    s  = \min \cu* { \frac{\kappa}{ \fc }, \frac{\kappa r}{2 c}, \varepsilon},
\end{equation}
let $\theta \in V^{r /2 , s}_{\cM , x }$,
 and let $\tau \in (0, \infty]$ satisfy $\tau = \inf \rbr{ \cu{t \in [0, \infty) \colon \Theta_t^\theta \notin V^{r  , s}_{\cM , x } }\cup \cu{\infty}}$. \Nobs that the assumption that for all $y \in U$ it holds that $ \cG ( y ) = ( \nabla f ) ( y )$ and the fact that $U \ni y \mapsto (\nabla f ) ( y ) \in \R^\fd$ is continuous assure that $[0, \tau) \ni t \mapsto \Theta_t^\theta \in \R^\fd$ is continuously differentiable and that for all $t \in [0, \tau)$ it holds that $\frac{\d}{\d t} \Theta_t^\theta = - ( \nabla f ) ( \Theta_t ^\theta )$. This, \cref{lem:gradient:distance}, and the chain rule show for all $t \in [0, \tau)$ that
\begin{equation} \label{prop:gf:eq:dist1} \cfclear \cfadd{def:norm}
    \frac{\d}{\d t} \dist  _\cM ( \Theta_t^\theta ) = - \spro*{ (\nabla f ) ( \Theta_t^\theta), (\nabla \dist _\cM ) ( \Theta_t^\theta ) } = - \spro*{ (\nabla f ) ( \Theta_t^\theta), \frac{\Theta_t^\theta - \projm _\cM ( \Theta_t^\theta)}{ \norm{\Theta_t^\theta - \projm _\cM ( \Theta_t^\theta)}} }
\end{equation}
\cfload.
Next \nobs that \cref{prop:gf:eq:hess}, \cref{prop:gf:eq:defdelta},
\cref{prop:gf:eq:dist1}, and \cref{prop:hessian:decomp} demonstrate for all $t \in [0, \tau)$ that
\begin{equation}
    \begin{split}
    \frac{\d}{\d t} \dist _\cM ( \Theta_t^\theta ) 
    &= - \spro*{ (\Hs f ) ( \projm _\cM (\Theta_t^\theta ) ) (\Theta_t^\theta - \projm _\cM ( \Theta_t^\theta ) ) , \frac{\Theta_t^\theta - \projm _\cM ( \Theta_t^\theta)}{ \norm{\Theta_t^\theta - \projm _\cM ( \Theta_t^\theta)}} } \\
    & \quad -  \spro*{ (\nabla f ) ( \Theta_t^\theta) -  (\Hs f ) ( \projm _\cM (\Theta_t^\theta ) ) (\Theta_t^\theta - \projm _\cM ( \Theta_t^\theta ) ) , \frac{\Theta_t^\theta - \projm _\cM ( \Theta_t^\theta)}{ \norm{\Theta_t^\theta - \projm _\cM ( \Theta_t^\theta)}} } \\
    &\leq - 2 \kappa \norm{\Theta_t^\theta - \projm _\cM ( \Theta_t^\theta)} + \fc ( \dist _\cM ( \Theta_t^\theta  ) )^2 \\
    &= - 2 \kappa  \dist _\cM ( \Theta_t^\theta  ) + \fc ( \dist _\cM ( \Theta_t^\theta ) )^2 \leq - \kappa  \dist _\cM ( \Theta_t^\theta ).
    \end{split}
\end{equation}
Hence, we obtain for all $t \in [0, \tau)$ that
\begin{equation} \label{prop:gf:eq:distance}
    \dist _\cM ( \Theta_t^\theta ) \leq e^{-\kappa t} \dist _\cM ( \Theta_0^\theta ) = e^{-\kappa t} \dist  _\cM ( \theta ).
\end{equation}
It remains to prove that $\tau = \infty$. To this end, \nobs that the chain rule and \cref{lem:unique:projection} imply for all $t \in [0, \tau)$ that
\begin{equation}
    \tfrac{\d }{\d t} \projm _\cM ( \Theta_t^\theta ) = - (D \projm _\cM ) ( \Theta_t^\theta) (\nabla f ( \Theta_t^\theta ) ) .
\end{equation}
Combining this, \cref{prop:gf:eq:projection}, and \cref{prop:gf:eq:distance} ensures for all $t \in [0, \tau)$ that
\begin{equation}
    \norm*{ \tfrac{\d }{\d t} \projm _\cM ( \Theta_t^\theta )} \leq c \dist _\cM ( \Theta_t^\theta ) \leq c e^{-\kappa t} \dist _\cM ( \theta ) \leq c s e^{-\kappa t}.
\end{equation}
This and \cref{prop:gf:eq:defdelta} show for all $t \in [0, \tau)$ that
\begin{equation} \label{prop:gf:eq:dist:proj}
    \norm{\projm _\cM ( \Theta_t^\theta) - \projm _\cM ( \theta ) } \leq c s \int_0^t e^{-\kappa u} \, \d u \leq \frac{\kappa r}{2} \int_0^\infty e^{-\kappa u} \, \d u = \frac{r}{2}.
\end{equation}
Furthermore, \nobs that the assumption that $\theta \in V^{r /2 , s}_{\cM , x }$ assures that there exists $\delta \in (0, \infty)$ which satisfies that $\theta \in V^{r /2 - \delta , s}_{\cM , x } $. Combining this with \cref{prop:gf:eq:dist:proj} establishes for all $t \in [0, \tau)$ that $\Theta_t^\theta \in V^{r - \delta , s}_{\cM , x }$. Consequently, we must
have that $\tau = \infty$.
\end{cproof}

\subsection{Convergence rates for GF in the training of ANNs} \label{subsection:flow:convergence}

\begin{lemma} \label{lem:gradient:measurable}
	Assume \cref{setting:snn}.
	Then $\cG$ is locally bounded and measurable.
\end{lemma}
\begin{cproof}{lem:gradient:measurable}
	\Nobs that, e.g., \cite[Corollary 2.4]{JentzenRiekertFlow}
	demonstrates that $\cG$ is locally bounded and measurable.
\end{cproof}

\cfclear
\begin{prop} \label{theo:flow:convergence}
Assume \cref{setting:snn}, let $N \in \N \cap [1 , \width ]$,
$\fx_0, \fx_1, \ldots, \fx_N, \alpha_1, \alpha_2, \ldots, \alpha_N \in \R$ satisfy $a = \fx_0 < \fx_1 < \cdots < \fx_N = b$, assume for all $i \in \{1, 2, \ldots, N\}$, $x \in [\fx_{i-1}, \fx_i ]$ that $f(x) = f(\fx_{i-1} ) + \alpha_i ( x - \fx_{i-1} )$,
and let $\Theta^\theta \in C([0, \infty), \R^{\fd })$, $\theta \in \R^{\fd }$, satisfy for all $\theta \in \R^{\fd }$, $t \in [0, \infty)$ that 
\begin{equation} \cfadd{lem:gradient:measurable}
\Theta_t^\theta = \theta - \int_0^t \cG ( \Theta_s^\theta) \, \d s     
\end{equation}
\cfload.
Then there exist $\fc , \fC \in (0, \infty)$ 
and a non-empty open $U \subseteq \R^{\fd }$ such that for all $\theta \in U$, $t \in [0, \infty)$ it holds that $\cL ( \Theta_t^\theta ) \leq \fC e^{- \fc t }$.
\end{prop}

\cfclear
\begin{cproof}{theo:flow:convergence}
Throughout this proof let $\cM \subseteq \R^\fd$ satisfy $\cM = \{ \theta \in \R^\fd \colon \cL ( \theta ) = 0 \}$.
\Nobs that \cref{cor:minima:manifold:1} proves that there exist $\mdim \in \N \cap [ 1 , \fd)$ and an open $U \subseteq \R^\fd$ which satisfy $U \subseteq \fV$,
which satisfy that $\cL |_U$ is twice continuously differentiable,
 which satisfy that $(\Hs  \cL)|_U$ is locally Lipschitz continuous,
  which satisfy that $\cM \cap U$ is a non-empty $\mdim$-dimensional $C^2$-submanifold of $\R^\fd$,
   and which satisfy for all $\theta \in \cM \cap U$ that $\rank ( ( \Hs  \cL ) ( \theta ) ) = \fd - \mdim$.
Combining this, \cref{lem:gradient:measurable}, 
\cref{prop:loss:differentiable},
\cref{lem:local:tubular:nbhd},
and \cref{prop:tubular:nbhd} with \cref{prop:gradient:flow:abstract} 
ensures that there exist $\fm \in \cM \cap U$, $\fc \in (0, \infty)$, $V , \cV \in \cu{A \subseteq U \colon A \text{ is compact}}$ which satisfy that
\begin{enumerate} [label = (\roman*)]
	\item it holds that $\fm \in V^ \circ \subseteq V \subseteq \cV$,
	\item it holds for all $\theta \in \cV$ that $\dist_{\cM \cap U } ( \theta ) = \dist_{\cM \cap U \cap \cV}$,
	\item it holds for all $\theta \in V$, $t \in [0, \infty)$ that $\Theta_t^\theta \in \cV$,
	and
	\item it holds for all $t \in [0, \infty )$ that $\dist_{\cM \cap U } ( \Theta_t^\theta ) \le e^{- \fc t } \dist_{\cM \cap U } ( \theta )$
\end{enumerate} \cfadd{def:norm}\cfadd{def:dist:set}\cfload. Furthermore, \nobs that the fact that $\cL |_U$ is twice continuously differentiable proves that there exists $ \fC \in (0, \infty)$ which satisfies for all $\theta, \vartheta \in \cV$ that $\abs{ \cL ( \theta ) - \cL ( \vartheta ) } \leq \fC \norm{\theta - \vartheta }$. This assures that for all $\theta \in V^\circ$, $t \in [0, \infty)$ we have that
\begin{equation}
\begin{split} 
    \cL ( \Theta_t^\theta ) &= \inf\nolimits_{\vartheta \in \cM \cap U \cap \cV} \abs{\cL ( \Theta_t^\theta ) - \cL ( \vartheta ) } \leq \fC \br[\big]{\inf\nolimits_{\vartheta \in \cM \cap U \cap \cV} \norm{ \Theta_t^\theta - \vartheta } } \\
    &= \fC \br[\big]{ \dist _{\cM \cap U} ( \Theta_t^\theta ) } \leq \fC e^{-\fc t} \dist _{\cM \cap U} ( \theta ).
    \end{split}
\end{equation}
\end{cproof}

\subsection{Convergence rates for GF with random initializations in the training of ANNs} \label{subsection:flow:random}

\cfclear
\begin{cor} \label{cor:flow:convergence:random} 
Assume \cref{setting:snn}, let $N \in \N \cap [1 , \width ]$, $\fx_0, \fx_1, \ldots, \fx_N, \alpha_1, \alpha_2, \ldots, \alpha_N \in \R$
satisfy $a = \fx_0 < \fx_1 < \cdots < \fx_N = b$, 
 assume for all $i \in \{1, 2, \ldots, N\}$, $x \in [\fx_{i-1}, \fx_i ]$ that $f(x) = f(\fx_{i-1} ) + \alpha_i ( x - \fx_{i-1} )$,
let $(\Omega, \cF, \P)$ be a probability space,
let $\Theta \colon [0, \infty ) \times \Omega \to \R^\fd$ be a stochastic process with continuous sample paths, assume that $\Theta_0$ is standard normally distributed, and assume for all $t \in [0, \infty) $, $\omega \in \Omega $ that
\begin{equation} \cfadd{lem:gradient:measurable}
    \Theta_t ( \omega ) = \Theta_0 ( \omega ) - \int_0^t \cG (\Theta_s ( \omega ) ) \, \d s
\end{equation}
\cfload.
Then there exist $\fc , \fC \in (0, \infty)$ such that $ \P ( \forall \, t \in [0, \infty) \colon \cL ( \Theta_t ) \leq \fC e^{- \fc t} ) > 0$.
\end{cor}
\begin{cproof2}{cor:flow:convergence:random}
\Nobs that \cref{theo:flow:convergence} ensures that 
there exist $\fc , \fC \in (0, \infty)$ and a non-empty open $U \subseteq \R^{\fd }$ 
which satisfy for all $t \in [0, \infty)$, $\omega \in \Omega$ 
with $\Theta_0 ( \omega ) \in U$ that $\cL ( \Theta_t ( \omega ) ) \leq \fC e^{- \fc t}$. 
\Nobs that 
the fact that $ U $ is a non-empty open set and 
the assumption that $\Theta_0$ is standard normally distributed 
imply that $ \P( \Theta_0 \in U ) > 0 $.
\end{cproof2}

\section{Local convergence to the set of global minima for gradient descent (GD)} \label{section:gd:convergence}

In this section we employ \cref{cor:minima:manifold:1} from \cref{section:risk:global:min} 
to establish in \cref{theo:gd:convergence} in \cref{subsection:gd:convergence}, 
\cref{cor:gd:convergence:random} in \cref{subsection:gd:random}, 
and \cref{cor:gd:random:init} in \cref{subsection:gd:random} 
under the assumption that the target function is piecewise affine linear 
that the risk of certain GD processes converges to zero. 
Our proofs of \cref{cor:gd:convergence:random,cor:gd:random:init} 
are based on an application of \cref{theo:gd:convergence} and 
our proof of \cref{theo:gd:convergence} uses the abstract local convergence result 
for GD processes in \cref{prop:gradient:descent:abstract} 
in \cref{subsection:gd:abstract} below. 
\cref{prop:gradient:descent:abstract} and its proof 
are strongly inspired by Fehrman et al.~\cite[Proposition 17]{FehrmanGessJentzen2020}.
Our proof of \cref{prop:gradient:descent:abstract} employs the elementary uniform estimate for certain exponential sums in \cref{lem:sum:estimate} in \cref{subsection:gd:abstract}. For completeness we include in this section also a detailed proof for \cref{lem:sum:estimate}.

\subsection{Abstract convergence result for GD to a submanifold of global minima} \label{subsection:gd:abstract}

\begin{lemma} \label{lem:sum:estimate}
Let $\rho \in [0,1)$, $c, \fg \in (0, \infty)$. Then there exists $\fC \in \R$ such that for all $\gamma \in (0, \fg ]$ it holds that
\begin{equation}
    \sum_{k=1}^\infty \gamma k^{- \rho} \exp \rbr*{ - c \gamma (k - 1 )^{1 - \rho } } \leq \fC.
\end{equation}
\end{lemma}
\begin{cproof}{lem:sum:estimate}
First \nobs that for all $\gamma \in (0, \fg]$ it holds that
\begin{equation} \label{eq:sum:estimate:1}
\begin{split}
         \sum_{k=1}^\infty \gamma k^{- \rho} \exp \rbr*{ - c \gamma (k - 1 )^{1 - \rho } } 
         &\leq \gamma + \sum_{k=2}^\infty  \gamma (k-1)^{- \rho} \exp \rbr*{ - c \gamma (k - 1 )^{1 - \rho } } \\
         &\leq \fg + \sum_{n=1}^\infty \gamma n^{-\rho} \exp \rbr*{ - c \gamma n^{1 - \rho } } \\
         &\leq 2 \fg + \sum_{n=2}^\infty \gamma n^{-\rho} \exp \rbr*{ - c \gamma n^{1 - \rho } }.
\end{split}
\end{equation}
Next \nobs that the fact that for all $\gamma \in (0, \infty)$ it holds that $[1, \infty) \ni x \mapsto x^{-\rho} \exp ( - c \gamma x^{1 - \rho} ) \in \R $ is continuous and non-increasing assures that for all $\gamma \in (0, \fg]$ it holds that
\begin{equation} \label{eq:sum:estimate:2}
\begin{split}
    \sum_{n=2}^\infty \gamma n^{-\rho} \exp \rbr*{ - c \gamma n^{1 - \rho } } 
    & \le \sum_{n=2}^\infty \br*{ \int_{n-1}^n \gamma x^{- \rho } \exp \rbr*{ - c \gamma x^{1 - \rho } } \, \d x } \\
    &= \int_1^\infty \gamma  x^{-\rho} \exp \rbr*{ - c \gamma x^{1 - \rho } } \, \d x.
\end{split}
\end{equation}
Moreover, \nobs that the integral transformation theorem proves for all $\gamma \in (0, \fg]$ that
\begin{equation} \label{eq:sum:estimate:3}
    \begin{split}
        & \int_1^\infty \gamma  x^{-\rho} \exp \rbr*{ - c \gamma x^{1 - \rho } } \, \d x =  \int_{\gamma^{1/(1-\rho)}}^\infty \gamma^{1 + \frac{\rho}{1 - \rho } } x^{-\rho}  \exp \rbr*{ - c x^{1 - \rho } } \gamma^{-\frac{1}{1 - \rho } } \, \d x \\
        &\leq \int_0^\infty x^{-\rho}  \exp \rbr*{ - c x^{1 - \rho } } \, \d x 
        \leq \int_0^1 x^{-\rho} \, \d x + \int_1^\infty \exp \rbr*{ - c x^{1 - \rho } } \, \d x \\
        &= \frac{1}{1-\rho} + \int_1^\infty \exp \rbr*{ - c x^{1 - \rho } } \, \d x.
    \end{split}
\end{equation}
Furthermore, \nobs that the assumption that $c \in (0, \infty )$ and the assumption that $\rho \in [0,1)$ ensure that $\int_1^\infty \exp \rbr*{ - c x^{1 - \rho } } \, \d x < \infty$. Combining this, \cref{eq:sum:estimate:1}, \cref{eq:sum:estimate:2}, and \cref{eq:sum:estimate:3} establishes for all $\gamma \in (0, \fg]$ that
\begin{equation}
     \sum_{k=1}^\infty \gamma k^{- \rho} \exp \rbr*{ - c \gamma (k - 1 )^{1 - \rho } } \leq 2 \fg + \frac{1}{1-\rho} + \int_1^\infty \exp \rbr*{ - c x^{1 - \rho } } \, \d x < \infty.
\end{equation}
\end{cproof}

\cfclear
\begin{prop} \label{prop:gradient:descent:abstract}
Assume \cref{setting:mfld}, assume for all $x \in \cM $ that $\rank ( (\Hs  f) (x ) ) = \fd-n$, let $\cG \colon \R^\fd \to \R^\fd$ satisfy for all $x \in U$ that $\cG ( x ) = (\nabla f ) ( x )$, let $x \in \cM$, $\rho \in [0,1)$,
 and let $\Theta^{\theta, \gamma} \colon \N_0 \to \R^\fd$, $\theta \in \R^\fd$, $\gamma \in \R$, satisfy for all $\theta \in \R^\fd$, $\gamma \in \R$, $n \in \N$ that $\Theta_0^{\theta, \gamma} = \theta $ and 
\begin{equation}
\Theta_n^{\theta, \gamma } = \Theta_{n-1}^{\theta, \gamma } - \tfrac{\gamma}{n^\rho} \cG ( \Theta_{n-1}^{\theta, \gamma } ).    
\end{equation}
 Then there exist $r , s \in (0, \infty)$ such that 
\begin{enumerate} [label = (\roman*)] \cfadd{def:minimal:eigenvalue} \cfadd{def:local:tubular:nbhd}
    \item it holds for all $\theta \in V^{r/2 , s }_{\cM , x } $, $\gamma \in (0, \min \{ [\sup_{y \in \cM \cap V^{r,s}_{\cM , x }} \Lambda ((\Hs  f) ( y ) )]^{-1}, 1 \} ]$, $n \in \N_0$ that $\Theta_n^{\theta, \gamma} \in V^{r , s }_{\cM , x } $,
    \item it holds that $ \inf_{y \in \cM \cap V^{r,s}_{\cM , x }} \br[\big]{ \sigma ((\Hs  f) ( y ) ) } > 0$, and
    \item it holds for all $\theta \in V^{r/2 , s }_{\cM , x } $, $\gamma \in (0, \min \{ [\sup_{y \in \cM \cap V^{r,s}_{\cM , x }} \Lambda ((\Hs  f) ( y ) )]^{-1}, 1 \} ]$, $n \in \N_0$ that 
    \begin{equation}
        \dist _\cM ( \Theta_n^{\theta, \gamma } ) \leq \exp \rbr*{- \tfrac{\gamma}{2(1 - \rho)} \br*{\inf\nolimits_{y \in \cM \cap V^{r,s}_{\cM , x }} \br[\big]{ \sigma ((\Hs  f) ( y ) ) } } n ^{1 - \rho}} \dist _\cM ( \theta )
    \end{equation}
\end{enumerate}
\cfload.
\end{prop}
\cfclear
\begin{cproof}{prop:gradient:descent:abstract}
\Nobs that \cref{prop:tubular:nbhd,lem:second:order:approx} prove that 
there exist $r, \varepsilon, \fc \in (0, \infty)$ 
which satisfy 
$\overline{V^{r, \varepsilon}_{\cM , x } }  \subseteq U $, 
which satisfy 
    \begin{equation} \cfadd{def:unique:projection}
      V^{r , s}_{\cM , x } = \cu*{y \in \tubn_\cM \colon \br[\big]{( \norm{x - \projm_\cM ( y ) } \le r ) , \, ( \norm{y - \projm_\cM ( y ) } < s)}},
    \end{equation}
and which satisfy for all $y \in \overline{V^{r, \varepsilon}_{\cM , x } }$ that
\begin{equation} \label{prop:gd:eq:hess}
  \norm*{ (\nabla f) ( y ) - (\Hs  f) ( \projm _\cM ( y ) ) ( y - \projm _\cM ( y ) )} \leq \fc ( \dist _\cM ( y ) )^2
\end{equation}
\cfload.
In the following let $\kappa \in \R $ satisfy 
$\kappa = \inf_{y \in \cM \cap V^{r , \varepsilon} _{\cM , x }} \br[\big]{ \sigma ((\Hs  f) ( y ) ) } $. 
\Nobs that the fact that $ U \ni y \mapsto (\Hs f)(y) \in \R^{ \fd \times \fd } $ 
is locally Lipschitz continuous and the fact 
that the eigenvalues are continuous functions of a matrix (cf., e.g., Kato~\cite[Section 2.5.1]{Kato1966}) 
prove that $ \kappa > 0 $. 
Next \nobs that the fact that $\overline{V^{r, \varepsilon}_{\cM , x } }$ is compact and 
the fact that $ U \ni y \mapsto (\nabla f)(y) \in \R^{ \fd } $ is continuously differentiable demonstrate that there exists $c \in (0, \infty)$ which satisfies for all $y \in \overline{V^{r, \varepsilon}_{\cM , x } }$ that
\begin{equation} \label{prop:gd:eq:gradient} \cfadd{def:norm} \cfadd{def:subset:unique:projection}
  \norm{( \nabla f ) ( y ) } = \norm{ ( \nabla f ) ( y ) - (\nabla f ) ( \projm _\cM ( y ) ) } \leq c \norm{y - \projm _\cM ( y ) } = c \dist _\cM ( y ) 
\end{equation}
     \cfload.
In the following let $\fC \in (0, \infty)$ satisfy for all $\gamma \in (0, 1]$ that
\begin{equation} \label{prop:gd:eq:sum:estimate}
    \sum_{k=1}^\infty \gamma k^{-\rho} \exp \rbr*{ - \tfrac{\kappa \gamma}{2(1 - \rho)} (k - 1 )^{1 - \rho} } \leq \fC
\end{equation}
(cf.\ \cref{lem:sum:estimate}), let $s \in (0, \infty)$ satisfy
\begin{equation} \label{prop:gd:eq:defdelta}
    s  = \min \cu* { \frac{\kappa}{ 2 \fc }, \frac{r}{2(2 +  c \fC )}, \varepsilon},
\end{equation}
let $\theta \in V^{r/2 , s }_{\cM , x } $ and $\gamma \in  (0, \min \{ [\sup_{y \in \cM \cap V^{r,s}_{\cM , x }} \Lambda ((\Hs  f) ( y ) )]^{-1}, 1 \} ]$ be arbitrary, and let $\tau \in \N \cup \{ \infty\}$ satisfy $\tau = \inf \{n \in \N_0 \colon \Theta_n^{\theta , \gamma} \notin V^{r , s }_{\cM , x }  \}$. \Nobs that the fact that for all $n \in \N \cap (0, \tau]$ it holds that $\Theta_n^{\theta, \gamma} \in V^{r , s }_{\cM , x } $ proves that for all $n \in \N \cap (0, \tau]$ we have that
\begin{equation}
    \begin{split}
        \dist _\cM ( \Theta_n^{\theta, \gamma} ) & \leq \norm{\Theta_n^{\theta, \gamma} - \projm _\cM ( \Theta_{n-1}^{\theta , \gamma } )} \\
        &= \norm*{\Theta_{n - 1 }^{\theta, \gamma } - \projm _\cM ( \Theta_{n-1}^{\theta , \gamma } ) - \tfrac{ \gamma }{n^\rho} (\nabla f ) ( \Theta_{n-1}^{\theta, \gamma})} \\
        &\leq \norm*{\Theta_{n - 1 }^{\theta, \gamma} - \projm _\cM ( \Theta_{n-1}^{\theta , \gamma } ) - \tfrac{ \gamma }{n^\rho} (\Hs  f )( \projm _\cM ( \Theta_{n-1}^{\theta, \gamma }) )  ( \Theta_{n-1}^{\theta, \gamma } - \projm _\cM ( \Theta_{n-1}^{\theta, \gamma} ) ) } \\
        &\quad + \tfrac{\gamma}{n^\rho}\norm*{ \rbr[\big]{(\Hs  f )( \projm _\cM ( \Theta_{n-1}^{\theta , \gamma }) ) }  (  \Theta_{n-1}^{\theta , \gamma } - \projm _\cM( \Theta_{n-1}^{\theta , \gamma } ) ) -  (\nabla f ) ( \Theta_{n-1}^{\theta , \gamma})}.
    \end{split}
\end{equation}
Combining this, \cref{prop:hessian:decomp}, and \cref{prop:gd:eq:hess} demonstrates for all $n \in \N \cap (0, \tau]$ that
\begin{equation}
    \begin{split}
          \dist  _\cM ( \Theta_n^{\theta , \gamma } ) & \leq  \rbr*{1 - \tfrac{\kappa \gamma }{n^\rho}} \dist  _\cM ( \Theta_{n - 1 }^{\theta , \gamma } ) + \tfrac{ \fc \gamma }{n^\rho} ( \dist _\cM ( \Theta_{n - 1 }^{\theta , \gamma } ) ) ^2.
    \end{split}
\end{equation}
This, the fact that for all $n \in \N \cap (0, \tau]$ it holds that $\dist _\cM ( \Theta_{n - 1 }^{\theta , \gamma }) \leq s \leq \frac{\kappa}{2 \fc }$, and \cref{prop:gd:eq:defdelta} imply that for all $n \in \N \cap (0, \tau]$ it holds that
\begin{equation}
      \dist _\cM ( \Theta_n^{\theta, \gamma } ) \leq \rbr*{1 - \tfrac{\kappa \gamma }{2 n^\rho}} \dist _\cM ( \Theta_{n - 1 }^{\theta , \gamma } ).
\end{equation}
By induction, we therefore obtain for all $n \in \N \cap (0, \tau]$ that
\begin{equation} \label{prop:gd:eq:induction}
      \dist _\cM ( \Theta_n^{\theta , \gamma } ) \leq \br*{\textstyle\prod_{k=1}^n \rbr*{1 - \tfrac{\kappa \gamma }{2 k^\rho}} } \dist _\cM ( \theta ).
\end{equation}
Next \nobs that the assumption that $\gamma \leq  [\sup_{y \in \cM \cap V^{r,s}_{\cM , x }} \Lambda ((\Hs  f) ( y ) )]^{-1} \leq \kappa^{-1}$ shows for all $k \in \N$ that $\tfrac{\kappa \gamma }{2 k^\rho} \in (0,1)$. This and the fact that for all $u \in (0,1)$ it holds that $\ln (1-u) \leq -u $ prove that for all $n \in \N$ we have that
\begin{equation}
    \ln \br*{\textstyle\prod_{k=1}^n \rbr*{1 - \tfrac{\kappa \gamma }{2 k^\rho}} } = \textstyle \sum_{k=1}^n \ln \rbr*{1 - \tfrac{\kappa \gamma }{2 k^\rho}} \leq - \frac{\kappa \gamma }{2} \sum_{k=1}^n k^{-\rho} \leq - \frac{\kappa \gamma }{2} \int_0^n u^{-\rho} \, \d u = \frac{\kappa \gamma }{2(1 - \rho)} n^{1 - \rho}.
\end{equation}
Combining this with \cref{prop:gd:eq:induction} demonstrates for all $n \in \N \cap ( 0 , \tau ]$ that
\begin{equation} \label{prop:gd:eq:dist:estimate}
      \dist _\cM ( \Theta_n^{\theta, \gamma } ) \leq \exp \rbr*{ - \tfrac{\kappa \gamma }{2(1 - \rho)} n^{1 - \rho} } \dist _\cM ( \theta ).
\end{equation}
It only remains to show that $\tau = \infty$. \Nobs that \cref{prop:gd:eq:gradient} assures for all $n \in \N \cap (0, \tau]$ that
\begin{equation}
    \norm{\Theta^{\theta, \gamma}_n - \Theta_{n-1}^{\theta, \gamma } } = \tfrac{\gamma }{n^\rho} \norm{ (\nabla f ) ( \Theta_{n-1}^{\theta, \gamma } )  } \leq \tfrac{c \gamma }{n^\rho} \dist _\cM ( \Theta_{n-1}^{\theta, \gamma } )
\end{equation}
This, \cref{prop:gd:eq:dist:estimate}, the fact that $\gamma \leq 1$, \cref{prop:gd:eq:sum:estimate},
and the triangle inequality establish for all $n \in \N \cap (0, \tau]$ that
\begin{equation}
\begin{split}
    \norm{\Theta_n^{\theta, \gamma}  - \theta } 
    &\leq \sum_{k=1}^n c \gamma k^{- \rho } \exp \rbr*{ - \tfrac{\kappa \gamma }{2(1 - \rho)} (k - 1 )^{1 - \rho} } \dist _\cM ( \theta ) \\
    &\leq c s \sum_{k=1}^\infty  \gamma k^{- \rho } \exp \rbr*{ - \tfrac{\kappa \gamma }{2(1 - \rho)} (k - 1 )^{1 - \rho} } \leq c s \fC.
    \end{split}
\end{equation}
Combining this with \cref{prop:gd:eq:dist:estimate}, \cref{prop:gd:eq:defdelta}, and the triangle inequality proves for all $n \in \N \cap (0, \tau]$ that
\begin{equation}
\begin{split}
    \norm{\projm _\cM ( \Theta_n^{\theta, \gamma } ) - \projm _\cM ( \theta ) } & \leq \dist _\cM ( \Theta_n^{\theta , \gamma } ) + \norm{\Theta_n^{\theta , \gamma }  - \theta } + \dist _\cM ( \theta ) \\
    &\leq s ( 2 + c \fC ) \leq \tfrac{r}{2}.
    \end{split}
\end{equation}
Furthermore, \nobs that the assumption that $\theta \in V^{r/2 , s }_{\cM , x } $ assures that there exists $\delta \in (0, \infty)$ which satisfies that $\theta \in V^{r/2 - \delta , s} _{\cM , x } $. Hence, we obtain for all $n \in \N \cap (0, \tau]$  that $\Theta_n^{\theta , \gamma } \in V^{r - \delta , s} _{\cM , x }$. This implies that $\tau = \infty$.
\end{cproof}

\subsection{Convergence rates for GD in the training of ANNs} \label{subsection:gd:convergence}

\begin{theorem} \label{theo:gd:convergence}
Assume \cref{setting:snn}, let $N \in \N \cap [1 , \width ]$, $\rho \in [0,1)$, $\fx_0, \fx_1, \ldots, \fx_N, \alpha_1, \alpha_2, \ldots, \allowbreak \alpha_N \in \R$
satisfy $a = \fx_0 < \fx_1 < \cdots < \fx_N = b$, 
 assume for all $i \in \{1, 2, \ldots, N\}$, $x \in [\fx_{i-1}, \fx_i ]$ that $f(x) = f(\fx_{i-1} ) + \alpha_i ( x - \fx_{i-1} )$,
 let $\tbound \in \R$ satisfy
\begin{equation}
    \tbound = 1 + \abs{f(a)} + ( 1 + 2 \max\nolimits_{j \in \{1, 2, \ldots, \width \} } \abs{\alpha_j} ) ( \abs{a} + \abs{b} + 1 ),
\end{equation}
and let $\Theta^{\theta , \gamma} \colon \N_0 \to \R^{\fd }$, $\theta \in \R^{\fd }$, $\gamma \in \R$, satisfy for all $\theta \in \R^{\fd }$, $\gamma \in \R$, $n \in \N$ that $\Theta_0^{\theta, \gamma} = \theta$ and
\begin{equation} \label{theo:gd:eq:def:descent}
\Theta_n^{\theta , \gamma} = \Theta^{\theta , \gamma}_{n-1} - \tfrac{\gamma}{n^{\rho } } \cG ( \Theta_{n-1}^{\theta , \gamma}).    
\end{equation}
Then there exist $\fc, \fC \in (0, \infty)$ 
and a non-empty open $U \subseteq (- \tbound , \tbound )^{\fd }$ such that for all $\theta \in U$, $\gamma \in (0, ((3 N + 1 ) \rbr{24 \tbound ^5 + 16 N \tbound ^7 } \rbr{ \sup\nolimits_{x \in [a,b]} \dens ( x ) } )^{-1}]$, $n \in \N_0$ it holds that $\cL ( \Theta_n^{\theta  , \gamma}) \leq \fC \exp (- \fc \gamma n^{1 - \rho } )$.
\end{theorem}
\cfclear
\begin{cproof}{theo:gd:convergence}
Throughout this proof let $\cM \subseteq \R^\fd$ satisfy $\cM = \{ \theta \in \R^\fd \colon \cL ( \theta ) = 0 \}$.
\Nobs that \cref{cor:minima:manifold:1} proves that there exist $\mdim \in \N \cap [ 1 , \fd)$ and an open $U \subseteq (- \tbound , \tbound )^\fd$ which satisfy $U \subseteq \fV$,
which satisfy that $\cL |_U$ is twice continuously differentiable,
which satisfiy for all $\theta \in U$ that $\Lambda ( ( \Hs \cL ) ( \theta ) ) \leq (3 N + 1 ) \rbr{24 \tbound ^5 + 16 N \tbound ^7 } \rbr{ \sup\nolimits_{x \in [a,b]} \dens ( x ) }$,
which satisfy that $(\Hs  \cL)|_U$ is locally Lipschitz continuous,
 which satisfy that $\cM \cap U$ is a non-empty $\mdim$-dimensional $C^2$-submanifold of $\R^\fd$, and which satisfy for all $\theta \in \cM \cap U$ that $\rank ( ( \Hs  \cL ) ( \theta ) ) = \fd - \mdim$.
Combining this, \cref{lem:gradient:measurable}, 
\cref{prop:loss:differentiable},
\cref{lem:local:tubular:nbhd},
and \cref{prop:tubular:nbhd}
 with \cref{prop:gradient:descent:abstract} 
shows that there exist $\fm \in \cM \cap U$, $\fc \in (0, \infty)$, 
$V, \cV \in \cu{A \subseteq U \colon A \text{ is compact}}$ such that
\begin{enumerate} [label = (\roman*)]
	\item it holds that $\fm \in V^\circ \subseteq V \subseteq \cV$,
	\item it holds for all $\theta \in \cV$ that $\dist_{\cM \cap U } ( \theta ) = \dist_{\cM \cap U \cap \cV} ( \theta )$,
	and
	\item it holds for all $\theta \in V$,
	$\gamma \in (0, \min \{ ( \sup_{\vartheta \in \cM \cap V_2}  \Lambda ( ( \Hs f ) ( \vartheta ) ) )^{-1}, 1 \} ]$, $n \in \N_0$ that $\Theta_n^{\theta , \gamma } \in \cV$ and $\dist _{\cM \cap U} ( \Theta_n^{\theta, \gamma } ) \leq \exp (- \fc \gamma n^{1 - \rho } ) \dist _{\cM \cap U} ( \theta )$ 
\end{enumerate}
\cfadd{def:norm}\cfadd{def:dist:set}\cfload. 
In addition, \nobs that
\begin{equation}
\begin{split}
    \sup\nolimits_{\vartheta \in \cM \cap \cV}  \Lambda ( ( \Hs f ) ( \vartheta ) )
    &\leq \sup\nolimits_{\vartheta \in U}  \Lambda ( ( \Hs f ) ( \vartheta ) ) \\
    &\leq (3 N + 1 ) \rbr[\big]{24 \tbound ^5 + 16 N \tbound ^7 } \rbr[\big]{ \sup\nolimits_{x \in [a,b]} \dens ( x ) }.
    \end{split}
\end{equation}
Furthermore, \nobs that the fact that $\cL |_U$ is twice continuously differentiable implies that there exists $\fC \in (0, \infty)$ which satisfies for all $\theta, \vartheta \in \cV$ that $\abs{ \cL ( \theta ) - \cL ( \vartheta ) } \leq \fC \norm{\theta - \vartheta }$. This ensures that for all $\theta \in V^\circ$, $\gamma \in (0, ((3 N + 1 ) \rbr{16 \tbound ^5 + 8 N \tbound ^7 } \rbr{ \sup\nolimits_{x \in [a,b]} \dens ( x ) } )^{-1}]$, $n \in \N_0$ we have that
\begin{equation}
\begin{split}
  \cL ( \Theta_n^{\theta , \gamma } ) 
&
  = 
  \inf\nolimits_{\vartheta \in \cM \cap U \cap \cV} \abs{\cL ( \Theta_n^{\theta , \gamma } ) - \cL ( \vartheta ) } 
  \leq 
  \fC 
  \big[ 
    \inf\nolimits_{\vartheta \in \cM \cap U \cap \cV} \norm{ \Theta_n^{\theta , \gamma } - \vartheta } 
  \big] 
\\
& 
  = 
  \fC \big[ \dist_{ \cM \cap U }( \Theta_n^{ \theta , \gamma } ) \big]
  \leq 
  \fC 
  \exp(- \fc \gamma n^{1 - \rho } ) 
  \dist_{ \cM \cap U }( \theta ) 
  .
    \end{split}
\end{equation}
\end{cproof}

\subsection{Convergence results for GD with random initializations in the training of ANNs} \label{subsection:gd:random}

\begin{cor} \label{cor:gd:convergence:random} 
Assume \cref{setting:snn}, let $N \in \N \cap [ 1 , \width ]$, $\rho \in [0, 1)$,
$\fx_0, \fx_1, \ldots, \fx_N, \alpha_1, \alpha_2, \ldots, \allowbreak \alpha_N \in \R$
satisfy $a = \fx_0 < \fx_1 < \cdots < \fx_N = b$,
 assume for all $i \in \{1, 2, \ldots, N\}$, $x \in [\fx_{i-1}, \fx_i ]$ that $f(x) = f(\fx_{i-1} ) + \alpha_i ( x - \fx_{i-1} )$,
 let $\tbound \in \R$ satisfy
\begin{equation}
    \tbound = 1 + \abs{f(a)} + ( 1 + 2 \max\nolimits_{j \in \{1, 2, \ldots, \width \} } \abs{\alpha_j} ) ( \abs{a} + \abs{b} + 1 ),
\end{equation}
let $(\Omega, \cF, \P)$ be a probability space,
let $\Theta_n^\gamma \colon \Omega \to \R^\fd$, $\gamma \in \R$, $n \in \N_0$,
be random variables,
 assume for all $\gamma \in \R$ that $\Theta_0^\gamma $ is standard normally distributed, and assume for all $\gamma \in \R$, $n \in \N$, $\omega \in \Omega$ that
\begin{equation}
    \Theta_n^\gamma ( \omega ) = \Theta_{n-1}^\gamma ( \omega ) - \gamma n^{-\rho} \cG (\Theta_{n-1}^\gamma ( \omega ) ) .
\end{equation}
Then there exist $\fc, \fC \in (0, \infty)$ such that for all $\gamma \in (0, ((3 N + 1 ) \rbr{24 \tbound ^5 + 16 N \tbound ^7 } \rbr{ \sup\nolimits_{x \in [a,b]} \dens ( x ) } )^{-1}]$ it holds that $ \P ( \forall \, n \in \N_0 \colon \cL ( \Theta_n^\gamma ) \leq \fC \exp ( - \fc \gamma n^{1 - \rho } ) ) \geq \fc$.
\end{cor}
\begin{cproof2}{cor:gd:convergence:random}
\Nobs that \cref{theo:gd:convergence} ensures 
that there exist $\fc, \fC \in (0, \infty)$ 
and a non-empty open $ U \subseteq \R^{\fd } $ 
such that for all $\gamma \in (0, ((3 N + 1 ) \rbr{24 \tbound ^5 + 16 N \tbound ^7 } \rbr{ \sup\nolimits_{x \in [a,b]} \dens ( x ) } )^{-1}]$,
$\omega \in \Omega$, $n \in \N_0$ with $\Theta_0^\gamma ( \omega ) \in U$ it holds that 
\begin{equation}
\cL ( \Theta_n^\gamma ( \omega ) ) \leq \fC \exp ( - \fc \gamma n^{1 - \rho } ).    
\end{equation}
\Nobs that the fact that $ U $ is a non-empty open set 
and the assumption that for all $\gamma \in \R$ it holds that $\Theta_0^\gamma $ 
is standard normally distributed imply that there exists $\delta \in (0, \infty)$ such that for all $\gamma \in \R$ we have that $\P ( \Theta_0^\gamma \in U ) \geq \delta$.
\end{cproof2}

\begin{cor} \label{cor:gd:random:init}
Assume \cref{setting:snn}, let $N \in \N \cap [ 1 , \width ]$, 
$\fx_0, \fx_1, \ldots, \fx_N, \alpha_1, \alpha_2, \ldots, \allowbreak \alpha_N \in \R$
satisfy $a = \fx_0 < \fx_1 < \cdots < \fx_N = b$, 
 assume for all $i \in \{1, 2, \ldots, N\}$, $x \in [\fx_{i-1}, \fx_i ]$ that $f(x) = f(\fx_{i-1} ) + \alpha_i ( x - \fx_{i-1} )$,
  let $\tbound \in \R$ satisfy
\begin{equation}
    \tbound = 1 + \abs{f(a)} + ( 1 + 2 \max\nolimits_{j \in \{1, 2, \ldots, \width \} } \abs{\alpha_j} ) ( \abs{a} + \abs{b} + 1 ),
\end{equation}
 let $\Theta^{k, \gamma}_n \colon  \Omega \to \R^\fd$, $k , n \in \N_0$, $\gamma \in \R$,
 and $\bfk^{ k , \gamma}_n \colon \Omega \to \N$, $k , n \in \N_0$, $\gamma \in \R$,  be random variables,
 assume for all $\gamma \in \R$ that $\Theta_0^{k, \gamma}$, $k \in \N$, are independent standard normal random variables, 
 and assume for all $k \in \N$, $\gamma \in \R$, $n \in \N_0$,
 $\omega \in \Omega$ that
 \begin{equation}
      \Theta_{n+1}^{ k , \gamma} ( \omega ) = \Theta_{n}^{ k , \gamma} ( \omega ) - \gamma \cG (\Theta_{n}^{ k , \gamma } ( \omega ) ) 
 \end{equation}
 and
\begin{equation} \label{cor:gd:random:eq:defk}
    \bfk^{ k , \gamma}_n ( \omega) \in \arg \min\nolimits_{\ell \in \{1, 2, \ldots, k \} } \cL ( \Theta_n^{ \ell , \gamma } ( \omega ) ) .
\end{equation}
Then it holds for all $\gamma \in (0, ((3 N + 1 ) \rbr{24 \tbound ^5 + 16 N \tbound ^7 } \rbr{ \sup\nolimits_{x \in [a,b]} \dens ( x ) } )^{-1}]$  that
\begin{equation}
  \liminf \nolimits_{K \to \infty} \P \rbr*{ \limsup\nolimits_{n \to \infty} \cL \rbr[\big]{ \Theta^{  \bfk^{ K , \gamma}_n, \gamma}_n } = 0  } = 1.
\end{equation}
\end{cor}
\begin{cproof} {cor:gd:random:init}
Throughout this proof let $\fg \in \R$ satisfy $\fg = ((3 N + 1 ) \rbr{24 \tbound ^5 + 16 N \tbound ^7 } \allowbreak \rbr{ \sup\nolimits_{x \in [a,b]} \dens ( x ) } )^{-1}$.
\Nobs that \cref{theo:gd:convergence} assures that there exist $\fc, \fC \in (0, \infty)$ and an open $U \subseteq (- \fD , \fD )^{\fd}$ 
such that for all $\gamma \in (0, \fg]$, $k \in \N$, $\omega \in \Omega$, $n \in \N_0$ with $\Theta_0^{ k, \gamma} ( \omega) \in U$ it holds that $\cL ( \Theta_n^{ k, \gamma} ( \omega ) ) \leq \fC \exp ( - \fc \gamma n )$. Hence, we obtain for all $\gamma \in (0, \fg]$, $k \in \N$, $\omega \in \Omega$ with $\Theta_0^{ k, \gamma} ( \omega ) \in U$ that $\limsup_{n \to \infty} \cL ( \Theta_n^{k, \gamma} ( \omega ) ) = 0$.
Next \nobs that \cref{cor:gd:random:eq:defk} ensures for all $K \in \N$,
$\gamma \in (0, \fg]$ that
\begin{equation} \label{cor:gd:random:proof:eq1}
    \P \rbr*{ \limsup\nolimits_{n \to \infty} \cL \rbr[\big]{ \Theta^{  \bfk^{ K , \gamma}_n, \gamma}_n } = 0 } \geq \P \rbr*{\exists \, k \in \{1, 2, \ldots, K \} \colon \br[\big]{ \limsup\nolimits_{n \to \infty} \cL ( \Theta_n^{k, \gamma} ) = 0 } }.
\end{equation}
Furthermore, \nobs that the fact that for all $\gamma \in (0, \fg]$, $k \in \N$, $\omega \in \Omega$ with $\Theta_0^{ k, \gamma} ( \omega ) \in U$ it holds that $\limsup_{n \to \infty} \cL ( \Theta_n^{k, \gamma} ( \omega ) ) = 0$ shows that for all $K \in \N$, $\gamma \in (0, \fg]$ it holds that
\begin{equation} \label{cor:gd:random:proof:eq2}
    \P \rbr*{\exists \, k \in \{1, 2, \ldots, K \} \colon \br[\big]{ \limsup\nolimits_{n \to \infty} \cL ( \Theta_n^{k, \gamma} ) = 0 } } \geq \P \rbr*{ \exists \, k \in \{1, 2, \ldots, K \} \colon \Theta_0^{k, \gamma} \in U}.
\end{equation}
In addition, \nobs that the fact that for all $\gamma \in \R$ it holds that $\Theta_0^{k, \gamma}$, $k \in \N$, are i.i.d.\ implies that for all $K \in \N$, $\gamma \in (0, \fg]$ it holds that
\begin{equation} \label{cor:gd:random:proof:eq3}
\begin{split}
    \P \rbr*{ \exists \, k \in \{1, 2, \ldots, K \} \colon \Theta_0^{k, \gamma} \in U} &= 1 - \P \rbr*{ \forall \, k \in \{1, 2, \ldots, K \} \colon \Theta_0^{k, \gamma} \in (\R^\fd \backslash U ) } \\
    &= 1 - \br[\big]{ \P \rbr[\big]{\Theta_0^{1, \gamma} \in ( \R^\fd \backslash U ) } } ^K.
    \end{split}
\end{equation}
Moreover, \nobs that the fact that $U$ is open
and the fact that for all $\gamma \in \R$ it holds that $\Theta_0^{1, \gamma}$ is standard normally distributed prove that for all $\gamma \in \R$ it holds that $\P \rbr[\big]{\Theta_0^{1, \gamma} \in ( \R^\fd \backslash U ) } < 1$. This and \cref{cor:gd:random:proof:eq3} demonstrate for all $\gamma \in (0, \fg]$ that
\begin{equation} 
    \liminf \nolimits_{K \to \infty} \P \rbr*{ \exists \, k \in \{1, 2, \ldots, K \} \colon \Theta_0^{k, \gamma} \in U} = 1.
\end{equation}
Combining this with \cref{cor:gd:random:proof:eq1,cor:gd:random:proof:eq2} shows for all $\gamma \in (0, \fg]$ that
\begin{equation}
  \liminf \nolimits_{K \to \infty} \P \rbr*{ \limsup\nolimits_{n \to \infty} \cL \rbr[\big]{ \Theta^{  \bfk^{ K , \gamma}_n, \gamma}_n } = 0  } = 1.
\end{equation}
\end{cproof}

\subsection*{Acknowledgements}
This work has been funded by the Deutsche Forschungsgemeinschaft
(DFG, German Research Foundation) under Germany’s Excellence Strategy EXC 2044-390685587, Mathematics Münster: Dynamics-Geometry-Structure.


\begin{thebibliography}{10}

\bibitem{AbsilMahonyAndrews2005}
P.-A. Absil, R.~Mahony, and B.~Andrews.
\newblock Convergence of the iterates of descent methods for analytic cost
  functions.
\newblock {\em SIAM J. Optim.}, 16(2):531--547, 2005.
\newblock \href {https://doi.org/10.1137/040605266}
  {\path{doi:10.1137/040605266}}.

\bibitem{AkyildizSabanis2021}
{\"O}mer~Deniz Akyildiz and Sotirios Sabanis.
\newblock Nonasymptotic analysis of {S}tochastic {G}radient {H}amiltonian
  {M}onte {C}arlo under local conditions for nonconvex optimization, 2021.
\newblock \href {http://arxiv.org/abs/2002.05465} {\path{arXiv:2002.05465}}.

\bibitem{AllenzhuLiLiang2019}
Zeyuan Allen-Zhu, Yuanzhi Li, and Yingyu Liang.
\newblock Learning and generalization in overparameterized neural networks,
  going beyond two layers.
\newblock In H.~Wallach, H.~Larochelle, A.~Beygelzimer, F.~d\textquotesingle
  Alch\'{e}-Buc, E.~Fox, and R.~Garnett, editors, {\em Advances in Neural
  Information Processing Systems}, volume~32, pages 6158--6169. Curran
  Associates, Inc., 2019.
\newblock URL:
  \url{https://proceedings.neurips.cc/paper/2019/file/62dad6e273d32235ae02b7d321578ee8-Paper.pdf}.

\bibitem{AllenzhuLiSong2019}
Zeyuan Allen-Zhu, Yuanzhi Li, and Zhao Song.
\newblock A convergence theory for deep learning via over-parameterization.
\newblock In Kamalika Chaudhuri and Ruslan Salakhutdinov, editors, {\em
  Proceedings of the 36th International Conference on Machine Learning},
  volume~97 of {\em Proceedings of Machine Learning Research}, pages 242--252.
  PMLR, 09--15 Jun 2019.
\newblock URL: \url{http://proceedings.mlr.press/v97/allen-zhu19a.html}.

\bibitem{AroraDuHuLiWang2019}
Sanjeev Arora, Simon Du, Wei Hu, Zhiyuan Li, and Ruosong Wang.
\newblock Fine-grained analysis of optimization and generalization for
  overparameterized two-layer neural networks.
\newblock In Kamalika Chaudhuri and Ruslan Salakhutdinov, editors, {\em
  Proceedings of the 36th International Conference on Machine Learning},
  volume~97 of {\em Proceedings of Machine Learning Research}, pages 322--332,
  Long Beach, California, USA, 6 2019. PMLR.
\newblock URL: \url{http://proceedings.mlr.press/v97/arora19a.html}.

\bibitem{AttouchBolte2009}
Hedy Attouch and J\'{e}r\^{o}me Bolte.
\newblock On the convergence of the proximal algorithm for nonsmooth functions
  involving analytic features.
\newblock {\em Math. Program.}, 116(1-2, Ser. B):5--16, 2009.
\newblock \href {https://doi.org/10.1007/s10107-007-0133-5}
  {\path{doi:10.1007/s10107-007-0133-5}}.

\bibitem{BachMoulines2013}
Francis Bach and Eric Moulines.
\newblock Non-strongly-convex smooth stochastic approximation with convergence
  rate {$O(1/n)$}.
\newblock In C.~J.~C. Burges, L.~Bottou, M.~Welling, Z.~Ghahramani, and K.~Q.
  Weinberger, editors, {\em Advances in Neural Information Processing Systems},
  volume~26, pages 773--781. Curran Associates, Inc., 2013.
\newblock URL:
  \url{http://papers.nips.cc/paper/4900-non-strongly-convex-smooth-stochastic-approximation-with-convergence-rate-o1n.pdf}.

\bibitem{BercuFort2013}
Bernard Bercu and Jean-Claude Fort.
\newblock {\em Generic Stochastic Gradient Methods}, pages 1--8.
\newblock American Cancer Society, 2013.
\newblock URL: \url{https://doi.org/10.1002/9780470400531.eorms1068}.

\bibitem{BertsekasTsitsiklis2000}
Dimitri~P. Bertsekas and John~N. Tsitsiklis.
\newblock Gradient convergence in gradient methods with errors.
\newblock {\em SIAM Journal on Optimization}, 10(3):627--642, 2000.
\newblock \href {https://doi.org/10.1137/S1052623497331063}
  {\path{doi:10.1137/S1052623497331063}}.

\bibitem{Bottou2018optimization}
Léon Bottou, Frank~E. Curtis, and Jorge Nocedal.
\newblock Optimization methods for large-scale machine learning, 2018.
\newblock \href {http://arxiv.org/abs/1606.04838} {\path{arXiv:1606.04838}}.

\bibitem{CheriditoJentzenRiekert2021}
Patrick Cheridito, Arnulf Jentzen, Adrian Riekert, and Florian Rossmannek.
\newblock A proof of convergence for gradient descent in the training of
  artificial neural networks for constant target functions, 2021.
\newblock \href {http://arxiv.org/abs/2102.09924} {\path{arXiv:2102.09924}}.

\bibitem{CheriditoJentzenRossmannek2020}
Patrick Cheridito, Arnulf Jentzen, and Florian Rossmannek.
\newblock Non-convergence of stochastic gradient descent in the training of
  deep neural networks.
\newblock {\em Journal of Complexity}, page 101540, 2020.
\newblock \href {https://doi.org/10.1016/j.jco.2020.101540}
  {\path{doi:10.1016/j.jco.2020.101540}}.

\bibitem{CheriditoJentzenRossmannek2021}
Patrick Cheridito, Arnulf Jentzen, and Florian Rossmannek.
\newblock Landscape analysis for shallow {R}e{LU} neural networks: complete
  classification of critical points for affine target functions, 2021.
\newblock \href {http://arxiv.org/abs/2103.10922} {\path{arXiv:2103.10922}}.

\bibitem{DereichKassing2021}
Steffen Dereich and Sebastian Kassing.
\newblock Convergence of stochastic gradient descent schemes for
  {L}ojasiewicz-landscapes, 2021.
\newblock \href {http://arxiv.org/abs/2102.09385} {\path{arXiv:2102.09385}}.

\bibitem{DereichMuller_Gronbach2019}
Steffen Dereich and Thomas M\"{u}ller-Gronbach.
\newblock General multilevel adaptations for stochastic approximation
  algorithms of {R}obbins-{M}onro and {P}olyak-{R}uppert type.
\newblock {\em Numer. Math.}, 142(2):279--328, 2019.
\newblock \href {https://doi.org/10.1007/s00211-019-01024-y}
  {\path{doi:10.1007/s00211-019-01024-y}}.

\bibitem{DuLeeLiWangZhai2019}
Simon Du, Jason Lee, Haochuan Li, Liwei Wang, and Xiyu Zhai.
\newblock Gradient descent finds global minima of deep neural networks.
\newblock In Kamalika Chaudhuri and Ruslan Salakhutdinov, editors, {\em
  Proceedings of the 36th International Conference on Machine Learning},
  volume~97 of {\em Proceedings of Machine Learning Research}, pages
  1675--1685, Long Beach, California, USA, 6 2019. PMLR.
\newblock URL: \url{http://proceedings.mlr.press/v97/du19c.html}.

\bibitem{DuZhaiPoczosSingh2018}
Simon~S. Du, Xiyu Zhai, Barnabas Poczos, and Aarti Singh.
\newblock Gradient descent provably optimizes over-parameterized neural
  networks.
\newblock In {\em International Conference on Learning Representations}, 2019.
\newblock URL: \url{https://openreview.net/forum?id=S1eK3i09YQ}.

\bibitem{EMaWojtowytschWu2020}
Weinan E, Chao Ma, Stephan Wojtowytsch, and Lei Wu.
\newblock Towards a mathematical understanding of neural network-based machine
  learning: what we know and what we don't, 2020.
\newblock \href {http://arxiv.org/abs/2009.10713} {\path{arXiv:2009.10713}}.

\bibitem{EMaWu2020}
Weinan E, Chao Ma, and Lei Wu.
\newblock A comparative analysis of optimization and generalization properties
  of two-layer neural network and random feature models under gradient descent
  dynamics.
\newblock {\em Sci. China Math.}, 63(7):1235--1258, 2020.
\newblock \href {https://doi.org/10.1007/s11425-019-1628-5}
  {\path{doi:10.1007/s11425-019-1628-5}}.

\bibitem{FehrmanGessJentzen2020}
Benjamin Fehrman, Benjamin Gess, and Arnulf Jentzen.
\newblock Convergence rates for the stochastic gradient descent method for
  non-convex objective functions.
\newblock {\em J. Mach. Learn. Res.}, 21:Paper No. 136, 48, 2020.

\bibitem{Ge2015}
Rong Ge, Furong Huang, Chi Jin, and Yang Yuan.
\newblock Escaping from saddle points --- online stochastic gradient for tensor
  decomposition.
\newblock In Peter Grünwald, Elad Hazan, and Satyen Kale, editors, {\em
  Proceedings of The 28th Conference on Learning Theory}, volume~40 of {\em
  Proceedings of Machine Learning Research}, pages 797--842, Paris, France,
  03--06 Jul 2015. PMLR.

\bibitem{GolubLoan2013}
Gene~H. Golub and Charles~F. Van~Loan.
\newblock {\em Matrix computations}.
\newblock Johns Hopkins Studies in the Mathematical Sciences. Johns Hopkins
  University Press, Baltimore, MD, fourth edition, 2013.

\bibitem{GoodfellowBengioCourville2016}
Ian Goodfellow, Yoshua Bengio, and Aaron Courville.
\newblock {\em Deep learning}.
\newblock Adaptive Computation and Machine Learning. MIT Press, Cambridge, MA,
  2016.

\bibitem{JentzenKroeger2021}
Arnulf Jentzen and Timo Kröger.
\newblock Convergence rates for gradient descent in the training of
  overparameterized artificial neural networks with biases, 2021.
\newblock \href {http://arxiv.org/abs/2102.11840} {\path{arXiv:2102.11840}}.

\bibitem{JentzenKuckuckNeufeldVonWurstemberger2021}
Arnulf Jentzen, Benno Kuckuck, Ariel Neufeld, and Philippe von Wurstemberger.
\newblock Strong error analysis for stochastic gradient descent optimization
  algorithms.
\newblock {\em IMA J. Numer. Anal.}, 41(1):455--492, 2021.
\newblock \href {https://doi.org/10.1093/imanum/drz055}
  {\path{doi:10.1093/imanum/drz055}}.

\bibitem{JentzenRiekertFlow}
Arnulf Jentzen and Adrian Riekert.
\newblock Convergence analysis for gradient flows in the training of artificial
  neural networks with {R}e{LU} activation, 2021.
\newblock \href {http://arxiv.org/abs/2107.04479} {\path{arXiv:2107.04479}}.

\bibitem{JentzenRiekert2021}
Arnulf Jentzen and Adrian Riekert.
\newblock A proof of convergence for stochastic gradient descent in the
  training of artificial neural networks with {R}e{LU} activation for constant
  target functions, 2021.
\newblock \href {http://arxiv.org/abs/2104.00277} {\path{arXiv:2104.00277}}.

\bibitem{JentzenvonWurstemberger2020}
Arnulf Jentzen and Philippe {von Wurstemberger}.
\newblock Lower error bounds for the stochastic gradient descent optimization
  algorithm: Sharp convergence rates for slowly and fast decaying learning
  rates.
\newblock {\em Journal of Complexity}, 57:101438, 2020.
\newblock \href {https://doi.org/10.1016/j.jco.2019.101438}
  {\path{doi:10.1016/j.jco.2019.101438}}.

\bibitem{Karimi2020linear}
Hamed Karimi, Julie Nutini, and Mark Schmidt.
\newblock Linear convergence of gradient and proximal-gradient methods under
  the {P}olyak-{L}ojasiewicz condition, 2020.
\newblock \href {http://arxiv.org/abs/1608.04636} {\path{arXiv:1608.04636}}.

\bibitem{Kato1966}
Tosio Kato.
\newblock {\em Perturbation theory for linear operators}.
\newblock Classics in Mathematics. Springer-Verlag, Berlin, 1995.
\newblock Reprint of the 1980 edition.

\bibitem{LeePanageasRecht2019}
Jason~D. Lee, Ioannis Panageas, Georgios Piliouras, Max Simchowitz, Michael~I.
  Jordan, and Benjamin Recht.
\newblock First-order methods almost always avoid strict saddle points.
\newblock {\em Math. Program.}, 176(1–2):311–337, July 2019.
\newblock \href {https://doi.org/10.1007/s10107-019-01374-3}
  {\path{doi:10.1007/s10107-019-01374-3}}.

\bibitem{LeeJordanRecht2016}
Jason~D. Lee, Max Simchowitz, Michael~I. Jordan, and Benjamin Recht.
\newblock Gradient descent only converges to minimizers.
\newblock In Vitaly Feldman, Alexander Rakhlin, and Ohad Shamir, editors, {\em
  29th Annual Conference on Learning Theory}, volume~49 of {\em Proceedings of
  Machine Learning Research}, pages 1246--1257, Columbia University, New York,
  New York, USA, 23--26 Jun 2016. PMLR.
\newblock URL: \url{http://proceedings.mlr.press/v49/lee16.html}.

\bibitem{LeiHuLiTang2020}
Y.~{Lei}, T.~{Hu}, G.~{Li}, and K.~{Tang}.
\newblock Stochastic gradient descent for nonconvex learning without bounded
  gradient assumptions.
\newblock {\em IEEE Transactions on Neural Networks and Learning Systems},
  31(10):4394--4400, 2020.
\newblock \href {https://doi.org/10.1109/TNNLS.2019.2952219}
  {\path{doi:10.1109/TNNLS.2019.2952219}}.

\bibitem{LiLiang2019}
Yuanzhi Li and Yingyu Liang.
\newblock Learning overparameterized neural networks via stochastic gradient
  descent on structured data.
\newblock In S.~Bengio, H.~Wallach, H.~Larochelle, K.~Grauman, N.~Cesa-Bianchi,
  and R.~Garnett, editors, {\em Advances in Neural Information Processing
  Systems}, volume~31, pages 8157--8166. Curran Associates, Inc., 2018.
\newblock URL:
  \url{https://proceedings.neurips.cc/paper/2018/file/54fe976ba170c19ebae453679b362263-Paper.pdf}.

\bibitem{LovasSabanis2020}
Attila Lovas, Iosif Lytras, Miklós Rásonyi, and Sotirios Sabanis.
\newblock Taming neural networks with {TUSLA}: Non-convex learning via adaptive
  stochastic gradient {L}angevin algorithms, 2020.
\newblock \href {http://arxiv.org/abs/2006.14514} {\path{arXiv:2006.14514}}.

\bibitem{LuShinSuKarniadakis2020}
Lu~Lu, Yeonjong Shin, Yanhui Su, and George~Em Karniadakis.
\newblock Dying {R}e{LU} and initialization: Theory and numerical examples.
\newblock {\em Communications in Computational Physics}, 28(5):1671--1706,
  2020.
\newblock \href {https://doi.org/10.4208/cicp.OA-2020-0165}
  {\path{doi:10.4208/cicp.OA-2020-0165}}.

\bibitem{BachMoulines2011}
Eric Moulines and Francis Bach.
\newblock Non-asymptotic analysis of stochastic approximation algorithms for
  machine learning.
\newblock In J.~Shawe-Taylor, R.~Zemel, P.~Bartlett, F.~Pereira, and K.~Q.
  Weinberger, editors, {\em Advances in Neural Information Processing Systems},
  volume~24, pages 451--459. Curran Associates, Inc., 2011.
\newblock URL:
  \url{https://proceedings.neurips.cc/paper/2011/file/40008b9a5380fcacce3976bf7c08af5b-Paper.pdf}.

\bibitem{Nesterov2015}
Yu~Nesterov.
\newblock Universal gradient methods for convex optimization problems.
\newblock {\em Math. Program.}, 152(1-2, Ser. A):381--404, 2015.
\newblock \href {https://doi.org/10.1007/s10107-014-0790-0}
  {\path{doi:10.1007/s10107-014-0790-0}}.

\bibitem{Nesterov2004}
Yurii Nesterov.
\newblock {\em Introductory lectures on convex optimization}, volume~87 of {\em
  Applied Optimization}.
\newblock Kluwer Academic Publishers, Boston, MA, 2004.
\newblock A basic course.
\newblock \href {https://doi.org/10.1007/978-1-4419-8853-9}
  {\path{doi:10.1007/978-1-4419-8853-9}}.

\bibitem{PanageasPiliouras2017}
Ioannis Panageas and Georgios Piliouras.
\newblock {Gradient Descent Only Converges to Minimizers: Non-Isolated Critical
  Points and Invariant Regions}.
\newblock In Christos~H. Papadimitriou, editor, {\em 8th Innovations in
  Theoretical Computer Science Conference (ITCS 2017)}, volume~67 of {\em
  Leibniz International Proceedings in Informatics (LIPIcs)}, pages 2:1--2:12,
  Dagstuhl, Germany, 2017. Schloss Dagstuhl--Leibniz-Zentrum fuer Informatik.
\newblock \href {https://doi.org/10.4230/LIPIcs.ITCS.2017.2}
  {\path{doi:10.4230/LIPIcs.ITCS.2017.2}}.

\bibitem{Panageas2019firstorder}
Ioannis Panageas, Georgios Piliouras, and Xiao Wang.
\newblock First-order methods almost always avoid saddle points: the case of
  vanishing step-sizes, 2019.
\newblock \href {http://arxiv.org/abs/1906.07772} {\path{arXiv:1906.07772}}.

\bibitem{Patel2021stopping}
Vivak Patel.
\newblock Stopping criteria for, and strong convergence of, stochastic gradient
  descent on {B}ottou-{C}urtis-{N}ocedal functions, 2021.
\newblock \href {http://arxiv.org/abs/2004.00475} {\path{arXiv:2004.00475}}.

\bibitem{Rakhlin2012}
Alexander Rakhlin, Ohad Shamir, and Karthik Sridharan.
\newblock Making gradient descent optimal for strongly convex stochastic
  optimization.
\newblock In {\em Proceedings of the 29th International Conference on Machine
  Learning}, page 1571–1578, Madison, WI, USA, 2012. Omnipress.

\bibitem{RotskoffEijnden2018}
Grant Rotskoff and Eric Vanden-Eijnden.
\newblock Parameters as interacting particles: long time convergence and
  asymptotic error scaling of neural networks.
\newblock In S.~Bengio, H.~Wallach, H.~Larochelle, K.~Grauman, N.~Cesa-Bianchi,
  and R.~Garnett, editors, {\em Advances in Neural Information Processing
  Systems}, volume~31. Curran Associates, Inc., 2018.
\newblock URL:
  \url{https://proceedings.neurips.cc/paper/2018/file/196f5641aa9dc87067da4ff90fd81e7b-Paper.pdf}.

\bibitem{Ruder2017overview}
Sebastian Ruder.
\newblock An overview of gradient descent optimization algorithms, 2017.
\newblock \href {http://arxiv.org/abs/1609.04747} {\path{arXiv:1609.04747}}.

\bibitem{SankararamanDeXuHuangGoldstein2020}
Karthik~A. Sankararaman, Soham De, Zheng Xu, W.~Ronny Huang, and Tom Goldstein.
\newblock The impact of neural network overparameterization on gradient
  confusion and stochastic gradient descent, 2020.
\newblock \href {http://arxiv.org/abs/1904.06963} {\path{arXiv:1904.06963}}.

\bibitem{SchmidtleRoux2013}
Mark Schmidt and Nicolas~Le Roux.
\newblock Fast convergence of stochastic gradient descent under a strong growth
  condition, 2013.
\newblock \href {http://arxiv.org/abs/1308.6370} {\path{arXiv:1308.6370}}.

\bibitem{Shamir2019}
Ohad Shamir.
\newblock Exponential convergence time of gradient descent for
  one-di\-men\-sion\-al deep linear neural networks.
\newblock In Alina Beygelzimer and Daniel Hsu, editors, {\em Proceedings of the
  Thirty-Second Conference on Learning Theory}, volume~99 of {\em Proceedings
  of Machine Learning Research}, pages 2691--2713, Phoenix, USA, 6 2019. PMLR.
\newblock URL: \url{http://proceedings.mlr.press/v99/shamir19a.html}.

\bibitem{Tu2011}
Loring~W. Tu.
\newblock {\em An introduction to manifolds}.
\newblock Universitext. Springer, New York, second edition, 2011.
\newblock \href {https://doi.org/10.1007/978-1-4419-7400-6}
  {\path{doi:10.1007/978-1-4419-7400-6}}.

\bibitem{Wojtowytsch2021stochastic}
Stephan Wojtowytsch.
\newblock Stochastic gradient descent with noise of machine learning type.
  {P}art {I}: Discrete time analysis, 2021.
\newblock \href {http://arxiv.org/abs/2105.01650} {\path{arXiv:2105.01650}}.

\bibitem{XuYin2013}
Yangyang Xu and Wotao Yin.
\newblock A block coordinate descent method for regularized multiconvex
  optimization with applications to nonnegative tensor factorization and
  completion.
\newblock {\em SIAM J. Imaging Sci.}, 6(3):1758--1789, 2013.
\newblock \href {https://doi.org/10.1137/120887795}
  {\path{doi:10.1137/120887795}}.

\bibitem{ZhangMartensGrosse2019}
Guodong Zhang, James Martens, and Roger~B Grosse.
\newblock Fast convergence of natural gradient descent for over-parameterized
  neural networks.
\newblock In H.~Wallach, H.~Larochelle, A.~Beygelzimer, F.~d\textquotesingle
  Alch\'{e}-Buc, E.~Fox, and R.~Garnett, editors, {\em Advances in Neural
  Information Processing Systems 32}, pages 8082--8093. Curran Associates,
  Inc., 2019.
\newblock URL:
  \url{http://papers.nips.cc/paper/9020-fast-convergence-of-natural-gradient-descent-for-over-parameterized-neural-networks.pdf}.

\bibitem{ZouCaoZhouGu2019}
Difan Zou, Yuan Cao, Dongruo Zhou, and Quanquan Gu.
\newblock Gradient descent optimizes over-parameterized deep {R}e{LU} networks.
\newblock {\em Machine Learning}, 109:467--492, 2020.
\newblock \href {https://doi.org/10.1007/s10994-019-05839-6}
  {\path{doi:10.1007/s10994-019-05839-6}}.

\end{thebibliography}

\end{document}